%%AMSTEX
\magnification1200
\overfullrule0pt
\input amstex
\documentstyle{amsppt}
\def\Proof{\demo{Proof}}
\def\endproof{\enddemo}
\def\vs{\vskip4pt plus2pt}
\def\ssk{\vskip2pt plus1pt}
\setbox1=\hbox{(iii)\ \enspace}
\def\Mat{\operatorname{Mat}}
\def\rank{\operatorname{rank}}
\def\Det{\operatorname{Det}}
\def\dim{\operatorname{dim}}
\def\codim{\operatorname{codim}}
\def\deg{\operatorname{deg}}
\def\Im{\operatorname{Im}}
\def\card{\operatorname{card}}
\def\top{{\text{\rm top}}}
\def\Hom{\operatorname{Hom}}
\def\cl{\operatorname{cl}}
\def\Cl{\operatorname{Cl}}
\def\Pic{\operatorname{Pic}}
\def\Ker{\operatorname{Ker}}
\def\Fl{\operatorname{Fl}}
\def\id{\operatorname{id}}
\def\pr{\operatorname{pr}}
\def\Nm{\operatorname{Nm}}
\def\Supp{\operatorname{Supp}}
\def\Gr{\operatorname{Gr}}
\let\scs\scriptstyle
\let\j=\i

\def\oms#1{\overline{\mathstrut #1}}
\def\Bb#1{{\Bbb #1}}
\loadeufm
\loadbold
\def\bfr{{\bold r}}
\def\bfD{{\bold D}}

\def\bfG{{\bold G}}
\def\bfX{{\bold X}}
\def\bfZ{{\bold Z}}

\def\hak {^{\vee}}
\def\kr{\! \cdot \!}
\def\M{_{*}}
\def\tilrel {\ _{\widetilde{ }}\ }

\let\ba=\backslash
\let\bi=\bigcirc
\let\lo=\lozenge
\let\tr=\vartriangle
\let\toto=\twoheadrightarrow
\let\hoto=\hookrightarrow

\let\dup=\diagup
\let\do=\diagdown
\let\ge=\geqslant
\let\le=\leqslant

\document
\topmatter
\title SYMMETRIC POLYNOMIALS AND DIVIDED DIFFERENCES 
IN FORMULAS OF INTERSECTION THEORY \endtitle
\rightheadtext{SYMMETRIC POLYNOMIALS AND DIVIDED
DIFFERENCES IN FORMULAS ...}
\author Piotr\ Pragacz 
\thanks During the preparation
of this paper, the author greatly benefited from the hospitality of 
the Max-Planck-Institut f\"ur Mathematik in Bonn and was partially
supported by the Alexander von Humboldt Stiftung and KBN grant No.
2 P301002 05. \endthanks \endauthor

\address Max-Planck-Institut f\"ur Mathematik,
Gottfried-Claren Strasse 26, \break
D-53225 Bonn, Germany \endaddress
\email pragacz\@mpim-bonn.mpg.de \endemail
\subjclass Primary 05E05, 14C15, 14C25, 14M15, 14N10,\break 57R20; 
Secondary 05E15, 14M12, 14H40, 14J60, 32S20, 55R40\endsubjclass
\endtopmatter
\vskip 20pt
\centerline{\eightrm To the memory of Constantin Banica (1942-1991)
and Krzysztof Jaczewski (1955-1994)}
\vskip30pt

\document
\head Introduction \endhead
The goal of this paper is at least two-fold. At first we attempt 
to give a survey of some recent 
(and developed up to the time of the Banach Center workshop
{\it Parameter Spaces}, February '$94$\footnote{Several notes have been
added on the galley proof; they provide information about some related 
results learned by the author after the above date.}) applications 
of the theory of symmetric polynomials and divided differences 
to intersection theory. 
Secondly, taking this opportunity, we complement
the story by either presenting some new proofs
of older results (and this takes place usually in the Appendices
to the present paper) or providing some new results which arose as by-products
of the author's work in this domain during last years.

Being in the past a good part of the classical algebraic knowledge
(related for instance to the theory of algebraic equations
and elimination theory), 
the theory of symmetric functions is rediscovered and developed nowadays 
(see, for example, the monograph [M1] of I.~G. Macdonald or the booklet
[L-S1] of A. Lascoux and M.-P. Sch\"utzenberger). 
Here, we discuss only some geometric applications of symmetric polynomials 
which are related to the present interest of the author. In
particular, the theory of {\it polynomials universally
supported on degeneracy loci} ([P3]) is surveyed in 
Section 1.

Divided differences appeared already in the interpolation
formula of I.~Newton [N, Liber III, p.~582, Lemma V: {\it ``Invenire lineam curvam
generis parabolici, qu\ae\ per data quotcunque puncta tran$\int\!$ibit."}]. 
Their appearance in intersection theory is about twenty years old starting 
with the papers [B-G-G] of I.~N. Bernstein, I.~M. Gelfand and S.~I. Gelfand 
and [D1,2] of M. Demazure.
A recent work [F2] of W. Fulton has illuminated the importance of divided
differences to flag degeneracy loci. This was possible thanks to
the algebraic theory of Schubert polynomials developed recently
by A. Lascoux and M.-P. Sch\"utzenberger ([L-S 2-6]).
\vs
 
The geometrical objects we study are: (ample) vector bundles, 
degeneracy loci of vector bundle homomorphisms, 
flag varieties, Grassmannians including isotropic 
Grassmannians, i.e. the parameter spaces for isotropic subspaces of a given 
vector space endowed with an antisymmetric or symmetric form, Schubert 
varieties and the parameter spaces of complete quadrics.

The algebro-combinatorial tools we use are: Schur polynomials including 
supersymmetric and $Q$-polynomials, binomial determinants and Pfaffians, 
divided differences,
Schubert polynomials of Lascoux and Sch\"utzenberger,
reduced decompositions in the Weyl groups and Young-Ferrers' diagrams.

The content of the article is as follows:
\ssk
{\parindent =20pt
\item{1.} Polynomials universally supported on degeneracy loci,
\ssk
\item{2.} Some explicit formulas for Chern and Segre classes of tensor 
bundles with applications to enumerative geometry,
\ssk
\item{3.} Flag degeneracy loci and divided differences,
\ssk
\item{4.} Gysin maps and divided differences,
\ssk
\item{5.} Fundamental classes, diagonals and
Gysin maps,
\ssk
\item{6.} Intersection rings of spaces $G/P$, divided differences
and formulas for isotropic degeneracy loci --- an introduction to [P-R 2-5],
\ssk
\item{7.} Numerically positive polynomials for ample vector
bundles with applications to Schur polynomials of Schur bundles 
and a vanishing theorem.
\par}

\vs

Apart of surveyed results, the paper contains also some new ones.
Perhaps the most valuable contribution, contained in Section 5, is provided by 
a method of computing
the fundamental class of a subscheme  using the class of the diagonal of the 
ambient scheme. The class of the diagonal can be determined with the help
of Gysin maps (see Section~5). This method has been applied successfully
in [P-R5] and seems to be useful also in other settings. Other results that
appear to be new are contained in Proposition 1.3(ii), Proposition 2.1
and Corollary 7.2. Moreover, the paper is accompanied by a series of appendices
which contain an original material but of more 
technical nature than the main text of the paper. 
Some proofs in the Appendices use an operator approach and the operators
involved are mostly divided differences. This point of 
view leads to more natural proofs of many results than the ones 
known before, and we hope to develop it in [L-L-P-T]. 

The following is the list of appendices:
\ssk
{\eightrm
\ssk
\noindent
A.1. \ Proof of Proposition 1.3(ii).
\ssk
\noindent
A.2. \ Proof of Proposition 2.1.
\ssk
\noindent
A.3. \ Recursive linear relations for $(\!(J)\!)$ and $[J]$.
\ssk
\noindent
A.4. \ A Gysin map proof of the formula from Example 3.5.
\ssk
\noindent
A.5. \ An operator proof of the Jacobi-Trudi identity.
\ssk
\noindent
A.6. \ A Schur complex proof of the Giambelli-Thom-Porteous formula.
\ssk
\noindent
A.$\tau$. \ Corrigenda and addenda to some former author's papers.
\par}

\vs

Several open problems are stated throughout the text. 

We use this opportunity to complete or correct some surveyed results.
Moreover, we give, in Appendix A.$\tau$,
an errata to some former author's papers.

\head Acknowledgement.\endhead
This paper, being a revised and substantially
extended version of the Max-Planck-Institut f\"ur Mathematik Preprint
No. 92-16 {\it Geometric applications of symmetric polynomials; some recent
developments}, is also an expanded version of two talks given by the
author during the workshop {\it Parameter Spaces; enumerative geometry,
algebra and combinatorics}. 
In fact, these talks were given in a linkage with two lectures by 
Alain Lascoux to whom the present paper owes a lot. 
At first, it was Lascoux who introduced me several years ago 
to this branch of mathematics and I wish to express to him my 
sincere gratitude.
Secondly, I~learned some material exposed in the present paper directly
from Lascoux. 
These are: the divided-differences interpretation of the symmetrizing 
operators in Propositions 4.1--4.4 and the content of Appendix A.3.
I thank for his permission to include this material here.

It is a pleasure to thank W. Fulton for his many years' encouragement
given to my work and for informing me about an argument reproduced
in Addenda to [DC-P] in the last Appendix.

The material surveyed in Section 6 stems mainly from a recent series of
papers written in collaboration with J. Ratajski.

I thank also I.~G. Macdonald for pointing me out some errors in the previous
version.

\head 1. Polynomials  universally  supported  on  degeneracy  loci.\endhead
This section summarizes mainly a series of results from  [P1-5],
[P-P1,2] and [P-R1].

Let  $\Mat_{m\times n} (K)$  be the affine space of $m\times n$ matrices 
over a field $K$. 
The subvariety $D_r$ of  $\Mat_{m\times n}(K)$  consisting of all matrices 
of rank $\le r$ is called a {\it determinantal variety (of order $r$)}. 
Algebro-geometric properties of these varieties were widely investigated 
in the seventies and eighties. 
The prototype of the results of this section is, however, an older result ---
a formula of Giambelli [G3] (1903) (see also [G1] and [G4]) 
for the degree of the  
{\it projective determinantal variety} (i.e. the class of  $D_r\setminus \{0\}$ 
in  $\Bb P(\Mat_{m\times n}(\Bb C))$ ). 
In order to perform his computations Giambelli used the machinery of 
symmetric polynomials developed mainly by 
the 18th- and 19th-century elimination theory.

Determinantal varieties are a particular case of {\it degeneracy loci}
$$
D_r(\varphi) = \bigl\{ x\in X\mid \rank\varphi(x) \le r \bigr\}
$$
$r=0,1,2,\dots$  associated with a homomorphism  $\varphi:F \to E$  
of vector bundles on algebraic (or differentiable) variety $X$. 
This concept overlaps many interesting situations like varieties
of special divisors (called also Brill-Noether loci) in Jacobians,
Thom-Boardmann singularities, variations of Hodge structures in
families of Riemann surfaces --- just to mention a few; for more details
and examples consult [Tu].
 
One of the fundamental problems in the study of concrete subscheme $D$ 
of a given (smooth) scheme $X$ is the computation 
of its fundamental class in terms 
of given generators of the cohomology or Chow ring of $X$. 
For instance, Giambelli's formula mentioned above gives the 
fundamental class of the (projective) determinantal variety in 
$H^*\bigl(\Bb P(\Mat_{m\times n}(\Bb C)),\Bb Z \bigr)$\footnote
{More precisely, Giambelli calculated the degree of 
$D_r(\varphi)$ for a general map 
$\varphi: \Cal O(m_1)\oplus \Cal O(m_2)\oplus \dots \to \Cal O(n_1)
\oplus \Cal O(n_2)\oplus \dots$ . His expression, in the 
notation introduced a bit later, is $\sum s_I(E) \cdot s_{\overline I}(F\hak)$ where 
the sum is taken over all partitions $I$ whose diagram is contained in the
rectangle $(n-r)\times (m-r)$ and the diagram of the partition
$\overline I$ complements the one of $I^{\sim}$ in the rectangle
$(m-r)\times (n-r)$; in today's language, explained in the sequel,
this expression equals $s_{(m-n)^{n-r}}(E-F)$.
We refer the reader to the article by D. Laksov [La] 
about Giambelli's work and life. 
This article contains also a complete bibliography of Giambelli overlapping 
his work on degeneracy loci formulas. 
Perhaps it is worth mentioning that several of Giambelli's
formulas have been recently rediscovered using the Gr\"obner bases
technique --- see e.g. [He-T].}. 
In 1957 R. Thom ([T]) proved that for sufficiently general homomorphisms  
$\varphi:F \to E$, there exists a polynomial, 
depending solely on the Chern classes $c_i(E), c_j(F)$ of $E$ and $F$, 
which describes the fundamental class 
of $D_r(\varphi)$\footnote{As Lascoux points out, there is a little step,
by combining this result of Thom and the 
above mentioned computation of Giambelli, 
to arrive at the formula from Theorem 1.0.}.
This polynomial has been found subsequently by Porteous:
$$
\Det\Bigl[ c_{n-r-p+q}(E-F)\Bigr]_{_{\scs 1\le p,q\le m-r}} 
\leqno(*)
$$
where  $c_k(E-F)$  is defined by:
$$
1+c_1(E-F)+c_2(E-F)+\dots = (1+c_1(E)+c_2(E)+\dots)/ (1+c_1(F)+c_2(F) + \dots) .
$$

Different variants and generalizations of (*) were considered later 
in [K-L], [L1], \hbox{[J-L-P],} [H-T1], [P5]  and recently in [F2]  
(compare Section 3). In particular, note that (*) can be rewritten 
using the {\it Segre classes} of $E$ and $F$ as the determinant:
$$
\Det\Bigl[ s_{m-r-p+q}(E-F)\Bigr]_{_{\scs 1\le p,q\le n-r}} 
$$
where  $s_k(E-F)$  is defined by:
$$
1+s_1(E-F)+s_2(E-F)+\dots = (1+s_1(E)+s_2(E)+\dots)/ (1+s_1(F)+s_2(F) + \dots) .
$$ 
and  $s_k(E)$ is here the $k$-th complete symmetric polynomial ([M1])
in the Chern roots of~$E$. (Note that this definition differs by a sign from
that for the Segre class of a bundle, used in [F1].) 

Today's formulation of the Giambelli-Thom-Porteous formula uses
much weaker assumptions than the transversality conditions in [Po]
thanks to the work of Kempf and Laksov [K-L] and Fulton-MacPherson's 
intersection theory [F1]. (We refer the reader to [F1] for
the notions of algebraic geometry used in the present
article.)

\proclaim{Theorem 1.0} If $X$ is a pure-dimensional 
Cohen-Macaulay scheme and the degeneracy
locus $D_r(\varphi)$, endowed with the scheme structure defined by the
ideal generated by $r+1$-minors of $\varphi$, is of pure codimension
$(m-r)(n-r)$ in $X$ or empty, then
$$
[D_r(\varphi)]=\Det\Bigl[ s_{m-r-p+q}(E-F)\Bigr]_{_{\scs 1\le p,q\le n-r}}\cap [X].
$$
\rm (In a modern treatment of intersection theory of [F1], one constructs,
for every vector bundle homomorphism $\varphi$ over a pure-dimensional
scheme $X$, a {\it degeneracy class} $\Bb D_r(\varphi)\in
A_{\dim X-(m-r)(n-r)}\bigl(D_r(\varphi)\bigr)$ whose image in $A_*(X)$
is given by the right-hand side of the formula of the theorem. If $D_r(\varphi)$
is of pure codimension $(m-r)(n-r)$ then $\Bb D_r(\varphi)$ is a positive
cycle whose support is $D_r(\varphi)$; if, moreover $X$ is Cohen-Macaulay
then $D_r(\varphi)$ is also Cohen-Macaulay and $\Bb D_r = [D_r(\varphi)]
$.)\footnote{%
Recall that if $D\subset X$ is a 
(closed) subscheme then $[D]\in A_*(X)$ is the class of the fundamental
cycle associated with $D$, i.e., if \ $D=D_1\cup \ldots \cup D_n$
 \ is a minimal decomposition into irreducible components then
$$
[D] = \sum_{i=1}^n (\operatorname{length} \Cal O_{D,D_{i}})[D_i],
$$
where $\Cal O_{D,D{_i}}$ is the local ring of $D$ along $D_i$.
Recall also that if $f:X\to Y$ is a proper morphism then it induces a morphism
of abelian groups $f_*:A_*(X)\to A_*(Y)$ such that 
$f_*[V]=\deg (f\big|_V)[f(V)]$ if $\dim f(V)=\dim V$ and 0 --- otherwise.
In particular, if $f$ establishes a birational isomorphism of $V$ and $f(V)$
then $f_*[V]=[f(V)]$.
If $X$ and $Y$ are nonsingular then a morphism $f:X\to Y$ induces a ring 
homomorphism $f^*:A^*(Y)\to A^*(X)$. If $X$, $Y$ are possibly singular
and $f$ is flat (or regular embedding) then there exists a group
homomorphism $f^* : A_*(Y)\to A_*(X)$.
}\endproclaim

The second domain of research concerning nonsingular degeneracy loci 
in nonsingular ambient spaces, is the calculation of 
their Chern numbers (see [H], [Na] and [H-T2]). 
Here, one deals with complex varieties and the problem is to find
expressions depending solely on $c_i(E)$, $c_j(F)$ and $c_k(TX)$ for such
numbers. A natural extension of this question is to ask about similar
universal formulas for the topological Euler-Poincar\'e characteristic
of $D_r(\varphi)$, or, even more, for the Chern-Schwartz-MacPherson
classes of these varieties, now without the smoothness assumption
on $X$ and $D_r(\varphi)$.

Finally, the third kind of problems stems from a study of different 
type homology of degeneracy loci (compare [Tu]).
\ssk
 
It turns out that all these questions are closely related with the following
problem whose investigation started with the author's papers [P1,3].

\proclaim{Problem}\rm Which polynomials in the Chern classes of $E$ and $F$ 
are universally supported on the $r$-th degeneracy locus?
\endproclaim

To state this problem precisely, assume that a homology theory $H(-)$
is given which is a covariant functor for proper morphisms
and is
endowed with Chern classes associated with vector bundles on a 
given variety $X$, acting as operators on $H(X)$. Then also the polynomials
in the Chern classes of vector bundles act as operators on $H(X)$.
For example, the Chow homology, the singular homology and the Borel-Moore
homology have these properties (see [P-R1] for more on that).

Let  $\iota_r:D_r(\varphi) \to X$  be the inclusion and let 
$(\iota_r)_*:H\bigl(D_r(\varphi)\bigr) \to H(X)$  be the induced morphism 
of the homology groups. 
Fix integers  $m,n > 0$  and  $r\ge 0$. 
Introduce $m+n$ variables  $c_1,\dots,c_n;\ c'_1,\dots,c'_m$  such that  
$\deg (c_i)=\deg (c'_i)=i$. 
Let  $\Bb Z[c.,c'.]=\Bb Z[c_1,\dots,c_n,c'_1,\dots,c'_m]$  be the polynomial 
algebra. 
Following [P1,3]  we say that $P\in \Bb Z[c.,c'.]$ 
{\it is universally supported on the $r$-th degeneracy locus} if 
$$
P\bigl(c_1(E),\dots,c_n(E),c_1(F),\dots,c_m(F)\bigr)\cap \alpha \in \Im (\iota_r)_*
$$
for any homomorphism $\varphi:F\to E$ of vector bundles on $X$ such that $n=\rank E$, 
$m=\rank F$ and any $\alpha \in H(X)$. 
Denote by  $\Cal P_r$ the set (which is, in fact, an ideal) of 
all polynomials universally supported 
on the $r$-th degeneracy locus. 
Of course, the Giambelli-Thom-Porteous polynomials (*) 
describing $D_i(\varphi)$ for $i\le r$  belong to 
$\Cal P_r$, but they do not generate this ideal if $r\ge 1$. 
An analogous problem can be stated for symmetric (resp. antisymmetric) 
morphisms: $F=E\hak,\ \varphi\hak=\varphi$ (resp. $\varphi\hak= -\varphi)$. 
In this case the corresponding ideal $\Cal P^s_r$  
(resp. $\Cal P^{as}_r$  $r$-even) is contained in 
$\Bb Z[c_1,\dots,c_n] = \Bb Z[c.]$.

It follows from the ``main theorem on symmetric polynomials" that for a 
sequence of variables $A=(a_1,\dots,a_n)$, where $\deg (a_i)=1$, the assignment
$$
c_i\mapsto (i\hbox{\rm -th elementary symmetric polynomial in\ } A)
$$
defines an isomorphism of $\Bb Z[c.]$  and  $S\Cal P(A)$ 
--- the ring of symmetric polynomials in $A$. 
Similarly, by considering an analogous assignment for the 
$c'_j$'s  and a second sequence of variables $B=(b_1,\dots,b_m)$, we get an 
isomorphism of $\Bb Z[c.,c'.]$ with $S\Cal P(A|B)=S\Cal P(A)\otimes 
S\Cal P(B)$ --- the ring of symmetric polynomials in $A$ and $B$ separately.  

A precise description of the ideals 
$\Cal P_r$, $\Cal P^s_r$  and $\Cal P^{as}_r$
requires two families of symmetric polynomials.
\hfill\break
(i)  Let $I=(i_1,\dots,i_k)$  be a sequence of integers. 
We define 
$$
s_I(A-B) = \Det \Bigl[s_{i_p-p+q}(A-B)\Bigr]_{_{\scs i\le p,q\le k}} ,
$$
where  $s_i(A-B)$  is a homogeneous polynomial of degree $i$ such that
$$
\sum\limits^{\infty}_{i=-\infty} s_i(A-B)=
\prod\limits^n_{i=1}(1-a_i)^{-1}\ \prod\limits^m_{j=1}(1-b_j) .
$$
Observe that the corresponding polynomials $s_I(c./c'.)$ 
in the variables $c.$ and $c'.$ are determined by
$$
s_i(c./c'.) = s_i - s_{i-1}c_1' + \dots + (-1)^{i-1}s_1c_{i-1}' + (-1)^i c_i', $$
where 
$$s_i=s_{i-1}c_1-s_{i-2}c_2+ \ldots +(-1)^{i-2}s_1c_{i-1}+(-1)^{i-1}c_i$$ 
for $i>0$ and $s_i=0$ for $i<0$, $s_0=1$. 
\ssk

Moreover, we put $s_I(A)=s_I(A-B)$ for 
$B=(0,\ldots,0)$ and similarly $s_I(c.)=s_I(c./c'.)$ 
for $c_j'=0$, $j=1,\ldots,m$.
\ssk

\noindent
(ii) Let $Q_i(A)$ be a symmetric polynomial defined by the expansion
$$
\sum\limits^{\infty}_{i=-\infty} Q_i(A) t^i =
\prod\limits^n_{i=1}(1+ta_i)(1-ta_i)^{-1}. 
$$  
Given nonnegative integers $i,j$, we set
$$
Q_{i,j}(A)=Q_i(A)\ Q_j(A) + 2\sum\limits^j_{p=1}(-1)^p Q_{i+p}(A) Q_{j-p}(A).
$$
Finally, if $I=(i_1,\dots,i_k)$  is a sequence of positive integers
then for odd $k$ we put
$$
Q_I(A) = \sum_{p=1}^k (-1)^{k-1}Q_{i_p}(A) \ Q_{i_1,\ldots,i_{p-1},i_{p+1},
\ldots,i_k}(A),
$$
and for even $k$,
$$
Q_I(A) = \sum_{p=2}^k (-1)^kQ_{i_1,i_p}(A) \ Q_{i_2,\ldots,i_{p-1},i_{p+1},
\ldots,i_k}(A).
$$

Observe that the corresponding polynomials $Q_I(c.)$ in the variables
$c.$ are determined by
$$
Q_i(c.) = s_i + s_{i-1}c_1 + \dots + s_1c_{i-1} + c_i
$$ 
for $i>0$ and $Q_0(c.)=1$, $Q_i(c.)=0$ for $i<0$.

The polynomials $s_I(A)$ and $s_I(c.)$ are called {\it Schur polynomials}
or  $S$-{\it polynomials}.
The polynomials $s_I(A-B)$ and $s_I(c./c'.)$ are often called  
{\it supersymmetric Schur polynomials} --- for an account 
to their properties we refer to [P4] and [P-T]. 
The polynomials $Q_I(A)$ are called  {\it Schur Q-polynomials}
--- for an account to their properties we refer to [H-H] and [P4]. For
another expression of $Q_I(A)$ in the form of a quadratic polynomial 
in the $s_J(A)$'s, see [La-Le-T1].

\vs

Now, let $E$ and $F$ be two vector bundles on $X$. 
Then  $s_I(E-F)$ is obtained from  
$s_I(c./c'.)$ via the specialization $c_i:= c_i(E)$, $i=1,\dots,n$; \ 
$c'_j:= c_j(F)$, $j=1,\dots,m$;
and $s_I(E)$ -- from $s_I(c.)$ by the substitution
$c_i:=c_i(E)$.
Similarly we define $Q_I(E)$ as the specialization of $Q_I(c.)$ with
$c_i:=c_i(E)$.

Recall that by a  {\it partition}  (of $n$) we understand a sequence of integers  
$I=(i_1,\dots,i_k)$, where  $i_1\ge i_2\ge \dots\ge i_k\ge 0$ and 
$\sum i_p=n$. 
A partition with strictly decreasing parts is called {\it strict}. 
For partitions $I,J$  we write  $I\supset J$  if $i_1\ge j_1$, $i_2\ge j_2$, 
\dots ; the partition $(i,\dots,i)$ ($r$-times) is denoted by $(i)^r$; 
finally the partition $(k,k-1,\dots,2,1)$ is denoted by $\rho_k$. 

Note that for every strict partition $I=(i_1>\dots>i_k>0)$, one has 
$Q_I(c.)= 2^k P_I(c.)$ for some polynomial $P_I(c.)$ 
with integer coefficients. 
These polynomials are called {\it Schur P-polynomials}. At first,
the ideals $\Cal P_r$, $\Cal P_r^s$ and $\Cal P_r^{as}$ were described
for the Chow homology. Let us give first a coarse description: 
  
\def\j#1;{\item{{\rm(#1)\ }}}
\proclaim{Theorem 1.1 {\rm[P1,3]}} Assume that $H(-)=A_*(-)$ is the
Chow homology theory. Then\par
{\parindent\wd1
\j i;The ideal $\Cal P_r$ 
is generated by  $s_I(c./c'.)$, where $I$ runs over all partitions\break
$I\supset (m-r)^{n-r}$.
\j ii;The ideal $\Cal P^s_r$ is generated by $Q_I(c.)$, where $I$ runs over 
all partitions  $I\supset \rho_{n-r}$.
\j iii;The ideal $\Cal P^{as}_r$ ($r$-even) is generated by 
$P_I(c.)$, where $I$ runs over all partitions $I\supset\rho_{n-r}$.\par}
\endproclaim

Observe that the ``positive" generator of the ideal $\Cal P_r$ agrees
with the Segre class version of the Giambelli-Thom-Porteous polynomial
from Theorem 1.0. The analogous generators of the ideals $\Cal P_r^s$
and $\Cal P_r^{as}$ are of different (Pfaffian) form than the determinantal
expressions given in [J-L-P], [H-T1] and [P5].
  
To prove that the quoted polynomials belong to $\Cal P_r$, $\Cal P^s_r$ 
and $\Cal P^{as}_r$, the key tools are certain factorization formulas 
and formulas for the Gysin map for Grassmannian bundles.

In the sequel, having two partitions $I$ and $J$ with $l(I)\le k$ and
$l({J^{\sim}})\le i$\footnote{%
For a given partition $I$, $l(I)=\card \{p:i_p>0$\} denotes its length and 
$I^{\sim}$ denotes the partition conjugate to $I$, i.e., 
$(h_1,h_2,\dots)$ where $h_p=\card \{q:i_q\ge p\}$.},
by \ $(i)^k+I,J$ \ we denote the partition $(i+i_1,\ldots ,i+i_k,j_1,j_2,
\ldots)$.    

%\eject
\proclaim {Proposition 1.2 \rm (Factorization Formula)} 
Let $I$, $J$  be two partitions such that $l(I)\le n$ and 
$l(J^{\sim})\le m$.  
Then\footnote{For a given partition $I$, we write $|I|:=\sum i_p$  
   --- the sum of parts of $I$, i.e. the number partitioned by $I$.}
$$
\leqalignno{
s_{(m)^n+I,J}(A-B)=&s_I(A)\ s_{(m)^n}(A-B)\ s_J(-B) &\text{\rm(i)}\cr
 =& (-1)^{|J|}\ s_I(c.)\ s_{(m)^n}(c./c.')\ 
s_{J^{\sim}}(c.')\cr
Q_{\rho_{n-1}+I}(A)=&Q_{\rho_{n-1}}(A)\ s_I(A).&\text{\rm(ii)}\cr}
$$
\endproclaim
\noindent\rm
(Formula (ii) is due to Stanley; 
for the history of (i) we refer to [L3].
Both the formulas are just special instances of much more general identities
for which we refer the interested reader to [P4, 1.3], [P-T] 
and [La-Le-T1].)

\proclaim{Proposition 1.3}  
Let  $\pi=\pi_E:\Cal G=G^q(E)\to X$ be the Grassmannian 
bundle parametrizing  $q$-quotients of E. Write $r=n-q$.
Let  $$0\leftarrow Q=Q_E\leftarrow E_{\Cal G}\leftarrow R=R_E\leftarrow 0$$ 
be the tautological exact sequence of vector bundles on $\Cal G$. 
Let $\alpha\in A_*(X)$.
{\parindent=\wd1
\item{\rm(i)} {\rm [J-L-P], [P3].} \ For every vector bundle $F$ on $X$ 
and any sequences of integers $I=(i_1,\ldots,i_q)$,
$J=(j_1,\ldots,j_r),$ 
$$
(\pi_E)_*\ \Bigl[\ s_I(Q-F_{\Cal G})\ s_J(R-F_{\Cal G})
\cap \pi_E^*\alpha \Bigr] 
= s_{i_1-r,\ldots,i_q-r,j_1,\ldots,j_r}(E-F)\cap \alpha,
$$
\item{\rm(ii)} Let $I=(i_1,\dots,i_k),\ J=(j_1,\dots,j_h)$ be two 
sequences of positive integers, $k\le q$, \ $h\le r$. Then 
$$
(\pi_E)_*\Bigl[c_\top(Q\otimes R)\ P_IQ\ P_JR\cap \pi_E^*\alpha \Bigr]= 
d P_{I,J}(E)\cap \alpha, $$ \par}
where $d$ is zero if $(q-k)(r-h)$ is odd, and $(-1)^{(q-k)r}
{[(n-k-h)/2]\choose[(q-k)/2]}$ 
--- otherwise\footnote{$I,J$ denotes here and in Appendix A.1
the juxtaposition of $I$ and $J$.}.  
\endproclaim

For a proof of (ii), see Appendix A.1. 

\vs

Propositions 1.2 and 1.3 allow one to prove the following algebraic 
result providing finite sets of generators of the ideals in question.  

\proclaim{Proposition 1.4 \rm [P3]}
\noindent
{\parindent=\wd1
\j i;  $\bigl(s_I(c./c'.),\ I\supset(m-r)^{n-r}\bigr)= 
\bigl(s_{(m-r)^{n-r}+I}(c./c'.),\ I\subset(r)^{n-r}\bigr)$,
\j ii; $ \bigl(Q_I(c.),\ I\supset \rho_{n-r}\bigr)=
\bigl(Q_{\rho_{n-r}+I}(c.),\ I\subset (r)^{n-r}\bigr)$,
\j iii; $ \bigl(P_I(c.),\ I\supset \rho_{n-r-1}\bigr)=
\bigl(P_{\rho_{n-r-1}+I}(c.),\ I\subset (r)^{n-r}\bigr)$, \ $r$ --- even.
\par}\endproclaim

\noindent
--- thus these ideals are generated by  $n \choose r$  elements.
\vs

Note that it is still an open problem to show that these sets form
minimal sets of generators of the corresponding ideals 
(in case (i), we assume that $m\ge n$).

For an explicit $\Bb Z$-basis of the ideal in (i), see [P3, 
Proposition 6.2]. It would be valuable to have a similar result
for the ideals in (ii) and (iii).
Moreover, the ideal in (i) is prime ([P2,4]), and is a set-theoretical
complete intersection (is equal to the radical of an ideal generated by 
a regular sequence of length $r+1$ (loc.cit.)).

As shown in loc.cit., the ideal in (i) gives
a generalization of the resultant of two polynomials in one variable. 
Let $$A(x)=x^n+\sum^n_{i=1}c_ix^{n-i} , \ \  
B(x)=x^m+ \sum^m_{j=1}c'_j x^{m-j}$$
be two polynomials in one variable with generic coefficients. 
It follows from the classical algebra, that there exists a polynomial 
in $\{c_i\}$, $\{c'_i\}$  called the {\it resultant}, whose vanishing 
(after a specialization of $\{c_i\}$, $\{c'_i\}$ to an algebraically closed 
field) implies that the corresponding polynomials have a common root 
(see, e.g., [L3] for an approach to the resultant via the symmetric
polynomials).  

Now, let $\Cal T_r$ be the ideal of all $P\in \Bb Z[c.,c'.]$, 
which vanish if, after a specialization of 
$\{c_i\}$, $\{c'_i\}$ to a field,  
$A(x)$ and $B(x)$  have $r+1$ common roots. 
Surprisingly (or not) we have

\proclaim{Theorem 1.5 \rm[P2,4]}
$\Cal T_r=\left(s_{(m-r)^{n-r}+I}(c./c'.),\ I\subset (r)^{n-r}\right).$
\endproclaim

In other words $\Cal T_r=\Cal P_r$ in the above notation. 
It would be interesting to have an {\it intrinsic} proof of this equality.
It is shown in [L-P] that an analogous ideal defined
in the ring $\Bb Z[AB]$ of all polynomials in $A$ and $B$
is just generated by $\Cal T_r\subset S\Cal P(A|B)
\subset \Bb Z[AB]$. 
A similar interpretation is given in [P4] (and correspondingly in [L-P]) 
for the ideals $\Cal P^s_r$ and 
$\Cal P_r^{as}$ generated by $Q$- and $P$-polynomials respectively.

\vs

Let us come back to Theorem 1.1.
The proof that the ideal $\Cal P_r$ is actually generated by the above 
polynomials is based on the investigation of the tautological 
determinantal variety  $\underline{D}_r \subset \underline{\Hom}(F,E)$ 
(the fibre of $\underline{D}_r$ over a point $x\in X$ is equal to 
$\{f\in \Hom\bigl(F(x),E(x)\bigr)\mid \rank (f)\le r \}$). 
The bundles $E$ and $F$ occurring in this construction are some 
``universal enough" vector bundles over the product $GG$ of two Grassmannians 
(see [P3]). 
In fact, in [P3], two proofs of this assertion are given. 
One of them [P3, pp.~441--445] is by induction on 
$r$ with the help of an exact sequence of Chow groups
$$
A_*(\underline{D}_{r-1})\to A_*(\underline {D}_r)\to 
A_*(\underline{D}_r\setminus \underline{D}_{r-1}) \to 0\,
$$
and a detailed analysis of $A_*\bigl(\underline{D}_r\setminus 
\underline{D}_{r-1}\bigr)$.
The second one [P3, pp.~428--432] uses a certain 
desingularization of $\underline{D}_r$ and has
been ameliorated in [P-R1] to give the assertion also for the Borel-Moore
homology and the singular homology.

\proclaim{Theorem 1.6 {\rm[P-R1]}}The statement of Theorem 1.1 is true
also for $H(-)$ being the Borel-Moore homology (both, the classical one
and that defined by Laumon in characteristic $p$) as well as for the
singular homology (with integer coefficients). 
\endproclaim

Since the same applies to Proposition 1.3, when appropriately formulated,
the proof that the quoted polynomials belong to $\Cal P_r$ is the same.
 
On the other hand, the proof that the ideal $\Cal P_r$ is generated by the 
above polynomials uses the following compactification of $\underline{D}_r$. 
Let us embed the above  $\underline{\Hom }(F,E)$  into a Grassmannian bundle 
{\bf X}$=G_m(F\oplus E)$ by assigning fibrewise to 
$f\in \Hom \bigl(F(x),E(x)\bigr)$ its 
(graph of $f$)$\in G_m\bigl(F(x)\oplus E(x)\bigr)$, $x$ belonging to 
the base space $GG$. 
On {\bf X} there exists a natural tautological extension of the universal 
homomorphism on $\underline{\Hom }(F,E)$ and its degeneracy loci 
serve to prove the assertion. 

An important advantage of the above compactification as well
as a certain natural desingularization $\bfZ$ of it is the vanishing of their 
odd homology groups --- this is not the case of $\underline{D}_r$  and its
analogous desingularization (see [P-R1]). Here $\bfZ$ is the subscheme
of zeros of the homomorphism $F_{\bfG}\to E_{\bfG}\to Q$ on ${\bfG}=G_r(E)$
where the first map is the pullback to $\bfG$ of $\varphi$.
Let
$\eta:\bfZ\to \bfD_r$ be the restriction of $\pi:\bfG\to X$ to $\bfZ$,
and let $j$ be the
closed immersion of $\bfZ$ into $\bfG$.
Then by using the rank-stratification $\{\bfD_k\setminus \bfD_{k-1}\}$ 
of $\bfD_r$,
the induced stratification $\bfZ^k\setminus\bfZ^{k-1}$ of {\bf Z}
$(\bfZ^k=\eta^{-1}\bfD_k)$, and proving (for the Borel-Moore homology)
that $\cl_{\bfD_k}$ and $\cl_{\bfZ^k}$
are isomorphisms, one shows that the induced push-forward map 
$\eta_*:H(\bfZ)\to H(\bfD_r)$ is surjective. Also, by analyzing the geometry
of {\bf Z}, one shows that $j^*$ is surjective. This implies, by the projection
formula, that $\Im j_*$ is a principal ideal in $H(\bfG)$ generated by the 
fundamental class [{\bf Z}]. It follows then, from the commutative diagram
$$
\CD
H(\bfZ) @>j_*>> H(\bfG)      \\
@V\eta_*VV         @V\pi_*VV \\
H(\bfD_r)     @>{\iota}_*>> H(\bfX),
\endCD
$$
that $\Im {\iota}_*=\pi_*\bigl([\bfZ]H(\bfG)\bigr)$.
This identity together with some algebra of symmetric polynomials
(which allows one to express explicitly $[\bfZ]H(\bfG)$) yields the desired
assertion about $\Im {\iota}_*$.
In this way we obtain a proof which is valid both for Chow homology and other
homology theories simultaneously.

In a similar way, though overcoming some additional difficulties, one can 
prove the analogous theorem in the symmetric and antisymmetric 
cases.
\vs

Theorem 1.1 allows us to calculate the Chow groups of 
some degeneracy loci.
A prototype of these results is the following result from [B]. 
Let $R$ be a normal noetherian ring, $M$ --- a $m\times n$ matrix of 
indeterminates, $\Cal I$ --- the ideal generated by all $(r+1)$-minors of $M$. 
Then, the divisor class groups satisfy:  
$\Cl(R[M]/\Cal I) \cong \Cl(R)\oplus\Bb Z$. 

The geometric analogue of $\Cl$ is $A^1$ (the Chow group of codimension 1 
cycles modulo rational equivalence).
Keeping the above notation for the tautological degeneracy loci in
$\underline{\Hom }$-bundle one has ($A^i(-)$ denotes below the Chow group
of codimension $i$ algebraic cycles modulo rational equivalence):

\proclaim{Theorem 1.7 {\rm[P3]}}If $m\ge n$ then the Chow group
of $\underline{D}_r$  is canonically isomorphic to the Chow group of
$G_r(E)$. Therefore, for every $i$, $A^i(\underline{D}_r)=\bigoplus A^{i-|I|}(X)$,
the sum over all partitions $I\subset (r)^{n-r}$, $|I|\le i$.
\endproclaim

Let $\Mat_{m\times n}(K)$ denote the affine space of $m\times n$ matrices over 
a field $K$ (assume $m\ge n$ without loss of generality)
and $D_r\subset \Mat_{m\times n}(K)$ be the subscheme 
defined by the ideal generated by all minors of order $r+1$.
The theorem implies, in particular, that  
for every  $K$-scheme $X$,  $A^i(X\times D_r)  \cong\bigoplus A^{i-|I|}(X)$, 
the sum as above. 
For $i=1$  this is a geometric analogue of the result from [B].
Note that the Chow group of $D_r$ is isomorphic to the Chow group 
of the Grassmannian $G_r(K^n)$. This could create an impression that
homologically $D_r$ behave like spaces which admit a cellular
decomposition. This is, however, {\bf not} the case --- see [P-R1] where
it is shown that complex determinantal varieties have nontrivial
Borel-Moore homology groups of {\it {odd}} degree. 

It would be interesting to find analogues of Theorem 1.7 for the tautological
degeneracy loci of homomorphisms with symmetries.
\vs

Finally we pass to perhaps the most spectacular application of the theory
of polynomials universally supported on degeneracy loci. This is
a formula for the Chern-Schwartz-MacPherson classes of degeneracy loci
associated with an $r$-general vector bundle homomorphism
$\varphi: F \to E$ over a (possibly singular) complex analytic
variety $X$. The Chern-Schwartz-MacPherson class $c_*(X)$ 
of a variety $X$ has its value in the Borel-Moore homology of $X$ 
and satisfies similar properties as the Chern class $c(TX)$ of
the tangent bundle of a complex manifold $X$. In particular,
for a possibly singular compact analytic variety $X$, we have the
following expression for the topological Euler-Poincar\'e 
characteristic: 
$$
\chi(X) = \int_X \, c\M (X)  %= \sum n_i\, \int_{V_i} \, c_M(V_i) 
$$ 
(see, for instance [F1, Chap.~19]).
Let us now fix a Whitney stratification $\Cal X$ of $X$.  Let $E$ be 
a holomorphic vector bundle on $X$ and $Z$ --- the variety of zeros 
of a holomorphic section $s$ of $E$. Assume that $s$ intersects,
on each stratum of $\Cal X$, the zero section of $E$ transversely. 
Let $\iota\colon Z\to X$ be the inclusion.  

\proclaim{Lemma 1.8 {\rm[P-P2]}}
$$
\iota_* (c\M (Z)) = c(E)^{-1}\kr c_\top (E) \cap  c\M(X) .
$$
In particular, for a compact analytic variety $X$,
$$
 \chi (Z) = \int_X c(E)^{-1}\kr c_\top (E) \cap  c\M(X) .
$$
\endproclaim

This is the simplest instance of the formula in question. 
To state the result in the most general form we need a notion
of $r$-generality of a vector bundle homomorphism.
We say that 
$\varphi$ is \  {\it $r$-general} \ if the section \ $s_{\varphi}\colon X \to 
\underline{\Hom }(F,E)$ \ induced by $\varphi$ intersects,
on each stratum of the Whitney stratification $\Cal X$, the subset
 \ $\underline{D}_k \setminus \underline{D}_{k-1}$  transversely 
for every $\, k = 0,1,\ldots,r$. 
For a pure-dimensional nonsingular $X$, this condition
can be expressed in a more transparent way: 
a morphism $\varphi$ is
$r$-general iff for every $k=0,1,\dots,r$, the subset $D_k(\varphi)\setminus 
D_{k-1}(\varphi)$ is nonsingular of pure dimension 
$\dim X-(m-k)(n-k)$ (here, $D_{-1}(\varphi)=\emptyset$).
\ssk

Let $m\wedge n$ denote the minimum of $m$ and $n$.
We now define the following element in $H_*(X)$. We set
$$
\Psi(k):= P_k(E,F)\cap c_*(X),
$$ 
where
$$
P_k(E,F):=\sum 
(-1)^{|I|+|J|} \,  D_{I,J}^{m-k,n-k} \, s_{(m-k)^{n-k}+I,J^{\sim}} (E-F).
$$ 
Here, the sum is over all
partitions $I,J$ such that $l(I)\le {m \wedge n}-k$, $l(J)\le {m \wedge n}-k$, 
and the numbers $ D_{I,J}^{m-k,n-k}$ are some binomial determinants which
will be defined in Theorem 2.4(i).

The following formula gives an explicit expression for the image of 
the Chern-Schwartz-MacPherson class of $D_r (\varphi)$ in the homology of $X$. 
Recall that $\iota\colon D_r(\varphi)\to X$ denotes the inclusion.

\proclaim{Theorem 1.9 \rm[P-P1,2]}
If $\varphi$ is $r$-general then one has in $H_*(X)$ 
$$
\iota_* (c_*(D_r(\varphi ))) =  \sum_{k=0}^r \,\, (-1)^k 
\binom {{m \wedge n}-r+k-1}k \Psi(r-k) . 
$$
In particular, if $X$ is a compact analytic variety, then
$$
\chi (D_r(\varphi )) = \int_X \, \sum_{k=0}^r \,\, (-1)^k \binom 
{m \wedge n-r+k-1}k \Psi(r-k) .
$$
\endproclaim 

Under the assumption $D_{r-1}(\varphi) = \emptyset$, the 
above formula reads $\chi (D_r(\varphi ))= \int_X \Psi(r)$. 
This result was established earlier in [P3] as a
particular case of an algorithm for computation 
the Chern numbers of nonsingular degeneracy loci.

The essence of the argument is to pass first to the above described
desingularization
of $D_r(\varphi)$ and calculate explicitly the image 
of the homology dual to its Chern class in the homology of $X$.
To this end, by using some algebra (of symmetric
polynomials and Gysin maps), we show that
this image has the form $P\cap c_*(X)$ where $P=P(\{c_i\},\{c'_j\})$ 
is a polynomial
universally supported on the $r$-th degeneracy locus
and not universally supported on the $(r-1)$-th one,
specialized 
by setting
$c_i=c_i(E)$, $c'_j=c_j(F)$. Thus ``morally", without changing 
the result of the computation,
we can assume that $D_{r-1}(\varphi)=\emptyset$. But then the desingularization
equals $D_r(\varphi)$ and the wanted class is known
by the result of [P3] quoted above.

Secondly, stratifying $D_r(\varphi)$ by the subsets where the rank of
$\varphi$ is constant, the desingularization turns out to be a Grassmannian
bundle over each stratum. This leads to an equation with the known
$H_*(X)$-image of the Chern class of the desingularization
on the one side and a linear combination of the unknown $H_*(X)$-images
of the Chern-Schwartz-MacPherson classes of $D_k(\varphi)$ \ ($k\le r$) --- 
on the other. By varying $r$, this leads to a system of linear equations
in the unknown $H_*(X)$-images
of the Chern-Schwartz-MacPherson classes of $D_r(\varphi)$  \ (and  
with known coefficients). Solving this system of equations with the help
of some algebra of binomial numbers, one gets the formula looked at.

\vs

As a by-product of our considerations, we also get a formula for 
the {\it Intersection
Homology-Euler characteristic} of $D_r(\varphi)$ associated with 
an \ $r$-general morphism $\varphi$:

\proclaim{Theorem 1.10 \rm[P-P2]}
If $X$ is nonsingular compact analytic variety 
and $\varphi$ is $r$-general, then
$$
\chi_{IH} (D_r(\varphi )) = \int_X \Psi(r) .
$$
\endproclaim

As an example of application of Theorem 1.9, we provide an expression for
the topological Euler-Poincar\'e characteristic of the Brill-Noether loci 
$W_d^r(C)=\{ L \in \Pic ^d(C) \mid h^0(C,L)>r \}$ parametrizing
all complete linear series of degree $d$ and dimension $r$ on a general
curve $C$ of genus $g$.
Let $ \rho:= \rho (r):= \rho (g,d,r):= g - (r+1)(g-d+r)$ be the
Brill-Noether number. 

For $\rho (r)\ge 0$, let 
$$\Phi(g,d,r)=(-1)^{\rho(r)} g! 
\sum D_{I,J}^{r+1,g-d+r}/h(I_{g,d,r}+I,J^{\sim}),
$$
where $I_{g,d,r}$ is the partition 
$(r+1)^{g-d+r}$, the sum is over 
partitions $I,J$ with
length $\le (r+1)\wedge (g-d+r)$ and such that $|I|+|J|=\rho(r)$. Moreover, for a partition $I$,
$h(I)$ denotes the product of all hook lengths associated with the
boxes in the Ferrers' diagram of $I$ (see [M1, Chap.~I]).  
We set $\Phi(g,d,r)=0$ if $\rho (r)<0$. 

\proclaim{Theorem 1.11 \rm[P-P2]}Assume that a curve $C$ of genus $g$ 
is general. 
Let $d$, $r$ be integers as above and such that 
$\rho(r)\ge 0$.
Then one has
$$
\chi (W_d^r(C)) = \sum_{k\ge r} (-1)^{k-r} {k \choose k-r}\Phi(g,d,k). 
$$
\endproclaim

From this formula, one deduces the following corollary. If we fix 
$g,d,r$ such that $\rho(r) \ge 0$ and the nonnegative numbers 
$\rho(r)$, $\rho(r+1), \dots$
change successively the parity, then $\chi(W_d^r(C))>0$ (resp. 
$\chi(W_d^r(C))<0 )$
iff $\rho(r)$ is even (resp. $\rho(r)$ is odd). Observe that 
the above numbers
change successively the parity if $r+1$ and $r+g-d$ are
of the same parity. The latter condition holds iff $g\not\equiv d\pmod 2$.
Thus we get, in the situation of the theorem, the following result.

\proclaim{Corollary 1.12 \rm[P-P2]}Assume 
$g\not\equiv d\pmod 2$. Then one has $\chi(W_d^r(C))<0$ $($resp. 
$\chi(W_d^r(C))>0)$
iff $g\equiv r\pmod 2$ \ $($resp. $g\not\equiv r\pmod 2)$.
\endproclaim

For example, if $\rho=0$, we have 
$$
\Phi(g,d,r)=\card (W^r_d(C))=g!/h((r+1)^{g-d+r}),
$$
which is the classical Castelnuovo formula, expressed here using the hook
number of $I_{g,d,r}$.

\head 2. Some  explicit  formulas  for  Chern  and  Segre  classes  of  
tensor  bundles  with  applications  to  enumerative  geometry.\endhead
In this paper, by $S^IE$ we will denote the {\it Schur bundle} associated
with a bundle $E$ and partition $I$ (whenever we speak about Schur bundles, we 
assume, for simplicity, that the ground field is of characteristic zero).
Recall that if $|I|=n$ and if $S_n$
stands for the symmetric group with $n!$ elements, then $S^IE=
\Hom _{S_n}(\Sigma^I, E^{\otimes n})$ where $\Sigma^I$ is the corresponding
irreducible representation of the group $S_n$ and this group
acts on $E^{\otimes n}$ via the permutations of the factors. 
Thus in particular $S^{(n)}E=S_n(E)$, the $n$-th symmetric power; and
$S^{(1)^n}E=\Lambda^n(E)$, the \hbox{$n$-th} exterior power.
In other words, $S^IE$ is the tensor bundle of $E$ associated with 
the irreducible representation of $GL_n$ defined by~$I$.

The problem of determining the Schur polynomials decomposition of
$s_I(S^JE)$ is very far of being solved. The present section and Section 7
provide some partial information related to this question.
   
Throughout this paper, for a vector bundle $E$, 
we write $c_\top (E)$ instead of $c_{\rank E}(E)$. We show first that
the Schur polynomials decomposition of $c_\top (S^JE)$
determines the one of $c(S^JE)$. 

\proclaim{Proposition 2.1}If $c_\top (S^JE)=\sum_K m_K \ s_K(E)$,
the sum over partitions $K$, \ then
$$
c(S^JE)= |J|^{-\rank (S^JE)} \sum_{K} \sum_{L\subset K}|J|^{|L|} \ m_K 
\ d_{KL} \ s_L(E).
$$
where the sum is over partitions $K=(k_1,\ldots,k_n)$, $L=(l_1,\ldots,l_n)$,
$n=\rank E$, and 
$$
d_{KL} = \Det\left[{k_p+n-p \choose l_q+n-q}\right]_{1\le p,q\le n}.
$$
\endproclaim

For a proof see Appendix A.2.

\vs
(In particular, note that if $c_\top (S^JE)$ and the Segre classes
$s_i(S^JE)$, $i\le p$, are known, then the remaining Segre classes 
$s_i(S^JE)$, $i>p$, are also determined.)

\vs
Recall that the Schur polynomials 
decompositions of $c(S^2E)$ and $c(\Lambda^2(E))$ 
are known.

\proclaim{Proposition 2.2 \rm [L2]}If $\rank E=n$ then 
$$c_\top (S^2E)=2^n s_{\rho_n}(E) \quad \hbox {and} \quad 
c_\top (\Lambda^2(E))=s_{\rho_{n-1}}(E).$$
\endproclaim

\noindent
(Note that [L2] also contains a formula for the 
decomposition of $c(E\otimes F)$ into Schur polynomials.) 

\proclaim{Example 2.3}\rm In the following, $s_I=s_I(E)$.\hfill\break
If $\rank E = 4$ then 
$$
c_\top (\Lambda^3(E)) = s_{3,1} + s_{2,2} +s_{2,1,1} +s_{1,1,1,1}.
$$
If $\rank E = 5$ then 
$$\eqalign{
c_\top (\Lambda^3(E)) &= 9s_{3,3,2,1,1} + 3s_{3,3,2,2} + 2s_{3,3,3,1} + 
9s_{4,2,2,1,1} + 3s_{4,2,2,2}+ 6s_{4,3,1,1,1} \cr
& + 9s_{4,3,2,1} + 3s_{4,3,3} + 3s_{4,4,1,1} + 3s_{4,4,2} 
+ 4s_{5,2,1,1,1} + 4s_{5,2,2,1} \cr
&+ 4s_{5,3,1,1} + 4s_{5,3,2} + 2s_{5,4,1} + s_{6,2,1,1} 
+ s_{6,2,2} + s_{6,3,1} \cr
&+ 6s_{3,2,2,2,1} + s_{3,3,3,1} + s_{6,1,1,1,1}.\cr}
$$
If $\rank E=2$ then
$$
c_{n+1}(S^nE)= \prod\limits_{j=0}^{(n-1)/2}[j(n-j)s_2 + (n^2-3j(n-j))s_{1,1}]
$$
for $n$ odd, and
$$
c_{n+1}(S^nE)= (n/2)s_1\cdot \prod\limits_{j=0}^{n/2-1}[j(n-j)s_2 
+ (n^2-3j(n-j))s_{1,1}]
$$
for $n$ even.
\endproclaim
   
The rest of this section summarizes some results from [La-La-T] and [P3].

Let $E,F$ be vector bundles of ranks $n$ and $m$ respectively. 
Assume $m\ge n$. 
We state

\proclaim{Theorem 2.4}\parindent=\wd1
\j i;{\rm[L-S1], [La-La-T]} 
The total Segre class of the tensor product $E\otimes F$ 
is  given by
$$
s(E\otimes F) = \sum D^{n,m}_{I,J}\ s_I(E)\ s_J(F),
$$
where the sum is over partitions $I,J$  of length $\le n$ and 
$$
D^{m,n}_{I,J}=\Det\left[{i_p+j_q+m+n-p-q\choose i_p+n-p}\right]_{1\le p,q\le n}.
$$
\ssk
\j ii;{\rm [La-La-T] \& [P3]} The total Segre class of the second symmetric 
power $S^2E$ is given by
$$
s(S^2E)=\sum (\!(I+\rho_{n-1})\!)\ s_I(E),
$$
where the sum is over all partitions $I$ and the definition of   $(\!(J)\!)$,
for $J=(j_1>\ldots>j_n\ge 0)$, is
as follows. 
If $n$ is even, define $(\!(J)\!)$ to be the Pfaffian of the $n\times n$ 
antisymmetric matrix  $\left[a_{p,q}\right]$  where for $p<q$,
$$
a_{p,q} = \sum {j_p+j_q\choose j}\qquad 
(\hbox{the sum over } j_q  < j \le j_p  ) ,
$$
and if $n$ is odd, then 
$(\!(J)\!):= \sum (-1)^{p-1}\ 2^{j_p}\ (\!(J\setminus\{j_p\})\!)$.
\ssk
\j iii;{\rm [La-La-T] \& [P3]} The total Segre class of the second exterior 
power $\Lambda^2(E)$ is given by
$$
s\bigl(\Lambda^2(E)\bigr) = \sum \ \left[I+\rho_{n-1}\right]\ s_I(E) 
$$
where the sum is over all partitions $I$ and the definition of $[J]$,
for $J=(j_1>\ldots>j_n\ge 0)$ is as follows. 
If $n$ is even, define $[J]$ to be the Pfaffian of the 
$n\times n$-antisymmetric matrix  
$\left[(j_p+j_q-1)!(j_p-j_q)/j_p! j_q!\right]$; 
if $n$ is odd then $[J]=0$ unless  $j_n=0$  where\ $[J] = [j_1,\dots,j_{n-1}]$.
\endproclaim

\proclaim{Remark 2.5 \rm(Background)}\rm  The history of the above formulas for 
$s(S^2E)$ and $s\bigl(\Lambda^2(E)\bigr)$ is as follows. 
At first, one of the authors of [La-La-T] has informed the author 
about recursive formulas for $(\!(J)\!)$ and $[J]$, in the form of linear
equations, obtained with the help of divided differences. 
(We will explain and use this extremely powerful technique 
in Sections 3, 4 and 6.) 
Using this recursion the author has found and proved the above 
Pfaffian formulas in [P3]. 
Finally, the authors of [La-La-T] managed to give a self-contained and 
elegant account 
of different formulas for $s(S^2E)$ and $s\bigl(\Lambda^2(E)\bigr)$ 
based on an interplay 
between the recursive formulas, Pfaffian expressions from [P3]  
and formulas which present $(\!(J)\!)$ and $[J]$ as sums of minors in some 
matrices of binomial numbers. 
Consequently, there are no divided differences 
in the final version of [La-La-T].
(``The power was eliminated by the elegance"!\footnote{The proof
of the linear equations for $(\!(J)\!)$ and $[J]$ 
via the divided differences is reproduced in Appendix A.3.})
\endproclaim

As it was mentioned in Section 1, the coefficients 
$D_{I,J}^{m,n}$ appearing in Theorem 2.4(i) are needed for the
expression of the Chern-Schwartz-MacPherson classes 
of $D_r(\varphi)$ associated with an $r$-general morphism.
 
The analogue of Theorem 1.9 for homomorphisms
with symmetries is not known yet; let us state, however, a weaker
result using, this time, the numbers $(\!(J)\!)$ and $[J]$
from Theorem 2.4(ii) and (iii). By $\iota$ we understand the inclusion
$D_r(\varphi)\to X$. 

\proclaim{Theorem 2.6 \rm [P3]}Assume that a (possibly singular)
complex analytic variety $X$ is compact, 
$\varphi$ is $r$-general and $D_{r-1}(\varphi)=\emptyset$.
{\parindent=\wd1
\j i;If $\varphi$ is symmetric then
$$
\iota_*(c_*(D_r(\varphi))) =
\sum (-1)^{|I|}(\!(I+\rho_{n-r-1})\!) Q_{\rho_{n-r}+I}(E)\cap c_*(X),
$$
the sum over all partitions of $I$ of length  $\le n-r$. 
\j ii;If $\varphi$ is antisymmetric, $r$ even, then
$$
\iota_*(c_*(D_r(\varphi)))=
\sum (-1)^{|I|}\ [I+\rho_{n-r-1}]\ P_{\rho_{n-r-1}+I}(E)\cap c_*(X),
$$
the sum over all partitions of $I$ of length  $\le n-r$.
\par}\endproclaim

Taking the degree of the expression on the right-hand side gives the topological Euler-Poincar\'e 
characteristic of $D_r(\varphi)$.
It would be valuable to extend the theorem to $r$-general morphisms
without the assumption of the emptiness of $D_{r-1}(\varphi)$.
\vs

Another application of Theorem 2.4 was given in [La-La-T] to the enumerative
properties of complete correlations and quadrics. 
Let us limit ourselves to the latter case. Here we assume
that the ground field is of characteristic different from 2.

Let us fix a positive integer $r$ and a projective space $\Bb P$. 
By a {\it complete quadric of rank} $r$ we understand a sequence  
$Q_{\bullet}: Q_1 \subset Q_2 \subset \dots\subset Q_n$  
($n$ can vary) of quadrics in $\Bb P$, such that
\smallskip
{\parindent=20pt
\item{1)} $Q_1$  is nonsingular,
\item{2)} the linear span $L(Q_i)$ of $Q_i$ is the vertex of $Q_{i+1}$, 
$i=1,\dots,n-1$,
\item{3)} $\dim L(Q_n) = r-1$.
\par}
\ssk
\noindent
There exists a natural structure of a nonsingular algebraic projective variety 
on $CQ(r)$ --- the set of all rank $r$ complete quadrics (see [La-La-T]
and the references therein). 
Let $\mu_i \in A_*(CQ(r))$ $(i=1,\dots,r)$ be the class of the locus of 
all complete quadrics $Q_{\bullet}$ such that $Q_n$ is tangent to a given 
(codimension $i$)-plane in $\Bb P$. 

Now let $\Cal G=G_r(\Bb P)$ be the Grassmannian parametrizing $(r-1)$-dimensional 
linear subspaces of $\Bb P$. 
Fix a sequence $I=(1\le i_1<i_2<\dots<i_r\le \dim \Bb P)$ of integers and 
consider the flag  $L_{\bullet}: L_1\subset L_2\subset \dots\subset L_r$ 
of linear subspaces in $\Bb P$ where $\dim L_j=i_j$, $j=1,\dots,r$. 
Let $\Omega (I)$ be the class in $A_*(\Cal G)$ of the Schubert cycle
$$
\{L\in \Cal G \mid \dim (L\cap L_j)\ge j-1,\ j=1,\dots,r\}.
$$ 
We have a map  $f: CQ(r) \to \Cal G$ such that
$f(Q_{\bullet}) = L(Q_n)$. 
Let $\omega (I) := f^*\Omega (I)$. 

Classics of enumerative geometry like Schubert, Giambelli \dots were interested 
in the computation of the number of complete quadrics $Q_{\bullet}$ such that $Q_n$ 
is tangent to $m_j$ fixed planes of codimension $j$ in general position 
in $\Bb P$ and such that 
$\dim (L(Q_n)\cap L_j)\ge j-1$ for each member of the above flag $L_{\bullet}$.
This question makes sense if \ $i_1+\dots+i_r+r-1 = m_1+\dots+m_r$  \ because then 
$\mu_1^{m_1} \mu_2^{m_2} \dots \mu_r^{m_r}\cdot\omega(I)$ is in $A_0(CQ(r))$. 
The answer to the question (under the above assumption) needs besides the 
numbers $(\!(J)\!)$ defined at the beginning of this section, also the function
$\alpha (p;k,j)$ defined by 
$$
\alpha(p;k,j):= \cases
\displaystyle{{k\choose 0} + {k\choose 1} p + \dots + {k\choose j} p^j}
& \hbox{if }j\ge 0,\\
0 & \hbox{--- otherwise.}\\
\endcases
$$
In fact, the following result answers a more general question:

\proclaim{Theorem  2.7 \rm [La-La-T]}Assume 
that $p$ is a number such that  $0\le p<r$  and 
$m_1+\dots+m_q > i_r+i_{r-1}+\dots+i_{r-q+1}+q-1$  for $q=1,\dots,p-1$. 
Then 
$$\displaylines{
\quad \mu_1^{m_1}\ \mu_2^{m_2}\ \dots\ \mu_{p+1}^{m_{p+1}}\cdot \omega (I) = 
1^{m_1}2^{m_2}\dots p^{m_p}\times \hfill\cr
\hfill \times \left[ (p+1)^{m_{p+1}}(\!(I)\!) - 
\sum \alpha \bigl(p;m_{p+1},m_{p+1}-|J|-(r-p)\bigr)\varepsilon_J 
(\!(J)\!) (\!(J')\!)\right],\cr}
$$
where the sum is over all subsequences $J$ in $I$ of cardinality $r-p$; 
$J'=I\setminus J$ and $\varepsilon_J=sign(J,J')$.
\endproclaim

This theorem generalizes and offers a ``modern treatment" of the results
of Schubert [S] and Giambelli [G2] from the end of the previous and
the beginning of the present century.
For more on this subject, consult also the paper [Th] by A. Thorup in the 
present volume. 

There is a similar formula for complete correlations which, in turn,  
uses the numbers $D^{n,m}_{I,J}$   (see [La-La-T]).

\head 3. Flag  degeneracy  loci and divided differences.\endhead
This section summarizes mainly some of the results of [F2].
Let
$$
F_{\bullet} : \ F_1\subset F_2\subset \dots \subset F_m= F \qquad\hbox{ and } \qquad
E^{\bullet}: \ E=E^n\twoheadrightarrow\dots\twoheadrightarrow E^2\twoheadrightarrow E^1
$$
be two flags of vector bundles
over a variety $X$ and let $\varphi:F\to E$ be a vector bundle homomorphism. 
Assume that a function  $\bfr :\{1,\dots,n\}\times\{1,\dots,m\}\to \Bb N$ 
is given (we will refer to $\bfr$ as to a {\it rank function}). 
Define 
$$
D_{\bfr}(\varphi)=\{x\in X\mid \rank \bigl(F_q(x)\to E^p(x)\bigr)\le 
\bfr(p,q) \ \forall p,q \}.
$$
In [F2] the author gives conditions on  $\bfr$ which guarantee that for 
a ``generic" $\varphi$, $D_{\bfr}(\varphi)$ is irreducible. 
Then, a natural problem arises, to find for such an  $\bfr$ and $\varphi$ 
a formula expressing $[D_{\bfr}(\varphi)]$ in terms of the Chern 
classes of $E^{\bullet}$ and $F_{\bullet}$ . 

It turns out that the crucial case is the case of complete flags, i.e. 
$\rank E^i=\rank F_i=i$  and  $m=n$. 
The desired formula in all other cases can be deduced from that one. 
In this situation, the degeneracy loci  $D_{\bfr}(\varphi)$  are 
parametrized by permutations $\mu\in S_n$, and 
$$
\bfr_{\mu}(p,q) = \card \{i\le p \mid \mu(i)\le q\}.
$$

Let $ \Omega_{\mu}(E^{\bullet},F_{\bullet}) = D_{\bfr_{\mu}}(\varphi)$. 
Then the expected (i.e. the maximal one, if the locus is nonempty) 
codimension of $\Omega_{\mu}(E^{\bullet},F_{\bullet})$ is $l(\mu)$ (the length of $\mu$). 
In order to describe a formula for the fundamental class of 
$\Omega_{\mu}(E^{\bullet},F_{\bullet})$ associated with a generic $\varphi$ we need some 
algebraic tools developed in  [B-G-G], [D1,2]  and [L-S 2,3] (for an
elegant account of this theory, see [M2]).
 
Let $A=(a_1,\dots,a_n),\ B=(b_1,\dots,b_n)$  be two sequences of independent and 
commuting variables. 
We have divided differences
$$
\partial_i:\Bb Z[AB]\to \Bb Z[AB]\qquad (\hbox{of degree }-1)
$$
defined by 
$$
\partial_i(f)=(f - \tau_if)/(a_i- a_{i+1})\qquad i=1,\dots,n-1,
$$
where $\tau_i=(1,\ldots,i-1,i+1,i,i+2,\ldots,n)$ 
denotes the $i$-th simple transposition. 
For every reduced decomposition  $\mu=\tau_{i_1}\cdot \dots\cdot \tau_{i_k}
$\footnote{This --- most common --- notation means that 
$\mu=(\mu(1),\dots,\mu(n))\in S_n$ is obtained from $(1,\dots,n)$ by the sequence 
of simple transpositions of components, where one performs first 
$\tau_{i_1}$, then $\tau_{i_2}$  etc.}  
one defines $\partial_{\mu}:= \partial_{i_1} \circ \dots \circ \partial_{i_k}$ 
--- an operator on $\Bb Z[AB]$ of degree $-l(\mu)$. 
In fact $\partial_{\mu}$ does not depend on the reduced decomposition chosen. 
Finally, for a permutation $\mu\in S_n$, we give, following [L5]
(see also [M2]):

\proclaim{Definition 3.1}\rm \underbar
{\it (Double) Schubert polynomials of 
Lascoux and Sch\"utzenberger}.\hfill\break  We set
$$
\frak S_{\mu}(A/B)=\partial_{\mu^{-1}\omega}\prod\limits_{i+j\le n}(a_i-b_j),
$$
where $\omega$ is the permutation of biggest length in $S_n$.
\endproclaim

Equivalently, the polynomials $\frak S_{\mu}(A/B)$ are defined inductively
by the equation
$$
\partial_i \bigl(\frak S_{\mu}(A/B)\bigr) =
\frak S_{\mu \tau_i}(A/B)
$$
if $\mu(i)>\mu(i+1)$, the top polynomial $\frak S_{\omega}(A/B)$ being 
$\prod\limits_{i+j\le n}(a_i-b_j)$.
\ssk

Note that the operators act here on the $A$-variables; 
however, it can be shown ([L5], [M2]) that 
$$
\frak S_{\mu}(A/B)=(-1)^{l(\mu)} \frak S_{\mu^{-1}}(B/A).
$$

Specialize now
$$
a_i:= c_1\bigl(\Ker  (E^i\to E^{i-1})\bigr) 
\quad\hbox{\rm and }\quad b_i:=c_1(F_i/F_{i-1}).
$$
 
Then we have

\proclaim{Theorem 3.2  \rm [F2]}Assume that $X$ is a pure-dimensional
Cohen-Macaulay scheme
and $\Omega_{\mu}(E^{\bullet},F_{\bullet})$ is of 
pure codimension $l(\mu)$ in $X$ or empty. Then the following equality 
holds in $A_*(X)$,
$$
[\Omega_{\mu}(E^{\bullet},F_{\bullet})]= \frak S_{\mu}(A/B)\cap [X].
$$
\endproclaim

The key point of the proof of Theorem 3.2 in [F2] is a geometric 
interpretation of the divided differences with the help of some
correspondences in flag bundles.
More precisely, assume, for simplicity, that $E\to X$ is a vector bundle
over a nonsingular variety $X$ and let $\Cal F\to X$ be the flag bundle
parametrizing the flags of quotients of $E$ of successive ranks
$n,n-1,\dots,2,1$.
Denote by
$$
E^{\bullet}: E=E^n\toto E^{n-1}\toto\dots\toto E^2\toto E^1
$$
the tautological flag on $\Cal F$.
It is well known that for $a_i=c_1$ $\bigl( \Ker (E^i\toto E^{i-1})\bigr)$,
$i=1,\dots,n$, $A^*(\Cal F)$ is a quotient ring of $A^*(X)[a_1,\dots,a_n]$.
Let $\Cal F(i)$ be the flag bundle parametrizing the
flags of quotients of \ $E$ \ of successive ranks \ \ 
$n,n-1,...,i+1,i-1,...,2,1$. Consider a fibre square
$$
\Cal F \times_{\Cal F(i)} \Cal F 
$$
Denote by
$$
p_1,p_2: \Cal F \times_{\Cal F(i)} \Cal F \to \Cal F
$$
the projections onto the successive factors.
These projections are $\Bbb P^1$-bundles.

\proclaim{Proposition 3.3 \rm [F2]}\parindent\wd1
\j i;The map $(p_1)_*\circ p_2^*:A^k(\Cal F)\to A^{k+1}(\Cal F)$
acts on polynomials in $a_1,\dots,a_n$ \ like the divided-differences 
operator $\partial_i$ does.
\j ii;Assume that a flag of subbundles
$$
F_{\bullet}:F_1\subset F_2\subset \dots\subset F_{n-1}\subset F_n=E
$$
is given on $X$.
Then, in $A^*(\Cal F)$,
$$
(p_1)_*\circ p_2^*\bigl[\Omega_{\mu}(E^{\bullet},(F_{\bullet})_{\Cal F})\bigr]=
\bigl[\Omega_{\mu \tau_i}(E^{\bullet},(F_{\bullet})_{\Cal F})\bigr]
$$
if $\mu(i)>\mu(i+1)$, and $0$ --- otherwise.
\endproclaim

The theorem generalizes in a uniform way the formulas for the fundamental 
classes of Schubert varieties in the flag varieties from [B-G-G] and [D2], 
and --- with the help of a rich algebra of Schubert polynomials (see [M2])
--- some other known before formulas for flag degeneracy loci like the 
Giambelli-Thom-Porteous formula (see Section 1) as well as determinantal 
formulas for flag degeneracy loci from [K-L], [L1] and [P3] which we recall in
the following examples. (Note that another approach to the 
Giambelli-Thom-Porteous formula, this time using a certain Schur complex, is
given in Appendix A.6.) 

\proclaim{Example 3.4 \rm [K-L]}\rm Assume that on $X$ a flag of vector
bundles
$$ 
B_1 \subsetneqq B_2 \subsetneqq \dots \subsetneqq B_k=B
$$
is given, with rank $B_i=m_i$. Moreover, let $\varphi: A\to B$ be a vector
bundle homomorphism where $\rank A=n$. Consider the locus:
$$
\Omega = \bigl\{ x\in X| \ \dim \Ker  (B_i(x)\hookrightarrow B(x)
@>\varphi (x) >> A(x))\ge i , \  i=1,\ldots,k \bigr\}.
$$
\endproclaim

Then, assuming that $X$ is a pure dimensional 
Cohen-Macaulay variety and $\Omega$ is of pure
codimension $\sum_i (n-m_i+i)$ in $X$ or empty, one has the equality
$$
[\Omega] = \Det \Bigl[ c_{n-m_i+j}(A-B_i) \Bigr]_{1\le i,j \le k}\cap [X]
$$
(compare also [L1]). The author of [F2] reports too restrictively on p.~417 that his approach
does not cover all the instances of this formula but only the cases of the
form $n-m_1+1>n-m_2+2>\ldots>n-m_k+k$. (A similar remark applies to the
formula treated in the next example.)

\proclaim{Example 3.5 \rm [P3, (8.3)]}\rm Assume that on $X$ two flags 
of vector bundles are given
$$ 
B_1\subset B_2\subset \dots \subset B_k=B, \qquad 
A=A_1 \twoheadrightarrow A_2 \twoheadrightarrow \dots \twoheadrightarrow 
A_{k-1} \twoheadrightarrow A_k
$$
with $\rank A_i=n_i$, $\rank B_i=m_i$. Moreover, let $\varphi: B\to A$
be a vector bundle homomorphism. Consider the locus
$$
\Omega = \bigl\{ x\in X \mid \dim \Ker \bigl(B_i(x)\hookrightarrow 
B(x)@> \varphi(x) >> A(x) 
\twoheadrightarrow A_i(x) \bigr)\ge i , \ i=1,\ldots,k \bigr\}.
$$

Then, assuming that $X$ is a pure-dimensional Cohen-Macaulay variety,
$m_i\ge i$, 
$$ 
n_1-m_1+1 \ge n_2-m_2+2 \ge \dots \ge n_k-m_k+k>0
$$
and $\Omega$ is of pure codimension $\sum (n_i-m_i+i)$ in $X$ or empty, one
has the equality
$$
[\Omega]=\Det \Bigl[ c_{n_i-m_i+j}(A_i - B_i) \Bigr]_{1\le i,j \le k}\cap [X].
$$

See Appendix A.4 for a proof of this formula with the use of Gysin maps.
\endproclaim

A combination of Theorem 3.2 with [G5] gives some interesting formulas
for specializations of indeterminates in Schubert polynomials (see [F2, p.419];
compare also some related computations in [He-T] using the Gr\"obner
bases technique).

\vs
Finally, note that Schubert polynomials are a useful tool
in the computation of Chern classes of the tangent vector bundles to the flag 
varieties --- for details see [L5]. 
For one more application of Schubert polynomials, this time to 
the cohomology rings of Schubert varieties, see [A-L-P].

\head 4. Gysin maps and divided differences.\endhead
As it was pointed out in [F2], the divided differences $\partial_i$ are
geometrically constructed from correspondences which are $\Bb P^1$-bundles
(see also the preceding section).

The aim of this section is to emphasize that some compositions of the
$\partial_i$'s can be interpreted geometrically as Gysin maps for
flag bundles.
Similar results are true for divided differences associated with other
semisimple algebraic groups. We also state some ``orthogonality" results
with respect to Gysin maps for flag bundles.

Let $\pi:G^1(E)\to X$ be the projective bundle parametrizing 1-quotients of $E$
where $E$ is a vector bundle on a variety $X$ of rank $n$.
Assume, for simplicity, that $X$ is smooth.
Let $A=(a_1,\dots,a_n)$ be a sequence of independent indeterminates and
$\alpha_1,\dots,\alpha_n$ --- the sequence of Chern roots of $E$.
One has the divided-differences
operators $\partial_i:\Bb Z[A]\to\Bb Z[A]$ ($i=1,\ldots,n-1$), associated with the simple transpositions, defined by the
formulas from the preceding section.

Denote by $\tr:\Bb Z[A]\to\Bb Z[A]$ the following composition of 
divided-differences operators
$$
\tr:=\partial_{n-1}\circ\dots\circ\partial_2\circ\partial_1
$$

We emphasize that while the divided-differences operators act on the 
whole polynomial ring $\Bb Z[A]$, the symmetrizing operators appearing
below are, in general, defined on proper subrings of $\Bb Z[A]$
(for more on symmetrizing operators, see [L-S4] and [P4]). 

\proclaim{Proposition 4.1}One has, for $P\in\Bb Z[A]^{S_1\times S_{n-1}}$,
$$
\tr\!\! P=\sum\limits_{\mathstrut\overline\sigma\in S_n/S_1\times S_{n-1}}
\sigma\left({P\over \prod\limits_{\beta\ge 2}(a_1-a_{\beta})}\right)
\tag i
$$
$$
\pi_*\bigl(P(\alpha_1,\dots,\alpha_n)\bigr)=(\tr\!\! P)(\alpha_1,\dots,\alpha_n).
\tag ii
$$
\endproclaim

A word about how one proves such a result (this method can also be 
applied to other results of this type stated in this section).
The equality (i) is just a straightforward verification.
In fact, since $\tr\!\!(\sigma P)=\sigma\tr\!\!(P)$ for $\sigma\in S_n$,
it suffices to check it for $P=a_1^k$, $k=0,1,\ldots,n-1$.
For the degree reasons it remains to show (i) for
$P=a_1^{n-1}$, the calculation of the expression on the right-hand side 
being essentially the Laplace development of the Vandermonde
determinant.
To show (ii) we can assume without loss of generality that $X$ is a ``big"
Grassmannian and $E$ is a universal bundle on it
with the Chern roots $a_1,\ldots,a_n$.
Then $\tr$  \  induces an $A^*(X)$-morphism.
Since $A^*(G^1E)=\bigoplus^{n-1}_{i=1} \pi^*A^*(X)\xi^i$,
where $\xi=c_1(\Cal O(1))$, the assertion follows from the equality 
$\tr\!\!(a_1^{n-1})=1$
because $\pi_*(\xi^{n-1})=1$ and $\pi_*(\xi^i)=0$ for $i<n-1$.

Let $\pi:G^q(E)\to X$ be the Grassmannian bundle parametrizing $q$-quotients 
of the bundle $E$ as above.
Denote by $\square:\Bb Z[A]\to\Bb Z[A]$ the following composition of
divided-differences operators ($r=n-q$): 
$$
\square:=(\partial_r\circ\dots\circ\partial_2 \circ\partial_1)\circ \dots\circ 
(\partial_{n-2}\circ\dots\circ\partial_q\circ\partial_{q-1})
\circ
(\partial_{n-1}\circ\dots\circ\partial_{q+1}\circ\partial_q) .
$$ 

\proclaim{Proposition 4.2}One has, for $P\in\Bb Z[A]^{S_q\times S_r}$,
$$
\square P=\sum\limits_{\mathstrut\overline{\sigma}\in S_n/S_q\times S_r}
\sigma\left({P\over\prod\limits_{\alpha\le q<\beta}(a_{\alpha}-a_{\beta})}
\right)
\tag i
$$
$$
\pi_*\bigl(P(\alpha_1,\dots,\alpha_n)\bigr)=(\square P)(\alpha_1,\dots,\alpha_n).
\tag ii
$$
\endproclaim

The operator $\square$ is known in interpolation theory as the
{\it Lagrange-Sylvester symmetrizer}.
We refer the reader to [L4] for an account of algebraic properties
of the Lagrange-Sylvester symmetrizer.

Let now $\tau=\tau_E:\Fl (E)=\Fl ^{n,n-1,\dots,1}(E)\to X$ be the flag bundle 
endowed with the tautological sequence of quotients
$$
E=Q^n\toto Q^{n-1}\toto\dots\toto Q^2\toto Q^1,
$$
where $\rank Q^i=i$.
Let $L_i=\Ker (Q^i\toto Q^{i-1})$ and $\alpha_i=c_1(L_i)$, $i=1,\dots,n$.

Denote by $\partial=\partial_A:\Bb Z[A]\to\Bb Z[A]$ the following
composition of divided differences:
$$
\partial:=(\partial_1\circ\dots\circ\partial_{n-1})\circ\dots\circ
(\partial_1\circ\partial_2\circ\partial_3)\circ(\partial_1\circ\partial_2)
\circ\partial_1
$$
In other words, $\partial=\partial_{\omega}$ is the operator associated with
the longest element $\omega\in S_n$ (in the notation of Section 3).

\proclaim{Proposition 4.3}One has, for $P\in\Bb Z[A]$,
$$
\partial P=\sum\limits_{\sigma\in S_n}\sigma
\left({P\over\prod\limits_{\alpha<\beta}(a_{\alpha}-a_{\beta})}\right) 
\tag i
$$
$$
\tau_*\bigl(P(\alpha_1,\dots,\alpha_n)\bigr)=(\partial P)(\alpha_1,\dots,\alpha_n).
\tag ii
$$
\endproclaim

Let us prove (ii).
Consider the factorization
$$
\tau=\tau_E: \Fl (E) \cong \Fl (Q^{n-1})@>\tau'=\tau_{Q^{n-1}}>>
G_1E @>\pi>> X.
$$
Hence $\tau_*=\pi_*\circ\tau_*'$.
Define the operator $\partial':\Bb Z[A]\to\Bb Z[A]$ by
$$
\partial'=(\partial_1\circ\dots\circ\partial_{n-2})\circ\dots
\circ(\partial_1\circ\partial_2)\circ\partial_1.
$$
Thus $\partial=\tr\!\!\circ\partial'$. Then assuming that equation (ii) is
valid for $\tau'$ and $\partial'$ instead of $\tau$ and $\partial$, and
invoking Proposition 4.1 we deduce, by induction on $n$, that equation (ii)
holds.
\vs

The operator $\partial$ is called the {\it Jacobi symmetrizer}.
A familiar Jacobi-Trudi identity can be restated as

\proclaim{Proposition 4.4}Let $I$ be a partition, $l(I)\le n$. Then
$$
\partial(a_1^{i_1+n-1}a_2^{i_2+n-2}\dots a_{n-1}^{i_{n-1}+1}a_n^{i_n})=s_I(A).
$$\endproclaim

For a simple operator proof %of the proposition, 
see Appendix A.5.
In fact this identity is valid for sequences $I\in \Bb Z^n$ such that
$i_1\ge-(n-1),\dots,i_{n-1}\ge-1$, $i_n\ge 0$ (notation: $I\ge-\rho_{n-1}$).
\vs

Using this proposition one can give a short operator proof of Proposition 
1.3(i).
Denoting $A_q=(a_1,\dots,a_q)$ and $A^r=(a_{q+1},\dots,a_n)$ and assuming
$\Bb Z^q\ni I\ge-\rho_{q-1}$, $\Bb Z^r\ni J\ge-\rho_{r-1}$, we obtain
$$\aligned
\square\bigl(s_I(A_q)\cdot s_J(A^r)\bigr)
&=\square\Bigl(\partial_{A_q}(a_1^{i_1+q-1}\cdot \dots\cdot a_q^{i_q})\cdot
\partial_{A^r}(a_{q+1}^{j_1+r-1}\cdot \dots \cdot a_n^{j_r})\Bigr)  \\
&=\partial_A(a_1^{i_1+q-1}\cdot \dots\cdot a_q^{i_q}
a_{q+1}^{j_1+r-1}\cdot \dots\cdot a_n^{j_r})  \\
&=\partial_A\bigl(a_1^{(i_1-r)+(n-1)}\cdot \dots\cdot a_q^{(i_q-r)+(n-q)}
a_{q+1}^{j_1+r-1}\cdot \dots\cdot 
a_n^{j_r}\bigr) \\
&=s_{i_1-r,\dots,i_q-r,j_1,\dots,j_r}(A),
\endaligned$$
which is the algebraic essence of the proposition (the second equality
being a straightforward verification).
\vs

Let us now pass to other root systems than the one of type $A_{n-1}$. 
We need more
notation which will also be used in Section 6.

Let $G$ be a semisimple algebraic group over $\Bbb C$, $B\subset G$ -- a Borel subgroup and 
$T\subset B$ -- a maximal torus.
With a character $\chi\in X(T):=\Hom (T,\Bb C\,\strut^*)$ one associates 
a line bundle $L_{\chi}$ over a generalized flag variety $G/B$.
The total space of $L_{\chi}$ is $G\times\Bb C/\tilrel$, where the 
relation ``$\tilrel$" is defined by:
$(g,z)\tilrel(gb,\chi(b^{-1})z)$ for $g\in G$, $z\in\Bb C$, $b\in B$;
recall that $X(T)=\Hom(B,\Bb C^*)$.
We have a map from $X(T)$ to $A^1(G/B)$ which assigns to $\chi$
the Chern class $c_1(L_{\chi})$.
It extends multiplicatively to the Borel characteristic map:
$$
c:S^{\bullet}(X(T))\to A^*(G/B).
$$

Let $R$ be the root system of $(G,T)$, endowed with the basis $\tr$ associated with $B$,
and $P\supset B$ be the parabolic subgroup associated with a subset 
$\theta\subset\tr$.
Denote by $W_{\theta}$ the subgroup of the Weyl group $W$
of $(G,T)$ or $R$
generated by all reflections $\{s_{\alpha}\}_{\alpha\in\theta}$.
Then the characteristic map restricted to the $W_{\theta}$-invariants 
factorizes through $A^*(G/P):$
$$
c:S^{\bullet}(X(T))^{W_{\theta}} @>>> A^*(G/P).
$$
Let $R_\theta$ (resp. $R^+_{\theta}$) denote the subset of the set of roots
(resp. positive roots with respect to $B$)
formed by linear combinations of roots in $\theta$ and
$d_{\theta}=\prod\limits_{\alpha\in R_{\theta}^+}\alpha$.
Consider a ``symmetrizing operator" \hbox{$\partial^{\theta}: S^{\bullet}(X(T))\to
S^{\bullet}(X(T))$} defined by
$$
\partial^{\theta}(f)
=\sum\limits_{w\in W_{\theta}}(-1)^{l(w)}w(f)/d_{\theta},
$$
where $l(w)$ is the length of $w$ with respect to $\tr$.

\proclaim{Theorem 4.5 \rm[A-C]}\par
{\parindent =\wd1
\line{\hbox to\wd1{\hfil\rm(a)\enspace}\hfill $\partial^{\theta}
(S^{\bullet}(X(T)))\subset 
S^{\bullet}(X(T))^{W_{\theta}}$.\hfill}
\item{\rm(b)} The diagram
$$\CD
S^{\bullet}(X(T)) @>c>> A^*(G/B) \\
@VV\partial^{\theta}V @V\pi_*VV \\
S^{\bullet}(X(T))^{W_{\theta}} @>c>> A^*(G/P)
\endCD $$
commutes, where $\pi_*$ is the Gysin map associated with $\pi:G/B\toto G/P$.
In other words, $\pi_*$ is induced by $\partial^{\theta}$.
\par}\endproclaim

Observe that for $\theta=\tr$, 
$\partial^{\tr}(f)=\sum(-1)^{l(w)}w\,f/d$ where $d$ is the product 
of all positive roots and the sum is over $w\in W$.
This formula gives, more generally, a symmetrizing operator
description of the Gysin maps associated with $G/B$-fibrations\footnote{Note 
added in proof: This has circulated
as a ``folklore" among specialists since some time.
A precise written account for $f_*: A_*(X/B) \to A_*(Y)$, where
$f: X\to Y$ is a principal $G$-bundle ($X$ and $Y$ are
nonsingular varieties, $G$ is a reductive group), is contained in
a paper by M.~Brion {\it The push-forward and Todd class
of flag bundles} -- in this volume. In this paper, the author
gives also a symmetrizing operator description of the 
Gysin map associated with a flag bundle $X/P\to Y$ where
$P\subset G$ is a parabolic subgroup (thus generalizing
the general linear group case from [P3] and the results of [A-C]).
Moreover, he computes the Todd class of the tangent bundle
of such a ``flag fibration" $X/P \to Y$. }  overlapping, e.g., the case of
Lagrangian and orthogonal flag bundles.

Let us recall another familiar interpretation of $\partial^{\tr}$.
Given a root $\alpha\in R$, one defines an operator
$\partial_{\alpha}:S^{\bullet}(X(T))\to S^{\bullet}(X(T))$ by
$$
\partial_{\alpha}f=(f-s_{\alpha}f)/\alpha,
$$
where $s_{\alpha}\in W$ is the reflection associated with $\alpha$.

\proclaim{Lemma 4.6 \rm [B-G-G], [D1]}If $w\in W$, 
$l(w)=k$ and $w=s_{\alpha_1}\cdot \ldots \cdot s_{\alpha_k}
=s_{\beta_1}\cdot \ldots \cdot s_{\beta_k}$ where $\alpha_i,\beta_i\in\tr$,
then
$$  
\partial_{\alpha_1}\circ \ldots \circ \partial_{\alpha_k}
=\partial_{\beta_1}\circ \ldots \circ \partial_{\beta_k}.
$$
\rm Thus the value of this operator can be denoted by $\partial_w$ without
ambiguity.
\endproclaim

It has been shown in loc.cit. that for the longest element $w_0\in W$,
$\partial_{w_0}=\partial^{\tr}$.
Since $\theta$ is a basis of the root (sub)system $R_{\theta}\subset R$,
the operator $\partial^{\theta}$ is similarly interpreted as the operator
associated with the longest word in $W_{\theta}$.
(Of course, the operators considered in Section~3
and in the beginning of this section
are special cases of the
ones here for type $A_{n-1}$, the Jacobi symmetrizer being 
$\partial^{\tr}=\partial_{(n,n-1,\dots,1)}$).

\vs

 Taking as 
the starting point the above mentioned 
symmetrizing operator description of the Gysin map 
associated with a $G/B$-fibration, the authors of [P-R5]
give some explicit formulas for the Gysin maps associated with Lagrangian
and orthogonal Grassmannian bundles. Consider, for
example, the Lagrangian case.
Let $V\to X$ be a vector
bundle of rank $2n$ endowed with a nondegenerate antisymmetric form.
Let $\tau: \operatorname{LFl}(V)\to X$ be the flag bundle parametrizing flags
of isotropic subbundles of ranks $1,2,\ldots,n$ in $V$ with respect to 
the above
form, and let $\pi: LG_n(V)\to X$ be the Lagrangian Grassmannian bundle
parametrizing top-dimensional isotropic subbundles in $V$. Let
$\alpha_1,\ldots,\alpha_n$ be the sequence of the Chern roots of the 
tautological Lagrangian (sub)bundle $R$ on $LG_n(V)$.
We know by the above that one has in $A^*(X)$,
$$
\tau_* \bigl(f(\alpha_1,\ldots,\alpha_n)\bigr)=
\bigl(\partial_{(\oms 1,\oms 2,\ldots,\oms n)}f\bigr)
(\alpha_1,\ldots,\alpha_n),
$$
where $f$ is a polynomial in $n$ variables.
(Recall that the symplectic Weyl group is
the group of barred permutations and $(\oms 1,\oms 2,
\ldots,\oms n)$ is the longest element of it.)

In connection with $\pi_*$, it is natural to associate with a vector bundle
$E$ the following polynomials in its Chern classes (consult [P-R5] and the
subsection {\it A connection with $Q$-polynomials} in Section 6). We set
$\widetilde Q_iE=c_i(E)$ and for a partition $I$, $l(I)
\ge 2$, to define $\widetilde Q_IE$ we mimic the definition
of the Schur polynomial $Q_I$ (see Section 1).

\proclaim{Proposition 4.7 \rm [P-R5]}One has in $A^*(X)$:
{\parindent\wd1
\j i;$\pi_*\bigl(f(\alpha_1,\ldots,\alpha_n)\bigr) = 
\Bigl(\partial_{(\oms n,\oms{n-1},\ldots, \oms 2, 
\oms 1)}f\Bigr)(\alpha_1,\ldots,\alpha_n)$ for a symmetric
polynomial $f$ in $n$ variables.
\j ii;The element $\widetilde Q_IR\hak$ has a nonzero
image under $\pi_*$ only if each number $p$, $1\le p\le n$,
appears as a part of $I$ with an odd multiplicity $m_p$.
If the latter condition holds then
$$
\pi_*\widetilde Q_I R\hak = \prod_{p=1}^{n} \bigl((-1)^p
c_{2p}V\bigr)^{(m_p-1)/2}.
$$
\j iii;The element $s_IR\hak$ has a nonzero image under
$\pi_*$ only if $I$ is the partition of the form $2J+\rho_n$
for some partition $J$. If $I=2J+\rho_n$ then
$$
\pi_*s_IR\hak = s_J^{[2]}V,
$$
where the right-hand side is defined as follows: if $s_J=P(e.)$
is a unique presentation of $s_J$ as a polynomial
in the elementary symmetric functions $e_i$, $E$ -- a 
vector bundle, then $s_J^{[2]}(E):= P$ with $e_i$
replaced by $(-1)^ic_{2i}(E)$ \ $(i=1,2,\ldots)$.
\par}
\endproclaim

We finish this section with some examples of the ``orthogonality" results
for Gysin maps for flag bundles. These results are crucial for
computing the classes of diagonals in flag bundles
and the knowledge of the classes of diagonals is useful in calculations of the classes of Schubert
varieties, following a procedure described in the 
next section.

\vs

In the next theorem we follow the previously introduced 
notation from this section.
Additionally, we put $\frak S_{\mu}(A)=\frak S_{\mu}(A/0,\ldots,0)$%
\footnote{Note added in proof: It is an interesting question to extend
the above definition of Schubert polynomials $\frak S_{\mu}(A)$ to other
semisimple groups. It appears that for the symplectic and orthogonal
groups, a satisfactory  algebro-combinatorial theory of this type has been
given by S.~Billey and M. Haiman in {\it Schubert polynomials for classical
groups}, J. Amer. Math. Soc. 8 (1995), 443--482. 
Compare also: S. Fomin and A.~N. Kirillov, {\it Combinatorial $B_n$ analogues
of Schubert polynomials} and S.~Billey, {\it Transition equations for isotropic
flag manifolds} --- preprints (1995). For a theory of Schubert polynomials, suited to algebraic geometry, 
which has grown up from a paper by P.~Pragacz and J.~Ratajski [P-R5], see a forthcoming article:
{\it Symplectic and orthogonal Schubert polynomials \`a la polonaise} --- in preparation in collaboration with
A.~Lascoux. }
for $\mu \in S_n$, in the notation of Section~3.
Moreover, we define
the polynomials $\widetilde Q_IA$ as follows.
We set $\widetilde Q_iA=e_i(A)$ (the $i$-th elementary symmetric polynomial
in $A$), and, for a strict partition $I$, $l(I)\ge 2$, we  mimic 
the definition of the Schur $Q$-polynomial $Q_I$
(thus the element $\widetilde Q_IE$ associated above
with a bundle $E$ is equal to $\widetilde Q_IA$
with $A$ specialized to the Chern roots of $E$).
Note that the algebro-geometric 
properties of these $\widetilde Q$-{\it polynomials} were
worked out in [P-R5] (see also [La-Le-T2] where an interesting specialization
of a Hall-Littlewood polynomial, denoted in loc.cit. by $Q'_I$, 
is studied, giving, for the specialization $q=-1$ of the 
parameter, the ``Young dual" of $\widetilde Q_I$). 

\proclaim{Theorem 4.8}{\parindent\wd1
\j i;For partitions $I,J\subset (r)^q$,
$$
\square\bigl(s_I(A_q)\cdot s_{\overline J}(-a_{q+1},\dots,-a_n)\bigr) = \delta_
{I,J} ,
$$
where the Ferrers' diagram of $\overline J$ is the complement of the one
of $J^{\sim}$ in $(q)^r$. Using a standard $\lambda$-ring notation,
this is equivalently rewritten as
$$
\square\bigl(s_I(A_q)\cdot s_{(r)^q/J}(-A^r)\bigr) = \delta_
{I,J} .
$$
\j ii;For permutations $\mu, \nu\in S_n$,
$$
\partial_{\omega}\bigl(\frak S_{\mu}(A)\cdot \frak S_{\nu \omega}
(-a_n,-a_{n-1},\ldots, -a_1) \bigr) = \delta_{\mu,\nu}.
$$
\j iii;For strict partitions $I,J\subset \rho_n=(n,n-1,\ldots, 2,1)$, one has the
following equality, in the symplectic case:
$$
\partial_{(\oms n,\oms{n-1},\ldots,\oms 2,
\oms 1)}\bigl(\widetilde Q_I(A)\cdot \widetilde Q_{\rho_n \smallsetminus 
J}(A)\bigr) = \delta_{I,J}.
$$\par}
Here, $\delta_{\cdot,\cdot}$ denotes the Kronecker delta and $\rho_n 
\smallsetminus J$ is the strict partition whose parts complement the
parts of $J$ in $\{1,\ldots,n\}$.
\endproclaim

Assertion (i) can be deduced from (the $\square$-version of) 
Proposition 1.3(i). For a proof of (ii), see [L-S6], [M2] or [L-P].
Assertion (iii) stems from [P-R5].

\head 5. Fundamental classes, diagonals and Gysin maps.\endhead
As it was pointed out in Section 1, one of the fundamental problems
in the study of a concrete (closed) subscheme of a given (smooth) scheme $X$ is
the computation of its fundamental class in terms of given 
generators of the Chow ring of $X$. 

The decisive role in the method described in this section 
is played by the diagonal of the ambient scheme 
or, more precisely, its class in the corresponding Chow group of a fibre
product\footnote{The material of this section is due to the author.
A discovery of this method is inspired by the
construction used in the proof of the main formula in the paper 
by G. Kempf and D.~Laksov [K-L].}.
As a matter of fact, we have already seen, in Section~3, 
one application of the diagonal
to the computation of fundamental classes.  
In the situation of Section~3, there is a vector bundle on the product
of flag bundles endowed with a section vanishing precisely on the diagonal.
The top Chern class of this bundle is represented by 
$\frak S_{\omega}(A/B)$. By applying
divided differences to this polynomial, one gets polynomials representing
the classes of other (i.e. higher dimensional) 
degeneracy loci in the product of flag 
bundles. (This generalizes the procedure from [B-G-G] and [D2]: starting
from the class of the point and applying divided differences one gets
the class of a curve, then --- the class of a surface, etc.)
The procedure given below is of different nature. By using a desingularization
of the subscheme whose class we want to compute, and the diagonal of the
ambient scheme, we replace the original problem by the one of computing
the image of the class of the diagonal under an appropriate Gysin map.   
Moreover, since the diagonal is not always represented as the scheme
of zeros of a vector bundle (this seems to happen, e.g., for flag bundles
for classical groups different from $SL_n$), we give a recipe allowing
to calculate the class of the diagonal of the fibre product with the help
of Gysin maps. 
\medskip

Let $S$ be a smooth scheme (over a field) and $\pi:X\to S$ 
a smooth morphism of schemes.
Suppose that $D\subset X$ is a (closed) subscheme whose class is to be
computed.
Let $p:Z\to S$ be a proper smooth morphism and $\alpha:Z\to X$ be an
$S$-morphism which maps birationally $Z$ onto $D$. Assume that $\alpha$ is
proper. 
Consider a commutative diagram:

\vs

$$\aligned
\Delta\hskip 5pt \hoto \hskip 5pt&X\times_S X  \\
&\hskip 14pt\vbox{\offinterlineskip
     \hbox{$\Big\uparrow \hskip 3pt {\scriptstyle 1\times\alpha}$} } \\
&\vtop{\vskip -3pt \hbox{$X\times_S Z$} }
 \hskip 8pt 
 \vtop{\offinterlineskip 
     \hbox{$ @<\ \ \sigma\ <<$}
     \hbox{$ @>>\ p_2\ >$} } 
 \hskip 8pt 
 \vtop{\vskip -3pt \hbox{$Z$} \vskip 11pt} \\
&\hskip 14pt
 \vtop{\offinterlineskip
     \vglue-7pt\hbox{$\Big\downarrow\hskip 3pt {}^{p_1}$}}
 \hskip 15pt     
 \vtop{\offinterlineskip
     \hbox{${}_\alpha\hskip 0.6pt \dup$}
     \hbox{$\hskip 0.1pt\swarrow$}  }
 \hskip 10pt    
 \vtop{\offinterlineskip
     \vglue-7pt\hbox{$\Big\downarrow\hskip 3pt {}^p$}%\vskip 7pt 
     }  \\
D\hskip 5pt\hoto\hskip 5pt &\hskip 12pt X\hskip 21pt @>>\ \pi\ \ >\hskip 4pt S
\endaligned$$

\vs
Here $p_1$ and $p_2$ are the projections, the section $\sigma \ ($of $p_2)$ 
equals
$\alpha\times_S\id$ and $\Delta$ is the diagonal in the fibre product 
$X\times_S X$.

\proclaim{Proposition 5.1}Suppose that the class 
of the diagonal $\Delta$ in $A^*(X\times_SX)$ is 
$[\Delta]=\sum \pr_1^*(x_i)\cdot \pr_2^*(y_i)$ where $\pr_i:X\times_SX\to X$ are
the projections and $x_i,y_i\in A^*(X)$. 
Then, in $A^*(X)$,
$$
[D]=\sum\limits_ix_i\cdot\bigl(\pi^*\,p_*\,\alpha^*(y_i)\bigr).
$$
\endproclaim

\Proof
By the assumption $[D]=\alpha_*([Z])$.
Since $\alpha=p_1\circ\sigma$, we have $\alpha_*([Z])=(p_1)_*[\sigma(Z)]$.
Now, the key observation is that, in the scheme-theoretic sense, one has 
the equality $\sigma(Z)=(1\times\alpha)^{-1}(\Delta)$.
Since $\Delta \cong X$ is smooth, this implies
$[\sigma(Z)]=(1\times\alpha)^*([\Delta])$ \ (see Lemma 9 in [K-L]).
We then have:
$$\aligned
[D] &=(p_1)_*\Bigl([\sigma(Z)]\Bigr)=
(p_1)_*\Bigl((1\times\alpha)^*([\Delta])\Bigr)  \\
&=(p_1)_*\Bigl((1\times\alpha)^*
\bigl(\sum\limits_i\pr_1^*(x_i)\cdot \pr_2^*(y_i)\bigr)\Bigr)  \\
&=(p_1)_*\Bigl(\sum\limits_i p_1^*(x_i)\cdot p_2^*(\alpha^*(y_i)
)\Bigr)  \\
&=\sum\limits_i x_i\cdot\bigl((p_1)_*\,p_2^*\,\alpha^*(y_i)\bigr)\\
&=\sum\limits_i x_i\cdot\bigl(\pi^*\,p_*\,\alpha^*(y_i)\bigr),
\endaligned$$
where the last equality follows from the above fibre product diagram and
[F1, Proposition 1.7].
\endproof

The next result shows how one can compute the fundamental class
of the diagonal $[\Delta]\in A^*(X\times_S X)$ using Gysin maps.

\proclaim{Theorem 5.2}Let $S$ be as above and
$\pi:X\to S$ be a proper smooth morphism such that 
$\pi^*$ makes $A^*(X)$ a free $A^*(S)$-module; 
$A^*(X)=\bigoplus_{\alpha \in \Lambda} A^*(S)\cdot a_{\alpha}$,
where $a_{\alpha}\in A^{n_{\alpha}}(X)$ and
$A^*(X)=\bigoplus_{\beta \in \Lambda} A^*(S)\cdot b_{\beta}$,
where $b_{\beta}\in A^{m_{\beta}}(X)$.
Suppose that for any $\alpha$ there is a unique $\beta=:\alpha'$ 
such that $n_{\alpha}+m_{\beta}=\dim X - \dim S$ and 
$\pi_*(a_{\alpha}\cdot b_{\beta})\ne 0$ $($assume 
$\pi_*(a_{\alpha}\cdot b_{\alpha'})=1)$.
Moreover, denoting by $\pr_i:X\times_{_{\scs S}}X\to X$ $(i=1,2)$ the projections,
suppose that the homomorphism
$A^*(X)\otimes_{_{\scs A^*(S)}}A^*(X)\to  
A^*(X\times_{_{\scs S}}X)$, defined by $g\otimes h\mapsto \pr_1^*(g)\cdot \pr_2^*(h)$, is an
isomorphism. Then 

{\parindent\wd1
\j i; The class of the diagonal $\Delta$ in $X\times_{_{\scs S}}X$ equals
 \ $\big[\Delta\big] = \sum_{\alpha,\beta} d_{\alpha\beta} 
a_{\alpha}\otimes b_{\beta},$ \ where, 
for any $\alpha,\beta$, \ $d_{\alpha\beta}=P_{\alpha\beta}
(\{\pi_*(a_{\kappa} \cdot 
b_{\lambda})\})$ 
for some polynomial $P_{\alpha\beta}\in \Bb Z[\{x_{\kappa\lambda}\}]$.

\j ii;The following conditions are equivalent: 
\item{}{\rm a)} One has $\pi_* (a_\alpha \cdot b_{\beta'})=\delta_{\alpha,\beta}$, 
the Kronecker delta. 
\item{}{\rm b)} The class of the diagonal $\Delta \subset X \times_{_{\scs S}} X$ equals 
 \ $\big[\Delta\big] = \sum_{\alpha} a_{\alpha}\otimes b_{\alpha'}.$\par}
\endproclaim

\Proof
Denote by $\delta:X\to X\times_{_{\scs S}} X$, $\delta':X\to X\times X$ 
($\times$ denotes the Cartesian product)  
the diagonal embeddings and by
$\gamma$ the morphism $\pi\times_{_{\scs S}}\pi:X\times_{_{\scs S}} X\to S$.
For $g,h\in A^*(X)$ we have
$$
\pi_*(g\cdot h) = \pi_*\Bigl((\delta')^*(g\times h)\Bigr) =
 \pi_*\Bigl(\delta^*(g\otimes h)\Bigr) = 
\gamma_*\delta_*\Bigl(\delta^*(g\otimes h)\Bigr)
= \gamma_*\Bigl([\Delta]\cdot(g\otimes h)\Bigr),   
$$
where all the equalities follow from the theory in [F1, Chap.~8]
by taking into account, for the second one, the commutative diagram

\vskip\abovedisplayskip
\centerline{\vbox{
\hbox{$X\times_S X \hskip 15pt \hookrightarrow \hskip 15pt X\times X$}
\medskip
\hbox{\vtop {\offinterlineskip
        \hbox{$\hskip 14.7pt\nwarrow$}\vskip-1pt
        \hbox{$\hskip13.6pt{}_\delta\hskip 3.2pt \do$}\vskip-1pt
        \hbox{$\hskip 32.2pt \do$}
            }      
     \hskip 25pt
     \vtop{\offinterlineskip
        \hbox{\hskip 16.9pt $\nearrow$}\vskip-1pt
        \hbox{$\hskip 9.4pt \dup \hskip2.7pt {}_{\delta'}$}\vskip-1pt
        \hbox{$\hskip 1.8pt \dup$}   
          }
      }
\medskip
\hbox{\hskip 50pt $X$}       
                  }
	    }
\vskip\belowdisplayskip

\noindent
and, for the third one, the equality $\pi=\gamma\circ\delta$.
 Hence, writing 
$[\Delta] = \sum d_{\mu\nu}b_{\mu}\otimes a_{\nu}$,
we get 
$$
\aligned
   \pi_*(a_{\alpha}\cdot b_{\beta})
&= \gamma_*\bigl([\Delta]\cdot(a_{\alpha}\otimes b_{\beta})\bigr)=
(\pi_*\otimes\pi_*)\Bigl(\bigl(\sum d_{\mu\nu}b_{\mu}\otimes 
a_{\nu}\bigr)\cdot(a_{\alpha}\otimes b_{\beta})\Bigr) \\
&= \sum_{\mu, \nu}
 d_{\mu\nu}\pi_*(b_{\mu}\cdot a_{\alpha})
\cdot \pi_*(a_{\nu}\cdot b_{\beta}).
\endaligned
\tag *
$$
\medskip
\noindent
(i) By the assumption and (*) with $\alpha$ replaced by $\alpha'$ and $\beta$
--- by $\beta'$, we get
$$
d_{\alpha\beta}=\pi_*(a_{\alpha'} \cdot b_{\beta'})-
\sum_{\mu \ne \alpha, \nu \ne \beta}
 d_{\mu\nu}\pi_*(b_{\mu}\cdot a_{\alpha'})
\cdot \pi_*(a_{\nu}\cdot b_{\beta'}).
\tag**
$$
where the degree of $d_{\mu\nu} \in A^*(S)$ such that
$\pi_*(b_{\mu}\cdot a_{\alpha'})\cdot
\pi_*(a_{\nu}\cdot b_{\beta'})\ne 0$ and
$\mu \ne \alpha$ or $\nu \ne \beta$, 
is smaller than the degree
of $d_{\alpha\beta}$. The assertion now follows by induction 
on the degree of $d_{\alpha\beta}$.

\noindent
(ii) a) $\Rightarrow$ b): By virtue of a), Equation (**) now reads \ 
$\pi_*(a_{\alpha'} \cdot b_{\beta'})= d_{\alpha\beta}$ \  and immediately
implies b).

\noindent
b) $\Rightarrow$ a):
Without loss of generality
we can assume that $\Lambda$ is endowed with a linear ordering
$\prec$ compatible with codimension, i.e. $n_{\alpha_1}<n_{\alpha_2}
\Rightarrow \alpha_1\prec \alpha_2$, $m_{\beta_1}<m_{\beta_2}\Rightarrow
\beta_1 \prec \beta_2$ and such that $\alpha_1\prec\alpha_2\Rightarrow
\alpha_2'\prec \alpha_1'$. The rows and columns of the matrices below
are ordered using the ordering $\prec$.  
Write $x_{\alpha \beta} = \pi_*(a_\alpha \cdot b_\beta)$. 
By virtue of b), Equation (*) gives us the following system of equations:
$$
x_{\alpha\beta} = \sum\nolimits_\mu x_{\alpha \mu} x_{\mu' \beta},
$$
where $\alpha, \beta\in \Lambda$.
Note that the antidiagonal of the matrix 
$M:=[x_{\alpha\beta}]_{\alpha,\beta\in \Lambda}$
is indexed by $\{(\alpha,\alpha')\mid\alpha\in \Lambda\}$. The assumption
implies that this antidiagonal consists of units. Moreover, because
of dimension reasons and the assumption again, we know that 
the entries above the diagonal are zero.
Let $P$ be the permutation matrix corresponding to the bijection 
$\alpha \mapsto \alpha'$ of $\Lambda$. 
The above system of equations is rewritten in the matrix form as: 
$$
MP=MP \cdot MP.
$$
Then $MP$ as a (lower) triangular matrix with the units on the diagonal,
must be the identity matrix. Hence $M=P^{-1}$ and this implies a).
\endproof

\proclaim{Remark 5.3}\rm A standard situation 
when the theorem can be applied is when $\pi: X\to S$
is a locally trivial fibration and $\{a_{\alpha}\}$, $\{b_{\beta}\}$
restrict to bases of the Chow ring of a fibre $F$ which are dual
under the Poincar\'e duality map: $(a,b)\mapsto \int_F a\cdot b$.
In such a situation the above method is successfully applied in [P-R5]
to solve a problem of J.~Harris --- for more on that, see the next section.
\endproclaim

\proclaim{Example 5.4}\rm Let $\pi: \Cal G\to X$ be a relative Grassmannian as in
Proposition 1.3. It is easy to see
the diagonal in $\Cal G_1\times_X \Cal G_2$, where $\Cal G_1
=\Cal G_2=\Cal G$, is given (in the
scheme-theoretic sense) by the vanishing of the entries of a matrix
of the homomorphism $R_{\Cal G_1}\to E_{\Cal G_1}=E_{\Cal G_2}
\to Q_{\Cal G_2}$. Hence, by the theorem and a formula for
the top Chern class of the tensor product [L2], we have that
$\pi_*\bigl(s_I(Q)\cdot s_{\overline J}(R\hak)\bigr)=\delta_{I,J}$ 
for $I,J\subset (r)^q$, 
where $\overline J$ is
the partition whose Ferrers' diagram complements the one of $J^{\sim}$ in 
the rectangle $(q)^r$. 
Equivalently, $\pi_*\bigl(s_I(Q)\cdot s_{(r)^q/J}(-R)\bigr) = \delta_{I,J}$.
This is coherent with Theorem 4.8(i) (invoking Proposition 4.2).
\endproclaim

In a similar way, using the calculation of the class of the diagonal
from [F2, Proposition 7.5] via the top Chern class of 
a suitable vector bundle, 
one can reprove the equality in Theorem 4.8(ii).
On the contrary, it appears that the equality in Theorem 4.8(iii) 
does not admit
a geometric interpretation in terms of the top Chern class
of a vector bundle on a fibre product of two Lagrangian Grassmannian
bundles.

\head 6. Intersection  rings  of spaces $G/P$, divided differences and 
formulas for isotropic degeneracy loci --- an introduction to {\rm [P-R 2-5]}.
\endhead
\vs
{\hfill {\sevenrm ``explicit" is not necessary  EXPLICIT!}}
\vs

In this section we summarize results and techniques mainly from
[P-R 2-5] and [R] as well as from preceding papers [H-B] and [P4, Sect.~6].
This section should serve as a ``friendly" introduction to the papers [P-R 2-5] 
by J. Ratajski and the author, and convince the reader that this series
of very technical --- at the very first glance --- papers relies on a childishly
simple idea! 
\vs

Let $G$ be a semisimple algebraic group, $B\subset G$ -- a Borel subgroup and 
$T\subset B$ a maximal torus. We have the characteristic map
$$
c:S^{\bullet}(X(T))\to A^*(G/B).
$$
(see Section 4) which allows one to study the multiplicative structure of
$A^*(G/B)$.
On the other hand, looking at the additive structure,
one has the Bruhat-Schubert cycles
$X_w=\Bigl[ \ \overline{B\strut^-wB/B} \ \Bigr]\in A^{l(w)}(G/B)$.
Here $w$ runs over the Weyl group $W$ and $B\strut^-$ denotes the Borel
subgroup opposite to $B$.
Since the Bruhat-Schubert cells $B\strut^-wB/B$ $(w\in W)$ form a cellular
decomposition of $G/B$, one has 
$A^*(G/B)=\bigoplus\limits_{w\in W} \Bb Z X_w$.

The key problem is now to understand the coefficients $c_w(f)$ appearing 
in the decomposition
$$
c(f)=\sum\limits_{w\in W} c_w(f) X_w,
$$
where $f\in S^{\bullet}(X(T))$.
This question is answered using the divided-differences operators
$\partial_w$ which were defined in Section~3.

\proclaim{Theorem 6.1 \rm [B-G-G], [D1,2]}If 
$f\in S^{\bullet}(X(T))$ is homogeneous then\break
$c(f)=\sum\limits_{l(w)=\deg f} \partial_w(f) X_w$.
Moreover, the kernel of $c$ is generated by the positive degree 
$W$-invariants and $c\otimes\Bb Q$ is surjective.
\endproclaim

Observe that $\partial_w(\Ker \,c)\subset \Ker \,c$, so $\partial_w$ acts also on
$A^*(G/B)$.
One has (loc.cit.):
$$
\partial_wX_v=\cases
X_{vw^{-1}} &\hbox{ if }\quad l(vw^{-1})=l(v)-l(w)\\
0 & \hbox{otherwise.}\\ \endcases
$$
In particular, $X_w=\partial_{w^{-1}w_0}X_{w_0}$ where $w_0\in W$ is
the longest element in $W$.
One has (loc.cit.):
$$
X_{w_0}=c\left(\prod\limits_{\alpha\in R^+}\alpha/|W|\right)
=c\bigl(\rho^N/N!\bigr),
$$
where $R^+$ is the set of positive roots with respect to $B$ of cardinality 
$N$ and  
$\rho$ denotes half of the sum of positive roots.

A similar theory works for the parabolic subgroups.
Let $\Delta$ be the basis of the root system of $(G,T)$, 
associated with $B$, and let $P\supset B$ be the parabolic 
subgroup associated with a subset $\theta\subset\Delta$.
Denote by $W^{\theta}$ the set
$$
W^{\theta}=\{w\in W \mid l(ws_{\alpha})
=l(w)+1\quad \forall\,\alpha\in\theta\}.
$$
The latter set is the set of minimal length left coset representatives of
$W_{\theta}$, the subgroup of $W$ generated by 
$\{s_{\alpha}\}_{\alpha\in\theta}$, in $W$.
(This follows, e.g., from the following fact:
for every $w\in W$, there exist unique $w^{\theta}\in W^{\theta}$, 
$w_{\theta}\in W_{\theta}$ such that $w=w^{\theta}w_{\theta}$ and
$l(w)=l(w^{\theta})+l(w_{\theta})$ --- see [Bou].)
The projection $G/B\to G/P$ induces an inclusion
$A^*(G/P)\to A^*(G/B)$ which (additively) identifies $A^*(G/P)$ with
$\bigoplus\limits_{w\in W^{\theta}}\Bb ZX_w$.
Multiplicatively, $A^*(G/P)$ is identified with the ring of invariants
${A^*(G/B)}^{W_{\theta}}$.
Indeed, $X_w\in A(G/B)^{W_{\theta}}$ iff $\partial_{\alpha}X_w=0$ for all 
$\alpha\in\theta$.
This takes place iff $l(ws_{\alpha})\ne l(w)-1$, i.e. $l(ws_{\alpha})=l(w)+1$,
for all $\alpha\in\theta$;
or equivalently $w\in W^{\theta}$.
The restriction
$$
c:S^{\bullet}(X(T))^{W_{\theta}}\to A^*(G/P)
$$
of the characteristic map satisfies, for a homogeneous 
$f$,
$$
c(f)=\sum\limits_{\textstyle{w\in W^{\theta} \atop l(w)=\deg f}} 
\partial_w(f)X_w.
$$

The series of papers [P-R 2-5] deals with the case when $G$ is a classical 
group and $P$ is a maximal parabolic subgroup, i.e. 
$\theta=\Delta\setminus\{\alpha\}$ where $\alpha$ is a simple root.
In every case like that, there exists a collection of ``special" Schubert
cycles $\{X_p\}$ forming a minimal set of multiplicative generators of
$A^*(G/P)$ and a collection $\{e_p\}$ of elements in $S^{\bullet}(X(T))$
such that $c(e_p)=X_p$ up to a scalar.
We specify these families ``type by type":
Let $S=S^{\bullet}(X(T))^{W_{\theta}}/\Ker c$.
We use here the ``Bourbaki notation" [Bou] for roots although denoting 
by $A_n$ the sequence of variables
$(a_1,a_2,\dots,a_n)$ where $(a_1,\ldots,a_m)$ is the sequence of ``basic
coordinates" and $n\le m$.
By $s_i$ we denote the simple reflection associated with the $i$-th
simple root $\alpha_i$ specified below (in all types); it should be not
confused with the Schur (or Segre) polynomial. Also, $W_{\alpha}:=W_{\Delta \smallsetminus \{\alpha\}}$.

\vs\noindent
{\bf A}$_{m-1}$ \quad Simple roots:
$\alpha_n=\varepsilon_n-\varepsilon_{n+1},\ n=1,\dots,m-1$.

\noindent
For a fixed $n$, $W_{\alpha_n}=S_n\times S_{m-n}\hookrightarrow S_m=W$;
\ssk
$S=S\Cal P(A_n)/\bigl(s_j(A_n)\mid m-n+1\le j\le m\bigr) ;\ e_p=e_p(A_n)$
\ssk
$X_p=X_{s_{n-p+1} s_{n-p+2}\ldots s_n}$  for $p=1,\dots,n$. 
One has  $c(e_p)=X_p$. 

\noindent
(Here, $G/P=G_n(\Bb C\,\strut^m)$.)

\vs\noindent
{\bf B}$_m$ \quad Simple roots:
$\alpha_n=\varepsilon_n-\varepsilon_{n+1}\ (n<m),\ \alpha_m=\varepsilon_m$.

\noindent
For a fixed $n$, $W_{\alpha_n}=
S_n\times\bigl(S_{m-n}\ltimes\Bb Z\,\strut^{m-n}_2\bigr)\hookrightarrow
S_m\ltimes\Bb Z\,\strut^m_2=W$
($W$ is identified with the group of ``barred permutations" --- see [H-B], 
[P-R2])
\ssk
$S=S\Cal P(A_n)/\bigl(s_j(a_1^2,\dots,a_n^2)\mid m-n+1\le j\le m\bigr);\ 
e_p=e_p(A_n),$
\ssk
$X_p=X_{s_{n-p+1}s_{n-p+2} \dots s_n}$ for $p=1,\dots,n.$
For $m=n$ one has $c(e_p)=2X_p$. 
\ssk
For $n<m$, $c(e_p)=X_p$. 
\ssk\noindent
(Here, $G/P=OG_n(\Bb C\,\strut^{2m+1})$ is the  
Grassmannian of isotropic
$n$-subspaces in $\Bb C\,\strut^{2m+1}$ endowed with a symmetric
nondegenerate form.)

\ssk
\noindent
{\bf C}$_m$ \quad Simple roots:
$\alpha_n=\varepsilon_n-\varepsilon_{n+1}\ (n<m),\ \alpha_m=2\varepsilon_m.$

\noindent
For a fixed $n$, $W_{\alpha_n}=S_n\times
\bigl(S_{m-n}\ltimes\Bb Z\,\strut^{m-n}_2\bigr)\hookrightarrow
S_m\ltimes\Bb Z\,\strut_2^m=W$;
\ssk
$S=S\Cal P(A_n)/\bigl(s_j(a_1^2,\dots,a_n^2)\mid m-n+1\le j\le m
\bigr)\,;\ e_p=e_p(A_n),$
\ssk
$X_p=X_{s_{n-p+1}s_{n-p+2} \dots s_n}$  for $p=1,\dots,n$.  
One has $c(e_p)=X_p$.

\ssk\noindent
(Here, $G/P=LG_n(\Bb C\,\strut^{2m})$ is the Lagrangian Grassmannian of
isotropic $n$-subspaces in $\Bb C\,\strut^{2m}$ endowed with an antisymmetric 
nondegenerate form.)

\vs\noindent
{\bf D}$_m$ \quad Simple roots:
$\alpha_n=\varepsilon_n-\varepsilon_{n+1}\ (n<m),\ \alpha_m=\varepsilon_{m-1}+
\varepsilon_m$.

\noindent
For a fixed $n$, $W_{\alpha_n}=S_n\times
\bigl(S_{m-n}\ltimes\Bb Z\,\strut^{m-n-1}_2\bigr)\hookrightarrow
S_m\ltimes\Bb Z\,\strut_2^{m-1}=W$, the latter being identified with
the group of ``even-barred permutations";

\item{}
For $n=m,\ e_p=e_p(A_n),\ p=1,\dots,n-1$, generate multiplicatively $S$; for
$n<m,$ $e_p(A_n),$ $p=1,\dots,n$ and $a_{n+1}\cdot\dots\cdot a_m$ generate
multiplicatively $S$.

\item{}
We propose the following choice of special Schubert cycles:\break
$X_p=X_{s_{n-p+1}s_{n-p+2} \dots s_n}$  for $p=1,\dots,n$ and 
additionally for $n<m$, $X_0=X_{s_{m-1}s_{m-2}\ldots s_n} - X_{s_ms_{m-2}\ldots s_n}$.
One has $c(e_p)=X_p$ for $p=1,\ldots,n$ and, for $n<m$, 
$c(a_{n+1}\cdot\dots\cdot a_m)=\pm X_0$. 

\noindent
(Here, $G/P=OG_n(\Bb C\,\strut^{2m})$ is the Grassmannian of
isotropic $n$-subspaces in $\Bb C\,\strut^{2m}$ endowed with a symmetric 
nondegenerate form. Recall that for $n<m$, $G/P=OG_n(\Bb C\,\strut^{2m})$
is irreducible, and for $n=m$, it has two connected components.)

\vs
Our main task here are Pieri-type formulas, i.e. we wish to give explicit 
formulas for the coefficients $m^v_{w,p}\in\Bb N$ appearing in the 
equations:
$$
X_w\cdot X_p=\sum m^v_{w,p} X_v
$$
where $w,v\in W^{\theta}$ $(\theta=\Delta\setminus\{\alpha\},\ \alpha$ 
a simple root).
For the spaces $G/B$, there exists a formula of this type due to 
C. Chevalley [Ch]\footnote{Note added in proof: 
This paper, written about 1958,
has been published only recently.}:
$$
X_w\cdot X_{s_{\alpha}}=\sum(\beta\hak,\omega_{\alpha})X_{ws_{\beta}},
$$
where $w\in W$ and $\beta$ runs over positive roots such that
$l(ws_{\beta})=l(w)+1$.
One can use the Chevalley formula %, ``in theory" say, 
to compute the Pieri-type formulas for $G/P$ in question.
This method leads, however, to very complicated computations (compare [H-B])
which make it practically unaccessible for the Grassmannians of
non-maximal subspaces of types $B$, $C$, $D$.

The strategy taken in [P-R 2-4] relies on a totally different method;
namely, it uses an iteration of the following Leibniz-type formula:
for $f,g\in S^{\bullet}(X(T))$ and a root $\alpha$, one has:
$$
\partial_{\alpha}(f\cdot g)=f\cdot(\partial_{\alpha}g)+
(\partial_{\alpha}f)\cdot (s_{\alpha}g).
$$

More precisely, the calculation goes as follows.
Let $f$ be a homogeneous element of $S^{\bullet}(X(T))\otimes {\Bb Q}$ 
such that $c(f)=X_w$.
We have for a reduced decomposition $v=s_{\alpha_1}\cdot\dots\cdot s_{\alpha_k}$,
$$
\aligned
m^v_{w,p} &=\partial_v(f\cdot e_p)
=\partial_{\alpha_1}\circ\dots\circ\partial_{\alpha_k}(f\cdot e_p)=\\
&=\sum \partial_I(f)\cdot\partial^I_{\alpha}(e_p),
\endaligned
$$
where the sum is over all subsequences $I=(i_1<\dots<i_l)\subset\{1,\dots,k\}$;
\ $\partial_I$ is $\partial_{r_I}$ where $r_I=s_{\alpha_{i_1}}\cdot\dots\cdot
s_{\alpha_{i_l}}$ and $\partial^I_{\alpha}$ is obtained by replacing in 
$\partial_{\alpha_1}\circ\dots\circ\partial_{\alpha_k}$ the subword 
$\partial_I$ through the subword $r_I$.
By the choice of $f$ we finally infer

\proclaim{Lemma 6.2}$m^v_{w,p}=\sum\partial^I_{\alpha}(e_p)$, \ 
where the sum runs over $I$ such that $r_I$ is a reduced decomposition of
$w$ $($notation: $r_I\in R(w))$.
\endproclaim

The content of papers [P-R 2-4] and [R] can be summarized, in a coarse way, 
in the following theorem.

\proclaim{Theorem 6.3 \rm [P-R 2-4], [R]}In the above notation, for types
$A_{m-1}$, $B_m$, $C_m$\footnote{Note added in proof:
The authors have recently extended the result also to type $D_m$. See a forthcoming paper:
P. Pragacz and J. Ratajski, {\it A Pieri-type theorem for
even orthogonal Grassmannians}.}
and any simple root $\beta$, the following assertions hold
$(W^{\beta}=W^{\Delta\setminus\{\beta\}})$:

{\parindent=\wd1
\item{\rm I.} For any $v\in W^{\beta}$ there exists a reduced decomposition 
$v=s_{\alpha_1}\cdot\dots\cdot s_{\alpha_k}$ such that for any 
$w\in W^{\beta}$ with $l(w)=k-p$, the cardinality of the set
$$
\Cal I=\bigl\{I\subset\{1,\dots,k\}\mid r_I\in R(w) \hbox{ and } \partial^I_{\alpha}(e_p)\ne 0\bigr\}
$$
is less than or equal to $1$.
\item{\rm II.} There exists an EXPLICIT combinatorial criterion for 
$\card \Cal I=1$.
\item{\rm III.} If $\card\Cal I=1$ then the unique $I=I(\alpha.,w)\in \Cal I$
is given by an EXPLICIT combinatorial algorithm.
\item{\rm IV.} The multiplicity is given by the formula
$m^v_{w,p}=\partial^I_{\alpha}(e_p)=2^{e(w,v)}$, where $e(w,v)$ is given
by an EXPLICIT combinatorial rule.
\par}
\endproclaim

We now illustrate the theorem for type $C_n$ and the simple root $\alpha_n$.
Let, for the rest of this section, $\Cal G=LG_n(\Bb C^{2n})$ denote 
the Lagrangian Grassmannian of $n$-dimensional isotropic subspaces in 
$\Bb C^{2n}$ with respect to a nondegenerate antisymmetric form 
on $\Bb C^{2n}$. 
Let $\Cal F$ denote the flag variety of (total) isotropic flags in $\Bb C^{2n}$ 
(with respect to the same antisymmetric form). 
By $\rho_n$ we denote the partition $(n,\dots,2,1)$. 
Let $I=(i_1>i_2>\dots>i_k>0) \subset \rho_n$ be a strict partition. 
We associate with $I$ the 
element $w_I$  of the symplectic Weyl group $W$\footnote{%
Footnote 9 applies here with $S_n$  replaced by $W$.}:
$$
w_I= (s_{n-i_k+1}\dots s_{n-1} s_n) \ \dots \ (s_{n-i_2+1}\dots s_{n-1} s_n) 
(s_{n-i_1+1}\dots s_{n-1} s_n).
$$

From the theory described above we get the {\it Schubert cycle} 
$X_{w_I}\in A^{|I|}(\Cal F)$, defined as the class of the closure of $B^-w_IB/B$  
where $B$ is the Borel subgroup of the symplectic group and 
$B^-$ --- its opposite. 
In fact, $X_{w_I}$ belongs to $A^{|I|}(\Cal G)\subset A^{|I|}(\Cal F)$. 
Denote this element in $A^{|I|}(\Cal G)$ by $\sigma(I)$, for brevity. 

As usual, we associate to a partition $I$ its Ferrers' diagram 
$D_I$ (see [M1]) and treat it as a subset of $\Bb Z \times \Bb Z$. 
A subset $D \subset \Bb Z\times \Bb Z$ is {\it connected} if
each of the sets $\{i: \exists j \ (i,j)\in D \}$ and $\{j: \exists
i \ (i,j)\in D \}$ is an interval in $\Bb Z$. 
This allows us to speak about the ``{\it connected components}" 
of skew diagrams, i.e. the differences between 
diagrams of partitions. 

The following result was given originally in [H-B] and reproved in [P-R2].

\proclaim{Theorem 6.4 \rm [H-B]}Let $I=(i_1,\dots,i_k)\subset\rho_n$ 
be a strict partition. The following equality holds 
in $A^*(\Cal G)$ \ $(p=1,\dots,n):$
$$
\sigma(I)\ \sigma(p) = \sum 2^{e(I,J)}\ \sigma(J),
$$
where the sum is over strict partitions $J$ such that $i_{h-1}\ge j_h\ge i_h$  
$(i_0=n, \ i_{k+1}=0)$, $|J|=|I|+p$ and $e(I,J)$ is the number of connected 
components of $D_J\setminus D_I$ not meeting the first column.
\endproclaim

\noindent
(This formulation of the theorem is slightly different from the one in
[H-B] and is suited to the techniques of [P-R2].)

\proclaim{Example 6.5}\rm $n=7$
$$\sigma(632)\ \sigma(5) = 2\sigma(763)+2^2\sigma(7531)+2\sigma(7621)+ 
2\sigma(7432)+\sigma(6532).$$
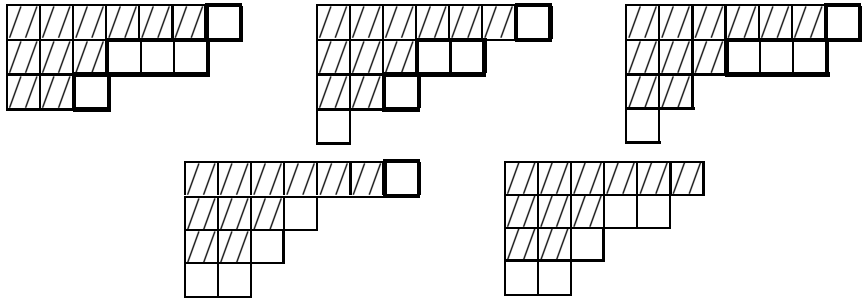
\centerline{%
% Diagram 1
\vtop{\offinterlineskip
    \hbox{\vbox{\hrule width 62.8pt} $\!\!\!\!$
          \vbox{\hrule width 11.3pt height 0.8pt depth 0.2pt} }   
    \hbox{\vrule //\vrule//\vrule//\vrule//\vrule//\vrule//\vrule 
             width 1pt \ \ \ \vrule width 1pt}
    \hbox{\vbox{\hrule width 31,6pt} $\!\!\!\!$
          \vbox{\hrule width 42.7pt height 0.8pt depth 0.2pt}  }
    \hbox{\vrule//\vrule//\vrule//\vrule width 1pt \quad\vrule\quad\vrule \quad
          \vrule width 1pt}
    \hbox{\vbox{\hrule width 21pt} $\!\!\!\!$
       \vbox{\hrule width 42.8pt height 0.7pt depth 0.3pt}  }
    \hbox{\vrule //\vrule//\vrule width 1pt \ \ \ \vrule width 1 pt}
    \hbox{\vbox{\hrule width 21pt} $\!\!\!\!$
       \vbox{\hrule width 11.7pt height 0.7pt depth 0.3pt}  }
     }\qquad
% Diagram 2     
\vtop{\offinterlineskip
    \hbox{\vbox{\hrule width 62.9pt} $\!\!\!\!$
       \vbox{\hrule width 11.3pt height 0.8pt depth 0.2pt} }   
    \hbox{\vrule//\vrule//\vrule//\vrule//\vrule//\vrule//\vrule 
                                      width 1pt \quad \vrule width 1pt}
    \hbox{\vbox{\hrule width 31.7pt} $\!\!\!\!$
       \vbox{\hrule width 22pt height 0.8pt depth 0.2pt} $\!\!\!\!$
       \vbox{\hrule width 9.1pt} $\!\!\!\!$
       \vbox{\hrule width 11.4pt height 0.8pt depth 0.2pt}     }
    \hbox{\vrule//\vrule//\vrule//\vrule width 1pt \quad\vrule\quad 
                                                       \vrule width 1pt}
    \hbox{\vbox{\hrule width 21.1pt} $\!\!\!\!$
       \vbox{\hrule width 32.2pt height 0.7pt depth 0.3pt}  }
    \hbox{\vrule //\vrule//\vrule width 1pt \quad \vrule width 1 pt}
    \hbox{\vbox{\hrule width 21.1pt} $\!\!\!\!$
       \vbox{\hrule width 11.5pt height 0.7pt depth 0.3pt}  }
    \hbox{\vrule height 7pt depth 3pt \quad \vrule}
    \hrule width 10.9pt
     }\qquad
% Diagram 3
\vtop{\offinterlineskip
    \hbox{\vbox{\hrule width 62.9pt} $\!\!\!\!$
       \vbox{\hrule width 11.3pt height 0.8pt depth 0.2pt} }   
    \hbox{\vrule //\vrule//\vrule//\vrule//\vrule//\vrule//\vrule 
                                      width 1pt \quad \vrule width 1pt}
    \hbox{\vbox{\hrule width 31.6pt} $\!\!\!\!$
          \vbox{\hrule width 42.7pt height 0.8pt depth 0.2pt}  }
    \hbox{\vrule//\vrule//\vrule//\vrule width 1pt \quad\vrule\quad\vrule\quad 
          \vrule width 1pt}
    \hbox{\vbox{\hrule width 31.6pt} $\!\!\!\!$
       \vbox{\hrule width 32.5pt height 0.7pt depth 0.3pt}  }
    \hbox{\vrule//\vrule//\vrule }
    \hbox{\vbox{\hrule width 21.5pt}  }
    \hbox{\vrule height 7.3pt depth 3pt \quad \vrule}
    \hrule width 11pt 
     }}
\vs
\centerline{
%Diagram 4
\vtop{\offinterlineskip
    \hbox{\vbox{\hrule width 62.9pt} $\!\!\!\!$
       \vbox{\hrule width 11.3pt height 0.8pt depth 0.2pt} }   
    \hbox{\vrule//\vrule//\vrule//\vrule//\vrule//\vrule//\vrule 
           width 1pt \ \ \ \vrule width 1pt}
    \hbox{\vbox{\hrule width 62.9pt} $\!\!\!\!$
       \vbox{\hrule width 11.3pt height 0.8pt depth 0.2pt}  }
    \hbox{\vrule//\vrule//\vrule//\vrule\quad\vrule }
    \hrule width 42pt
    \hbox{\vrule//\vrule//\vrule \quad \vrule}
    \hrule width 31.5pt
    \hbox{\vrule height 7.3pt depth 3pt \quad\vrule\quad \vrule}
    \hrule width 21pt 
     }     \qquad
%Diagram 5
\vtop{\offinterlineskip
    \hrule width 63.2pt
    \hbox{\vrule//\vrule//\vrule//\vrule//\vrule//\vrule//\vrule }
    \hrule width 63.2pt
    \hbox{\vrule//\vrule//\vrule//\vrule\quad\vrule\quad\vrule }
    \hrule width 52.4pt
    \hbox{\vrule//\vrule//\vrule \quad \vrule}
    \hrule width 31.5pt
    \hbox{\vrule height 7.3pt depth 3pt \quad\vrule\quad \vrule}
    \hrule width 21pt 
     }}\endproclaim

We now sketch a proof of this theorem, in the spirit of Theorem 6.3
above, stemming from [P-R2]. 
This proof is much simpler than the one in [H-B] based on the 
Chevalley formula.
\vs

Let $A=(a_1,\dots,a_n)$ be a sequence of independent variables. 
It follows from \hbox{[B-G-G]}
and [D1,2] that $A^*(\Cal F)$ is identified with $\Bb Z[A]/\Cal I$, where 
$\Cal I$ is the ideal generated by symmetric polynomials in $a^2_1,\dots,a^2_n$ 
without constant term. 
Also, $A^*(\Cal G)$ is identified with $(\Bb Z[A]/\Cal I)^{S_n}$, i.e. with 
the quotient of the symmetric polynomials modulo $\Cal I$ restricted to 
the ring of symmetric polynomials.

We have ``symplectic divided differences":
$\partial_i:\Bb Z[A]\to\Bb Z[A]$ \ (of degree $-1$), $i=1,\dots,n$,
\quad defined by
$$
\partial_i(f)=(f - s_if)/(a_i - a_{i+1}), \quad  i=1,\dots,n-1, \qquad
\partial_n(f) = (f-s_nf)/2a_n.
$$
 
Recall that there exists a surjective characteristic map 
$c:\Bb Z[A]\to A^*(\Cal F)$, and $c(e_p)=\sigma(p)=
X_p\in A^p(\Cal G)$ ([H-B]).

Let $f_I\in\Bb Z[A]$ be homogeneous such that $c(f_I) = \sigma(I)$. 
Then, for $w\in W$, $l(w)=|\lambda|$, we have $\partial_w(f_I)\not=0$ 
iff $w=w_I$ and $\partial_{w_I}(f_I)=1$. 
Following Theorem 6.3, our goal is to find the coefficients $m_J$  appearing in
$$
c(f_I\cdot e_p) = \sum m_J\ \sigma(J).
$$

Consider an arbitrary subset $D\subset D_J$. (Subsequences $I$ in Lemma 6.2
and Theorem 6.3 correspond here to subsets $D$ of $D_J$.) 
The boxes in $D_J$ which belong to $D$ will be called $D${\it -boxes}; 
the boxes in $D_J$$\setminus$$D$  will be called {\it non-$D$-boxes}. 
We associate with $D$ the %following 
word $r_D$ and the operator $\partial^D_J$,
given by Definitions 6.6 and 6.7. 
For technical reasons we will use, from now on, the ``matrix" coordinates 
for indexing boxes in $D_J\subset\rho_n$ but the columns are numbered from
left to right by $n,n-1,\ldots, 2,1$ successively.

In Definitions 6.6, 6.7  we read $D_J$ like books (e.g., in Europe), i.e.
row by row from left to right starting
from the top row.

\proclaim{Definition 6.6 \rm of $r_D$}\  \rm Read $D_J$. 
Every $D$-box in the $i$-th column gives us the $s_i$. 
Non-$D$-boxes have no influence on $r_D$. 
Then $r_D$ is the word obtained by writing the $s_i$'s from 
right to left.
\endproclaim

\proclaim{Definition 6.7 \rm of $\partial^D_J$}\ \rm
Read $D_J$. Every $D$-box in the $i$-th column gives us the $s_i$. 
Every non-$D$-box in the $i$-th column gives the $\partial_i$. 
Then $\partial^D_J$ is the composition of the  
$s_i$'s and $\partial_i$'s  (the composition written from right to left).
\endproclaim

\proclaim{Example 6.8}\rm $J=(763),\ n=7$.
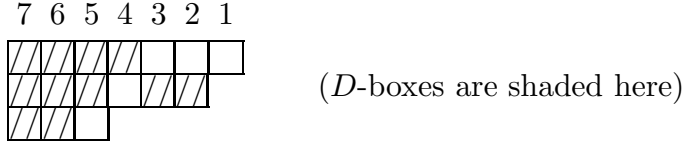
\vs

\hbox{
\qquad\qquad
\vtop{\hbox{\;7\hskip 5.6pt 6\hskip 5.6pt 5\hskip 5.6pt 4\hskip 5.6pt 
              3\hskip 5.6pt 2\hskip 5.6pt 1} 
    \vskip 5pt
    \hrule width 74.5pt
    \offinterlineskip
    \hbox{\vrule //\vrule//\vrule//\vrule//\vrule\hskip 10.4pt
          \vrule\hskip 10.4pt\vrule\hskip 10.4pt\vrule} 
    \hrule width 74.5pt
    \hbox{\vrule//\vrule//\vrule//\vrule\hskip 10.4pt\vrule//\vrule//\vrule}
    \hrule width 63.6pt 
    \hbox{\vrule //\vrule//\vrule \hskip 10.4pt \vrule}
    \hrule width 32pt
     }  
\qquad
\vtop{\vskip 15pt
    \hbox{($D$-boxes are shaded here)}
     }}
\vs
$r_D = s_6\cdot s_7\cdot s_2\cdot s_3\cdot s_5\cdot s_6\cdot 
s_7\cdot s_4\cdot s_5\cdot s_6\cdot s_7$ ,
\vs
$\partial^D_J=\partial_5\circ s_6\circ s_7\circ s_2\circ s_3\circ
\partial_4\circ s_5\circ s_6\circ s_7\circ\partial_1\circ\partial_2\circ
\partial_3\circ s_4\circ s_5\circ s_6\circ s_7$.
\endproclaim

In the above notation, Lemma 6.2 reads as follows:
$
m_J=\sum \partial^D_J(e_p),
$
where the sum is over all $D\subset D_J$ such that  $r_D\in R(w_I)$.

One proves that if $j_{h+1}>i_h$ for some $h$ (in particular if $l(J)>l(I)+1)$ 
then $\partial^D_J(\sum e_p)=0$  for every  $D\subset D_J$  such that  
$r_D\in R(w_I)$.

Moreover, fix a strict partition  $I\subset \rho_n$. 
Let $J$ be a strict partition such that  $I\subset J\subset\rho_n$, 
$j_{h+1}\le i_h$  for every $h$ (in particular $l(J)\le l(I)+1$). 
Then
there exists {\it exactly one} $D^{I,J}\subset D_J$  such that  
$r_D\in R(w_I)$  and  $\partial^D_J(\sum e_p)\not=0$ 
for $D=D^{I,J}$.

The idea of constructing such a $D^{I,J}$  can be easily explained using 
pictures. 
The boxes from $D_I\subset D_J$  are shaded in the pictures below. 
A part:
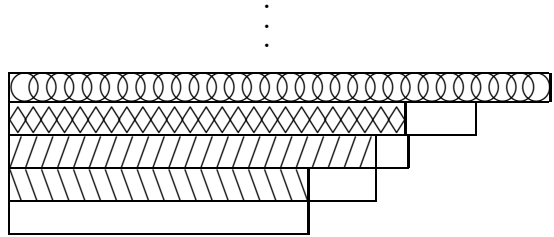
\vs

%  str. 15 diagram 1
\centerline{
\vbox{\offinterlineskip
      \hbox{\hskip 80pt .}\vskip 5pt
      \hbox{\hskip 80pt .}\vskip 5pt
      \hbox{\hskip 80pt .}\vskip 8pt
      \hrule width 170pt
      \hbox{\vrule $\!\hskip 1.4pt\bi\!\!\!\!\bi\!\!\!\!\bi
            \!\!\!\!\bi\!\!\!\!\bi\!\!\!\!\bi\!\!\!\!\bi\!\!\!\!\bi
	    \!\!\!\!\bi\!\!\!\!\bi\!\!\!\!\bi\!\!\!\!\bi\!\!\!\!\bi
	    \!\!\!\!\bi\!\!\!\!\bi\!\!\!\!\bi\!\!\!\!\bi\!\!\!\!\bi
	    \!\!\!\!\bi\!\!\!\!\bi\!\!\!\!\bi\!\!\!\!\bi\!\!\!\!\bi
	    \!\!\!\!\bi\!\!\!\!\bi\!\!\!\!\bi\!\!\!\!\bi\!\!\!\!\bi
	    \!\!\!\!\bi\!\!\!\bi\!\hskip 1.1pt$\vrule}
      \hrule width 170pt
      \hbox{\vrule height 8.3pt depth 1.7pt\hskip 0.9pt$\!\lo\!\lo\!\lo\!\lo
            \!\lo\!\lo\!\lo\!\lo\!\lo\!\lo\!\lo\!\lo\!\lo\!\lo\!\lo\!\lo\!\lo
	    \!\lo\!\lo\!\lo\!\lo\!\lo\!\lo\!\lo\!\lo\!$\hskip 
	    0.9pt\vrule\hskip 21.2pt\vrule}
      \hrule width 147.4pt
      \hbox{\vrule $\!\hskip 1pt$///////////////////////\hskip 0.4pt\vrule\hskip 10.4pt\vrule}
      \hrule width 126.1pt
      \hbox{\vrule $\!\hskip 1pt\ba\ba\ba\ba\ba\ba\ba\ba\ba\ba\ba\ba\ba\ba\ba\ba\ba
                    \ba\ba\!$\hskip 1pt\vrule\hskip 21.2pt\vrule}
      \hrule width 115.8pt
      \hbox{\vrule height 7pt depth 3pt \hskip 94pt\vrule}
      \hrule width 94.8pt
     } }

\vs\noindent
of the diagram $D_I\subset D_J$ is deformed to: 
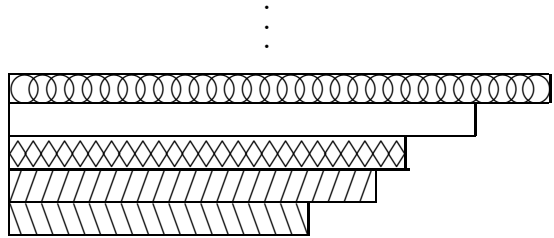
\vs

%  str. 15 diagram 2
\centerline{
\vbox{\offinterlineskip
      \hbox{\hskip 80pt .}\vskip 5pt
      \hbox{\hskip 80pt .}\vskip 5pt
      \hbox{\hskip 80pt .}\vskip 8pt
      \hrule width 170pt
      \hbox{\vrule $\!\hskip 1.4pt\bi\!\!\!\!\bi\!\!\!\!\bi\!\!\!\!\bi
            \!\!\!\!\bi\!\!\!\!\bi\!\!\!\!\bi\!\!\!\!\bi\!\!\!\!\bi\!\!\!\!\bi
	    \!\!\!\!\bi\!\!\!\!\bi\!\!\!\!\bi\!\!\!\!\bi\!\!\!\!\bi\!\!\!\!\bi
	    \!\!\!\!\bi\!\!\!\!\bi\!\!\!\!\bi\!\!\!\!\bi\!\!\!\!\bi\!\!\!\!\bi
	    \!\!\!\!\bi\!\!\!\!\bi\!\!\!\!\bi\!\!\!\!\bi\!\!\!\!\bi\!\!\!\!\bi
	    \!\!\!\!\bi\!\!\!\bi\!\hskip 1.1pt$\vrule}
      \hrule width 170pt
      \hbox{\vrule height 8.3pt depth 1.7pt\hskip 146.6pt\vrule}
      \hrule width 147.2pt
      \hbox{\vrule height 8.3pt depth 1.7pt\hskip 0.9pt$\!\lo\!\lo\!\lo\!\lo
            \!\lo\!\lo\!\lo\!\lo\!\lo\!\lo\!\lo\!\lo\!\lo\!\lo\!\lo\!\lo\!\lo
	    \!\lo\!\lo\!\lo\!\lo\!\lo\!\lo\!\lo\!\lo\!$\hskip 0.9pt\vrule}      
      \hrule width 126.3pt
      \hbox{\vrule $\!\hskip 1pt$///////////////////////\hskip 0.4pt\vrule}      
      \hrule width 115.9pt      
      \hbox{\vrule $\!\hskip 1pt\ba\ba\ba\ba\ba\ba\ba\ba\ba\ba\ba\ba\ba\ba\ba\ba\ba
                    \ba\ba\!$\hskip 1pt\vrule}
      \hrule width 94.8pt
     } }     

\vs\noindent
On the other side, a part: 
\bigskip

%  str. 15 diagram 3
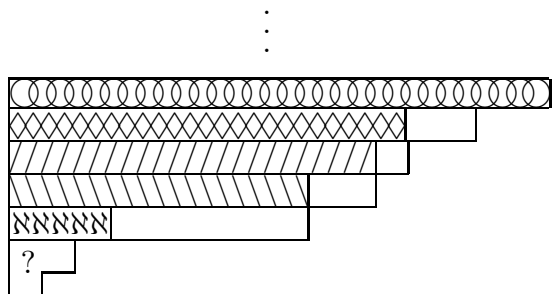
\centerline{
\vbox{\offinterlineskip
      \hbox{\hskip 80pt .}\vskip 5pt
      \hbox{\hskip 80pt .}\vskip 5pt
      \hbox{\hskip 80pt .}\vskip 8pt
      \hrule width 170pt
      \hbox{\vrule $\!\hskip 1.4pt\bi\!\!\!\!\bi\!\!\!\!\bi
            \!\!\!\!\bi\!\!\!\!\bi\!\!\!\!\bi\!\!\!\!\bi\!\!\!\!\bi
	    \!\!\!\!\bi\!\!\!\!\bi\!\!\!\!\bi\!\!\!\!\bi\!\!\!\!\bi
	    \!\!\!\!\bi\!\!\!\!\bi\!\!\!\!\bi\!\!\!\!\bi\!\!\!\!\bi
	    \!\!\!\!\bi\!\!\!\!\bi\!\!\!\!\bi\!\!\!\!\bi\!\!\!\!\bi
	    \!\!\!\!\bi\!\!\!\!\bi\!\!\!\!\bi\!\!\!\!\bi\!\!\!\!\bi
	    \!\!\!\!\bi\!\!\!\bi\!\hskip 1.1pt$\vrule}
      \hrule width 170pt
      \hbox{\vrule height 8.3pt depth 1.7pt\hskip 0.9pt$\!\lo\!\lo\!\lo\!\lo
            \!\lo\!\lo\!\lo\!\lo\!\lo\!\lo\!\lo\!\lo\!\lo\!\lo\!\lo\!\lo\!\lo
	    \!\lo\!\lo\!\lo\!\lo\!\lo\!\lo\!\lo\!\lo\!$\hskip 
	    0.9pt\vrule\hskip 21.2pt\vrule}
      \hrule width 147.4pt
      \hbox{\vrule $\!\hskip 1pt$///////////////////////\hskip 0.4pt\vrule\hskip 10.4pt\vrule}
      \hrule width 126.1pt
      \hbox{\vrule $\!\hskip 1pt\ba\ba\ba\ba\ba\ba\ba\ba\ba\ba\ba\ba\ba\ba\ba\ba\ba
                    \ba\ba\!$\hskip 1pt\vrule\hskip 21.2pt\vrule}
      \hrule width 115.8pt
      \hbox{\vrule height 8.3pt depth 1.7pt\hskip 0.6pt$\aleph\aleph\aleph\aleph
            \aleph$\hskip 0.4pt\vrule\hskip 62pt\vrule}
      \hrule width 94.8pt
      \hbox{\vrule height 10.2pt depth 10pt\hskip 3.4pt?\hskip 1.8pt\vrule 
             height 0pt depth 10pt$\!$\hskip 1.2pt
	     \vbox{\hrule width 10.4pt}\vrule height 10.2pt depth 0pt}
      \hrule width 10.6pt	     
     } }     

\vs\noindent
of the diagram $D_I\subset D_J$  is deformed to: 
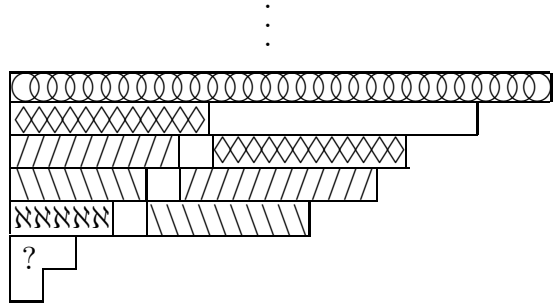
\vs

%  str. 15 diagram 4
\centerline{
\vbox{\offinterlineskip
      \hbox{\hskip 80pt .}\vskip 5pt
      \hbox{\hskip 80pt .}\vskip 5pt
      \hbox{\hskip 80pt .}\vskip 8pt
      \hrule width 170pt
      \hbox{\vrule $\!\hskip 1.4pt\bi\!\!\!\!\bi\!\!\!\!\bi
            \!\!\!\!\bi\!\!\!\!\bi\!\!\!\!\bi\!\!\!\!\bi\!\!\!\!\bi
	    \!\!\!\!\bi\!\!\!\!\bi\!\!\!\!\bi\!\!\!\!\bi\!\!\!\!\bi
	    \!\!\!\!\bi\!\!\!\!\bi\!\!\!\!\bi\!\!\!\!\bi\!\!\!\!\bi
	    \!\!\!\!\bi\!\!\!\!\bi\!\!\!\!\bi\!\!\!\!\bi\!\!\!\!\bi
	    \!\!\!\!\bi\!\!\!\!\bi\!\!\!\!\bi\!\!\!\!\bi\!\!\!\!\bi
	    \!\!\!\!\bi\!\!\!\bi\!\hskip 1.1pt$\vrule}
      \hrule width 170pt
      \hbox{\vrule height 8.3pt depth 1.7pt\hskip 0.7pt$\lo\!\lo\!\lo\!\lo
            \!\lo\!\lo\!\lo\!\lo\!\lo\!\lo\!\lo\!\lo$\hskip 
	    0.6pt\vrule\hskip 83.5pt\vrule}
      \hrule width 147.4pt
      \hbox{\vrule\hskip 1.3pt//////////\hskip 1.4pt\vrule\hskip 
            10.4pt\vrule\hskip 1.3pt$\!\lo\!\lo\!\lo\!\lo\!\lo\!\lo\!\lo\!\lo
	    \!\lo\!\lo\!\lo\!\lo\!$\hskip 1.4pt\vrule}
      \hrule width 126.1pt
      \hbox{\vrule \hskip 1.2pt$\ba\ba\ba\ba\ba\ba\ba\ba$\hskip 1.2pt\vrule
            \hskip 10.4pt\vrule\hskip 0.7pt////////////\hskip 0.6pt\vrule}
      \hrule width 115.8pt
      \hbox{\vrule height 8.3pt depth 1.7pt\hskip 0.8pt$\aleph\aleph\aleph
            \aleph\aleph$\hskip 0.6pt\vrule\hskip 10.4pt\vrule\hskip 
	    0.3pt$\ba\ba\ba\ba\ba\ba\ba\ba\ba\ba$\hskip 0.3pt\vrule}
      \hrule width 94.8pt
      \hbox{\vrule height 10.2pt depth 10pt\hskip 3.4pt?\hskip 1.8pt\vrule 
             height 0pt depth 10pt$\!$\hskip 1.2pt
	     \vbox{\hrule width 10.4pt}\vrule height 10.2pt depth 0pt}
      \hrule width 10.6pt	     
     } }

\vs\noindent
The deformations are performed in direction South $\to$ North.

Fix a strict partition $I\subset\rho_n$ and a number $p=1,\dots,n$. 
Let $J$ be a strict partition such that  $I\subset J\subset\rho_n$, 
$|J|=|I|+p$, $j_{h+1}\le i_h$  for every $h$. 
Let  $D=D^{I,J}$. 
Every $\partial_i$ involved in $\partial^D_J$ is associated with 
a box in $D_J\ba D$. 
It turns out that the connected components of $D_J\ba D$ play a crucial role 
in the computation of $\partial_J(e_p)$. 
Namely, in the above notation, 
$$
\partial^D_J(e_p) = 2^{e(I,J)},
$$
where $e(I,J)$ is the number of connected components of $D_J\ba D$ 
not meeting the $n$-th column. 
By changing the numbering of columns to the usual ordering, this 
can be easily restated as: $e(I,J)$ is the number of connected components of 
$D_J\ba D_I$ not meeting the first column. 

This finishes the sketch of the proof of Theorem 6.4 in the spirit of
Theorem 6.3 --- for more details see [P-R2].

\proclaim{Example 6.9}\rm The diagrams  $D^{(632),J}$  for partitions $J$ 
appearing in the decomposition  $\sigma(632)\ \sigma(5)$, are:
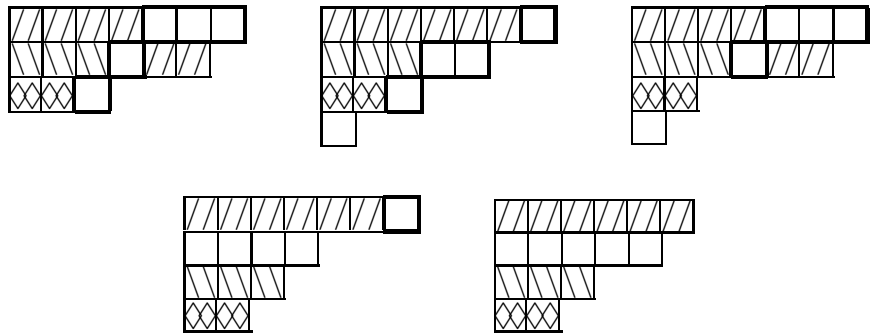
\vs

\centerline{
\vtop{\offinterlineskip
    \hbox{\vbox{\hrule width 43.8pt}$\!$\vbox{\hrule width 32.7pt
                height 0.8pt depth 0.2pt}}
    \hbox{\vrule//\vrule//\vrule//\vrule//\vrule width 1pt\hskip 
          10.1pt\vrule\hskip 10.4pt\vrule\hskip 10.1pt\vrule width 1pt}
    \hbox{\vbox{\hrule width 33pt}$\!$\vbox{\hrule width 43.6pt
                height 0.8pt depth 0.2pt}}
    \hbox{\vrule$\ba\ba$\vrule$\ba\ba$\vrule$\ba\ba\!$\hskip1.4pt\vrule 
          width 1pt\hskip 10.1pt\vrule width 1pt$\!$\hskip
	  1.4pt//$\!$\hskip 1.4pt\vrule//\hskip 0.2pt\vrule}
    \hbox{\vbox{\hrule width 22.5pt}$\!$\vbox{\hrule width 22.3pt
                height 0.8pt depth 0.2pt}$\!$\vbox{\hrule width 22.5pt}}
    \hbox{\vrule height 8.3pt depth 1.7pt$\!\hskip0.9pt\lo\!\!\hskip1.2pt\lo
          \!\hskip0.9pt$\vrule$\!\hskip0.9pt\lo\!\!\hskip1.3pt\lo
	  \!\hskip0.9pt$\vrule width 1pt\hskip 10.1pt\vrule width1pt}
    \hbox{\vbox{\hrule width 22.5pt}$\!$\vbox{\hrule width 11.3pt
                height 0.8pt depth 0.2pt}}
      }    \qquad
% str.16 diagram 2
\vtop{\offinterlineskip
    \hbox{\vbox{\hrule width 64.7pt}$\!$\vbox{\hrule width 11.5pt
                height 0.8pt depth 0.2pt}}
    \hbox{\vrule//\vrule//\vrule//\vrule//\vrule//\vrule//\vrule 
          width 1pt\hskip 10.1pt\vrule width 1pt}
    \hbox{\vbox{\hrule width 33.05pt}$\!$\vbox{\hrule width 22pt
                height 0.8pt depth 0.2pt}$\!$\vbox{\hrule width 
		12.9pt}$\!$\vbox{\hrule width 11.5pt height 0.8pt depth 0.2pt}}
    \hbox{\vrule$\ba\ba$\vrule$\ba\ba$\vrule$\ba\ba\!$\hskip1.4pt\vrule 
          width 1pt\hskip 10.1pt\vrule\hskip 10.1pt\vrule width 1pt}
    \hbox{\vbox{\hrule width 22.4pt}$\!$\vbox{\hrule width 32.7pt
                height 0.8pt depth 0.2pt}}
    \hbox{\vrule height 8.3pt depth 1.7pt$\!\hskip0.9pt\lo\!\!\hskip1.2pt\lo
          \!\hskip0.9pt$\vrule$\!\hskip0.9pt\lo\!\!\hskip1.3pt\lo
	  \!\hskip0.9pt$\vrule width 1pt\hskip 10.1pt\vrule width1pt}
    \hbox{\vbox{\hrule width 22.5pt}$\!$\vbox{\hrule width 11.3pt
                height 0.8pt depth 0.2pt}}
    \hbox{\vrule height 7pt depth3pt\hskip 10.4pt\vrule}
    \hrule width 11.2pt	
      }      \qquad
% str.16 diagram 3
\vtop{\offinterlineskip
    \hbox{\vbox{\hrule width 43.8pt}$\!$\vbox{\hrule width 32.7pt
                height 0.8pt depth 0.2pt}}
    \hbox{\vrule//\vrule//\vrule//\vrule//\vrule width 1pt\hskip 
          10.1pt\vrule\hskip 10.4pt\vrule\hskip 10.1pt\vrule width 1pt}
    \hbox{\vbox{\hrule width 33pt}$\!$\vbox{\hrule width 43.6pt
                height 0.8pt depth 0.2pt}}
    \hbox{\vrule$\ba\ba$\vrule$\ba\ba$\vrule$\ba\ba\!$\hskip1.4pt\vrule 
          width 1pt\hskip 10.3pt\vrule width 1pt$\!$\hskip
	  1.4pt//$\!$\hskip 1.4pt\vrule//\hskip 0.2pt\vrule}
    \hbox{\vbox{\hrule width 33.3pt}$\!$\vbox{\hrule width 11.4pt
                height 0.8pt depth 0.2pt}$\!$\vbox{\hrule width 22.5pt}}
    \hbox{\vrule height 8.3pt depth 1.7pt$\!\hskip1pt\lo\!\!\hskip1.3pt\lo
          \!\hskip1pt$\vrule$\!\hskip1pt\lo\!\!\hskip1.4pt\lo
	  \!\hskip1pt$\vrule }
    \hrule width 21.5pt
    \hbox{\vrule height 7pt depth3pt\hskip 10.4pt\vrule}
    \hrule width 11.2pt	    
      }
    }
\bigskip    
\centerline{
% str.16 diagram 4
\vbox{\offinterlineskip
    \hbox{\vbox{\hrule width 64.6pt}$\!$\vbox{\hrule width 11.55pt
                height 0.8pt depth 0.2pt}}
    \hbox{\vrule//\vrule//\vrule//\vrule//\vrule//\vrule//\vrule 
          width 1pt\hskip 10.1pt\vrule width 1pt}
    \hbox{\vbox{\hrule width 64.6pt}$\!$\vbox{\hrule width 11.55pt
                height 0.8pt depth 0.2pt}}
    \hbox{\vrule height 7pt depth 3pt\hskip 10.15pt\vrule\hskip 
           10.15pt\vrule\hskip 10.15pt\vrule\hskip 10.15pt\vrule}
    \hrule width 42.6pt 	   
    \hbox{\vrule$\ba\ba$\vrule$\ba\ba$\vrule$\ba\ba$\vrule}
    \hrule width 32.1pt
    \hbox{\vrule height 8.3pt depth 1.7pt$\!\hskip0.9pt\lo\!\!\hskip1.2pt\lo
          \!\hskip0.9pt$\vrule$\!\hskip0.9pt\lo\!\!\hskip1.5pt\lo
	  \!\hskip1pt$\vrule}    
    \hrule width 21.5pt
      }    \qquad 
% str.16 diagram 5
\vbox{\offinterlineskip
    \hrule width 63.1pt
    \hbox{\vrule//\vrule//\vrule//\vrule//\vrule//\vrule//\vrule} 
    \hrule width 63.1pt
    \hbox{\vrule height 7pt depth 3pt\hskip 10.15pt\vrule\hskip 
           10.15pt\vrule\hskip 10.15pt\vrule\hskip 10.15pt\vrule\hskip 
	   10.15pt\vrule}
    \hrule width 53.1pt 	   
    \hbox{\vrule$\ba\ba$\vrule$\ba\ba$\vrule$\ba\ba$\vrule}
    \hrule width 32.1pt
    \hbox{\vrule height 8.3pt depth 1.7pt$\!\hskip0.9pt\lo\!\!\hskip1.2pt\lo
          \!\hskip0.9pt$\vrule$\!\hskip1pt\lo\!\!\hskip1.5pt\lo
	  \!\hskip1pt$\vrule}    
    \hrule width 21.5pt
      }      
     }\endproclaim

A geometric interpretation of  $\sigma(I)$ is as follows. 
Let $V$ be a vector space of dimension $2n$ endowed with an  antisymmetric 
nondegenerate form  $\phi:V\times V\to\Bb C$. 
Let $(v_1,\dots,v_n)$ be a basis of an isotropic $n$-subspace of $V$. 
Let  $V_1\subset V_2\subset\dots\subset V_n$ be a flag of isotropic subspaces
where $V_i$ is 
spanned by the first $i$ vectors in the sequence  $(v_1,\dots,v_n)$. 
Then $\sigma(i_1,\dots,i_k)$ is the class in $A^{|I|}(\Cal G)$ of the cycle 
of all isotropic $n$-subspaces $L$ in $V$ such that 
$$
\dim (L\cap V_{n+1-i_h})\ge h, \ h=1,\dots,k.
$$

\vs\vs \goodbreak
{\bf A connection with Q-polynomials.}  
In [P4] the author has deduced from Theorem~6.4 the following 
result, where $Q_I$  denotes the Schur $Q$-polynomial (see Section~1):

\proclaim{Theorem 6.10 {\rm[P4, Sect.~6]}}The assignment 
$Q_I\mapsto\sigma(I)$ for $I\subset\rho_n$ --- zero otherwise, defines 
a \underbar{ring} homomorphism and allows one to identify $A^*(\Cal G)$  with 
the factor of the ring of $Q$-polynomials modulo the ideal $\bigoplus\Bb ZQ_I$, 
the sum over all $I$ not contained in $\rho_n$.
\endproclaim

In particular the following Giambelli-type formula is valid for the 
Schubert cycle  $\sigma(I)$. 

\proclaim{Theorem 6.11 \rm [P4, Sect. 6]}Let 
$I=(i_1,\dots,i_k)\subset\rho_n$ be a strict partition, 
$k$-even $($we can always assume it by putting $i_k=0$, if necessary\/$)$. 
Then
$$
\sigma(I)=\operatorname{Pfaffian} \bigl[\sigma(i_p,i_q)\bigr]_
{1\le p,q\le k},
$$
where, for $p<q$, $\sigma(i_p,i_q)=\sigma(i_p)\sigma(i_q) + 2\sum\limits^q_{h=1}\ 
(-1)^h\sigma(i_p+h)\sigma(i_q-h)$, and $\sigma(i_p,0)=\sigma(i_p)$.
\endproclaim

With the use of the tautological Lagrangian (sub)bundle $R$ on 
$\Cal G$, this result is rewritten as $\sigma(I)=
\widetilde Q_I R\hak$ (see Section 4 for the definition
of the right-hand side).

Schubert calculus for usual Grassmannians is based on three main 
theorems: Pieri's formula, Giambelli determinantal formula and the basis 
theorem. 
Analogues of the formulas are given in Theorems 6.4 and 6.11. 
A basis-type theorem for the Lagrangian Grassmannian $\Cal G$, 
can be formulated as 

\proclaim{Theorem 6.12}$A_*(\Cal G) = \bigoplus \Bb Z \sigma (I)$, \ 
where the sum is over all strict partitions $I\subset\rho_n$. 
\endproclaim

This result globalizes to Lagrangian Grassmannian bundles
$LG_n(V)\to X$ (the notation as in Section 4).
One has $A^*(LG_n(V)) = \bigoplus A^*(X) \widetilde Q_I R\hak$,
the sum over all strict partitions $I\subset \rho_n$.

Theorem 6.12 can be deduced from a general theory of the cellular 
Schubert-Bruhat decompositions of homogeneous spaces  (see, e.g.,
[Ch], [B-G-G], [D2]). 
The cellular decomposition in the case of $\Cal G$ was described in details in 
[P4, Sect.~6]. 
Another simple, conceptual proof of Theorem 6.12 is given in [P-R2]. 
\vs

Using similar arguments one can prove a Pieri-type formula for the 
Grassmannian of $n$-dimensional isotropic subspaces of $(2n+1)$-dimensional
vector space endowed with a symmetric nondegenerate form (a Pieri-type 
formula in this case was originally given in [H-B]). 
A connection with {\it $P$-polynomials} and a Giambelli-type formula were
originally established in [P4, Sect.~6].

Analogous results in the case of Grassmannian of $n$-dimensional isotropic 
subspaces in a $2n$-dimensional vector space endowed with a symmetric 
nondegenerate form  were worked out in [P4, Sect.~6]. 
\vs

The Giambelli-type formulas for isotropic Grassmannians described above
are globalized in [P-R5] to Lagrangian and orthogonal degeneracy loci.
In particular, the following problem is solved:

\proclaim{Problem \rm(J. Harris)}\rm Let $V$ be a vector bundle over $X$
equipped with a nondegenerate antisymmetric or symmetric form and let
$E$, $F$ be two maximal isotropic subbundles of $V$.
Express the fundamental class of the locus
$$
D^k=\bigl\{x\in X\mid\dim (E(x)\cap F(x))\ge k\bigr\}
$$
as a polynomial in the Chern classes of the bundles involved. 
(For a definition
of an appropriate scheme structure on $D^k$ --- see [P-R5].)
\endproclaim

Observe that if $X=LG_n(\Bb C\,\strut^{2n})$, $E$ is the tautological vector 
bundle on $X$ and $F$ is a trivial isotropic rank $n$ bundle on $X$, then the 
fundamental class of the above locus is the Schubert cycle
$\sigma(k,k-1,\dots,1)$.
By the above, this Schubert cycle is expressed by 
$\widetilde Q_{(k,k-1,\dots,1)}E\hak$. 

More generally, one has:

\proclaim{Theorem 6.13 \rm [P-R5]}Let $V$ be a vector bundle over a 
pure-dimensional Cohen-Macaulay variety 
endowed with a nondegenerate antisymmetric form.
Let $E,\,F\subset V$ be two maximal isotropic subbundles.
If $D^k$ is of pure codimension $k(k+1)/2$ in $X$ or empty, then
$$
[D^{k}]=\sum_{I\subset \rho_k} \widetilde Q_IE\hak\cdot
\widetilde Q_{\rho_k\smallsetminus I}F\hak,
$$
where $\rho_k=(k,k-1,\dots,1)$ and for a strict $I\subset\rho_k$, 
$\rho_k\smallsetminus I$ is the strict partition whose parts complement $I$ in
$\{k,k-1,\dots,1\}$\footnote{Note added in proof: 
A different (in its form) solution
to Harris' problem has been obtained independently by W. Fulton
in {\it Determinantal Formulas for Orthogonal and Symplectic Degeneracy Loci} --- to
appear in J. Differential Geom. and {\it Schubert Varieties in Flag Bundles for the
Classical Groups} --- to appear in ``Hirzebruch 65", Israel Math. Conf. Proc.
Both the papers depend on the paper by D. Edidin and W. Graham, {\it Characteristic
classes and quadric bundles} --- Duke Math. J. 78 (1995), 277--299.}.
\endproclaim

There exists ([P-R5]) an analogue of this theorem for the symmetric form;
this analogue has allowed recently a computation of the cohomology class of
Brill-Noether loci in Prym varieties.
Let $K$ be an algebraically closed field of characteristic different from $2$.
Let $\pi:\widetilde C \to C$ be an \'etale double cover of a nonsingular algebraic
curve $C$ over $K$ of genus $g$ (hence the genus of $\widetilde C$ is $2g-1$).
Let $\Nm: \Pic ^{2g-2}(\widetilde C)\to \Pic ^{2g-2}(C)$ be the norm map associated
with $\pi$; the scheme $\Nm ^{-1}(\omega_C)$, $\omega_C$ being the canonical
class, breaks up into two connected components $P^+$ and $P^-$ corresponding
to even or odd values of $h^0(-)$. These components, being translates of
the Prym variety associated with $\pi$, are irreducible of dimension $g-1$.
The following definition of the closed subsets $V^r$ in $P^{\pm}$ (called
``Brill-Noether loci for Pryms") is due to G. Welters.
For every integer $r\ge -1$ one sets:
$$
V^r = \{ L\in \Nm ^{-1}(\omega_C) | \ h^0(L)\ge r+1 \ 
\hbox{and} \ h^0(L)\equiv r+1 \pmod 2 \}. $$

\proclaim{Theorem 6.14 \rm [DC-P]}Assume 
that $V^r$ either is empty or has pure codimension in $P^{\pm}$
equal to $r(r+1)/2$. Then its class in the numerical equivalence ring
of $P^{\pm}$, or its cohomology class in $H^*(P^{\pm},\Bb C)$ \ 
(for $K=\Bb C$), is equal to
$$2^{r(r-1)/2}\prod\limits_{i=1}^{r} {(i-1)! \over (2i-1)!}\xi ^{r(r+1)/2}.$$
where $\xi$ is the class of the theta divisor on $P^{\pm}$.
\endproclaim

The assumptions of the theorem, as shown by Welters, are satisfied
for a general curve $C$ and any irreducible double cover $\pi:
\widetilde C \to C$ of it. 
For more details, concerning, e.g., a definition of an appropriate
scheme structure on $V^r$, see [DC-P] and the references
therein.
          
The paper [P-R5] gives also formulas for more general loci; in the following
$a.$ is a sequence 
$a.=(1\le a_1<\dots<a_k\le n)$ :
$$
D(a.)=\bigl\{x\in X\mid\dim (E(x)\cap F_{a_p}(x))\ge p,\ p=1,\dots,k\bigr\}
$$
where $F_{\bullet}=F_1\subset F_2\subset\dots\subset F_n=F$ is an isotropic flag.
The fundamental classes of these loci are given in loc.cit.
as quadratic expressions in $\widetilde Q$-polynomials of $E$ and flag
Schur polynomials of $F_{\bullet}$ \ (i.e. the determinants of the 
matrices of the form $[s_{i_p-p+q}(F_p)]_{p,q}$).

By adapting the technique of Hodge and Pedoe, S. Sert\"oz has obtained
in [Se] a ``triple Pieri intersection theorem" for isotropic 
Grassmannians in the orthogonal case; however, the main result in [P-R 3,4] 
cannot be deduced from [Se].

\head 7. Numerically  positive  polynomials  for  ample  vector  bundles
with  applications  to  Schur polynomials  of  Schur  bundles
and  a  vanishing  theorem.\endhead
\vs
\vbox{
\line {\hfill {\sevenrm Question: ``Which result of algebraic geometry 
is the most useful for combinatorists?"}} 
\line {\hfill {\sevenrm Answer: ``Perhaps \dots Hard Lefschetz!"}}}
\vs

In this section, we denote by  $X$  a nonsingular 
projective variety. Let $L=\Cal O_X(1)$. 
Suppose that $c_1,\dots,c_e$ are independent variables such that $\deg c_i=i$. 
We say, following [F-L], that a polynomial 
$P\in \Bb Z[c_1,\dots,c_e]$, $\deg P=d$, 
is {\it numerically positive for ample vector bundles} if for every $X$ 
of dimension $n\ge d$ and every ample vector 
bundle of rank $e$ on $X$, the number  
$P(c_1(E),\dots,c_e(E))\cdot c_1(L)^{n-d}\in A_0(X)=\Bb Z$, is positive.

Recall (see Section 1) that
the Schur polynomial associated with a partition $I$ is defined by 
$$
s_I(c.)=\Det \left[c_{j_p-p+q}\right]_{1\le p,q\le i_1},
$$
where $J=I^{\sim}$.
One has

\proclaim{Theorem 7.1 \rm[F-L]}Every numerically positive polynomial is 
a $\Bb Z$-combination\break               %%JKK
$\sum d_I\ s_I(c.)$, where all $d_I\ge 0$ and 
$\sum d_I>0$%
\footnote{Note added in proof: W. Fulton has obtained
recently an analogue of this result for filtered ample bundles where the role
of Schur polynomials is played by Schubert polynomials, recalled in
Sections 3 and 4. See: W. Fulton, {\it Positive polynomials for ample vector
bundles}, Amer. J. Math. 117 (1995), 627--633.}.
\endproclaim

This theorem,
based on the Hard Lefschetz theorem and the Giambelli-Kempf-Laksov formula
(i.e. the formula from Example 3.4),
has a nice algebraic consequence. We follow the notation
from Section 2.

The following observation is due to the author.

\proclaim{Corollary 7.2}Let $E_1,\dots,E_k$ be vector bundles. 
Then for partitions $I;$ $J_1,\dots,J_k;$ $K_1,\dots,K_k$, the coefficients  
$d(I;J_1,\dots,J_k;K_1,\dots,K_k)$ appearing in the decomposition
$$
s_I\bigl(S^{J_1}E_1\otimes\dots\otimes S^{J_k}E_k\bigr) = 
\sum\limits_{K.}\ d(I;J_1,\dots,J_k;K_1,\dots,K_k)\ s_{K_1}(E_1)\cdot\dots\cdot 
s_{K_k}(E_k)
$$
satisfy  $d(I;J_1,\dots,J_k;K_1,\dots,K_k)\ge 0$ and 
$\sum\limits_{K.}\ d(I;J_1,\dots,J_k;K_1,\dots,K_k) > 0$.
\endproclaim

Of course, the coefficients $d(I;J_1,\dots,J_k;K_1,\dots,K_k)$ 
are universal, i.e., they depend 
on $\rank (E_1)$,\dots, $\rank (E_k)$  only. 
Therefore, it is sufficient to show the assertion for ample bundles 
$E_1,\dots,E_k$ 
with algebraically independent Chern roots. 
This can be achieved by taking as $X$ the product of Grassmannians
$$
X = \prod\limits^k_{i=1}\ G^{e_i}(\Bb C^{n_i}),
$$
where $e_i= \rank  E_i$ and $n_i\gg 0$.  
Let $L=\Cal O_X(1)$. 
Moreover, denote by $p_i$ the
projection $X\to G^{e_i}(\Bb C^{n_i})$  and by $Q_i$ the tautological 
quotient bundle
on $G^{e_i}(\Bb C^{n_i})$. 
Define $E_i:=(p^*_iQ_i)\otimes L$ \ $(i=1,\dots,k)$. 
The bundles $E_i$  are ample.
Hence $S^{J_1}E_1\otimes\dots\otimes S^{J_k}E_k$ is ample (see [Ha]). 
The assertion now follows from Theorem 7.1.
\vs

For example, the numbers $D^{m,n}_{I,J}$, $(\!(J)\!)$ and $[J]$ from Section 2 
are nonnegative.

It would be interesting to have a purely algebro-combinatorial 
proof of Corollary 7.2. 
Note that a combinatorial interpretation of the numbers $D^{m,n}_{I,J}$  was
given in [G-V].

In [D-P-S], the authors show that the above {\it Fulton-Lazarsfeld 
inequalities}
for Chern classes and Schur polynomials 
of ample vector bundles still hold for {\it nef} vector
bundles on compact K\"ahler manifolds.  Recall that 
a vector bundle $E$ is {\it numerically effective} if the Grothendieck
line bundle $\Cal O(1)$ on $G^1(E)$ is nef, i.e. has a nonnegative 
degree on each effective curve in $G^1(E)$.

\ssk

As shown by L. Manivel in the next theorem, the numerical positivity
of Schur polynomials of vector bundles can be used to extend
vanishing theorems to wider classes of vector bundles.

\proclaim{Theorem 7.3 \rm[Ma]}Let $E$ be a numerically effective
vector bundle on a projective variety $X$ of dimension $n$.
Suppose that a line bundle $L$ on $X$ is nef and either there exists 
a partition $I$ such that $\int_X s_I(E)>0$ or $L$ is ``big"
$($or equivalently $\int_X c_1(L)^n>0)$. Then, for every partition $J$
of length $l(J)\ge l(I)$,
$$
H^q(X,K_X\otimes S^JE\otimes (\Det E)^l\otimes L)=0$$
for $q>0$ and $l\ge l(J)$. 
\endproclaim    

This result, extending the familiar Kawamata-Viehweg vanishing theorem,
is derived from an expression of the self-intersection
of a line bundle on a relative flag manifold. This expression is
a consequence of some expansion of the Chern character of symmetric
powers combined with the Hirzebruch-Riemann-Roch formula [H], 
and gives an insight into the corresponding Gysin map that seems to be
not reachable by the formulas of Sections 4 and~1.

\vskip12pt plus2pt minus2pt
\centerline{\bf APPENDICES}

\head A.1. Proof of Proposition 1.3(ii).\endhead
Let  $\Fl ^{q_1,\ldots,q_l}(E)\to X$ be the flag bundle
para\-metrizing flags of successive quotients of ranks $q_1,\ldots,q_l$ of $E$.
It is endowed with a sequence of tautological quotients $Q^i$ where
$\rank Q^i=i$.
If $(q_1,\ldots,q_l)=(l,l-1,\ldots,1)$ then the projection in the flag
bundle is denoted by $\tau_E^l$.
The proof goes as follows.
First, one shows the assertion for $J=\emptyset$ and any $I$.
We use a commutative diagram
\vs
$$
\CD
\Fl ^{k,k-1,\dots,1}(Q)@= \Fl ^{q,k,k-1,\dots,1}(E) @=  G^{q-k}(K)\\
@VV\tau^k_Q V&& @VV\pi'V \\
G^q(E) @>\pi_E>> X @<\tau^k_E<< \Fl ^{k,k-1,\dots,1}(E)
\endCD
$$
\vs
\noindent
where $K=\Ker (E\to Q^k)$. Invoking [P3, Corollary 2.7], we have:

\vs 

$$
\eqalign{
(\pi_E)_*&\bigl(c_\top (R_E\otimes Q_E)\cdot 
P_I(Q_E)\cap\pi_E^*\alpha\bigr)= \cr
&=(\pi_E \tau^k_Q)_*\biggl(a_1^{i_1}\dots a_k^{i_k}
\prod\limits_{\textstyle{1\le i<j\le q\atop i\le k}}(a_i+a_j)
\prod\limits_{\textstyle{1\le i\le q\atop q<j\le n}}(a_i+a_j)\cap
(\pi_E\tau_Q^k)^*\alpha\biggr)  \cr
&=(\tau_E^k)_*\pi_*'\biggl(c_\top (R_K\otimes Q_K)\cdot a_1^{i_1}\dots a_k^{i_k}
\prod\limits_{\textstyle{1\le i<j\le n\atop i\le k}}
(a_i+a_j)\cap(\tau_E^k\pi')^*\alpha\biggr)  \cr
&=\pi_*'c_\top (R_K\otimes Q_K)\cdot P_I(E)\cap\alpha.
}$$
\vs\vs
\noindent
Now, it follows easily from Proposition 1.3(i) that
$\pi_*'c_\top (R_K\otimes Q_K)$ equals
$$\card \bigl\{I\subset(q-k)^{n-q} :
|I|\hbox{ even}\bigr\} - \card \bigl\{I\subset(q-k)^{n-q} : 
|I|\hbox{ odd}\bigr\}$$
which is the requested multiplicity $d$ in this case.

\vs\vs
Passing to the dual Grassmannian, we prove the formula 
for $I=\emptyset$ and any $J$.
If $|I|\ge 1$ (and $J$ is arbitrary) we proceed as in [P4, pp.~154--155] with
the following changes:
\vs\vs
\noindent
p.~155 l. $-6$\quad should read: \quad `` \ $\ldots = dP_{I',J}(R')
\cap\pi_2^*\alpha$ \ " ,

\vs
\noindent
p.~155 l. $-4$\quad should read: \quad
`` \ $\ldots = \pi_{2*}\bigl[d\cdot e_2P_{I',J}(R')
\cdot\xi^{i_i}\cap\pi_2^*\alpha
\bigr]=dP_{I,J}(E)\cap\alpha$". 
\vs\vs
\noindent
(Observe that the multiplicity ``$d$" associated with $P_{I',J}(R')$ is equal
to the one associated with $P_{I,J}(E)$.) 

\vs\vs

\head A.2. Proof of Proposition 2.1.\endhead
We need the following lemma.

\proclaim{Lemma \rm [L2]}Let $A=(a_1,\ldots,a_n)$ be a sequence 
of independent variables. Then the following equality holds
$$
s_I(a_1+1,\ldots,a_n+1) = \sum_{J\subset I} d_{IJ} \ s_J(A),
$$
where $I$ and $J$ denote partitions and 
$$
d_{IJ} = \Det\left[{i_p+n-p \choose j_q+n-q}\right]_{1\le p,q\le n}.
$$
\rm (As Lascoux points out, the simplest proof of this identity uses
the familiar Binet-Cauchy formula --- see, e.g., [L-S1]. 
These binomial determinants appear as counting certain correlations in 
[G2].)
\endproclaim

We pass to the proof of Proposition 2.1.
Let us decompose formally $E=L_1\oplus\ldots\oplus L_n$ 
where $\rank L_i=1$. Then, according to Schur's Thesis,
if $S^IE=\bigoplus L_1^{t_1}\otimes\dots\otimes 
L_n^{t_n}$, then $s_I(E)=\sum a_1^{t_1}\dots a_n^{t_n}$,
where $a_i=c_1(L_i)$ and both the sums are taken over the same multiset of 
sequences $(t_1,\dots,t_n)$.
More precisely, 
$s_J(A)=
\sum a^T$ where the sum runs over all standard tableaux $T$ of shape
$J$ filled up with $1,\ldots,n$ and $a^T=a_1^{t_1}\cdot \ldots 
\cdot a_n^{t_n}$, $t_i$ being
the cardinality of boxes with ``$i$" in $T$ (see, e.g. [M1]). 
Hence we want to compute the decomposition into Schur polynomials of
$$
\prod\limits_T (1+\sum_i t_i a_i),
$$
where the product is taken over all standard tableaux $T$ of shape $J$
filled up with $1,\ldots,n$.  
Note that by the assumption we have the following equation of symmetric
polynomials in independent variables $x_1,\ldots,x_n$: 
$$
\prod\limits_T (\sum_i t_i x_i) = \sum_K m_K \ s_K(x_1,\ldots,x_n).
$$
Using this equation and $\sum t_i = |J|$, we get: 
$$\eqalign{
\prod\limits_T (1+ \sum_i t_ia_i) &= |J|^{-\rank (S^JE)} \prod\limits_T 
(|J|+\sum_i |J|t_ia_i)\cr
&= |J|^{-\rank (S^JE)} \prod\limits_T \ (\sum_i t_i(1+|J|a_i))\cr
&= |J|^{-\rank (S^JE)} \sum_K m_K \ s_K(\ldots,1+|J|a_i,\ldots)\cr
&= |J|^{-\rank (S^JE)} \sum_K \sum_{L\subset K} |J|^{|L|} \ m_K \ d_{KL} \ 
s_L(E)\cr}
$$
by the lemma  and the homogeneity of Schur polynomials.

\head A.3. Recursive linear relations for (\!(J)\!) and [J].\endhead
We give here an alternative proof, due to 
A. Lascoux, of the enumeration of
complete quadrics by Schubert [S] and Giambelli [G2] given in [La-La-T].

Let $A=(a_1,a_2,\ldots)$ be a sequence of independent degree 1 indeterminates
and let $A_r=(a_1,a_2,\ldots,a_r)$.
Define two power series in $A_r$:
$$\aligned
F_r & :=\prod\limits_{1\le i<j\le r}\ {1\over 1-(a_i+a_j)}  \\
\hbox{and\ \ } H_r & :=\prod\limits_{1\le i\le j\le r}\ {1\over 1-(a_i+a_j)}.
\endaligned$$

Let $\tr:\Bb Z[A_{r+1}]\to\Bb Z[A_{r+1}]$ denote, in this section,
the operator $(-1)^r\partial_1\circ\ldots\circ\partial_r$. 
Let us make the following change of variables: $b_1=-a_{r+1}$,
$b_2=-a_r,\ \ldots\ b_{r+1}=-a_1$.
Note that with respect to the $b$-variables the above operator becomes
$\partial_{\tau}\circ\ldots\circ\partial_1$ (i.e. is the operator $\tr$
studied in Section 4).
Denoting $B^r=(b_2,\ldots,b_{r+1})$, we have 
$\tr\!\!\bigl(b_1^i\cdot s_J(B^r)\bigr)=s_{(i-r)J}(B_{r+1})$ by Proposition 1.3(i) 
where $q=1$.

\proclaim{Lemma 1}The following equalities hold:
{\parindent=25pt
\item{\rm 1)} $\tr\!\!(F_r)=F_{r+1}$ \ if $r$ is even, \  and $0$ --- otherwise.
\item{\rm 2)} $\tr\!\!\bigl((a_1+\ldots+a_r)\cdot F_r\bigr)=
\bigl({r+1\over 2}-(a_1+\ldots+a_{r+1})\bigr)\cdot F_{r+1}$ \ if $r$ is odd,
 \ and ${r\over 2}\cdot F_{r+1}$ if $r$ is even.
\item{\rm 3)} $\tr\!\!(H_r)=\bigl((r+1)-2(a_1+\ldots+a_{r+1})\bigr)\cdot H_{r+1}$.
\par}
\endproclaim

\Proof
 1) We have 
$$\displaylines{
\qquad F_r=F_{r+1}(1+b_1+b_2)\ldots(1+b_1+b_{r+1})\hfill\cr
\hfill =F_{r+1}\Bigl((1+b_1)^r+(1+b_1)^{r-1}s_1(B^r)+\ldots
+(1+b_1)^0s_{(1)^r}(B^r)\Bigr)\qquad \cr}
$$
It follows from the remark preceding the lemma that
$$
\tr\!\!\bigl((1+b_1)^{r-j}s_{(1)^j}(B^r)\bigr)=\tr\!\!\bigl(b_1^{r-j}s_{(1)^j}
(B^r)\bigr)=(-1)^j
$$
for $0\le j\le r$.
Therefore, $\tr\!\!(F_r)=(1-1+1-1+\ldots)F_{r+1}$, where the first factor 
contains  $r+1$ times \ ``$\pm 1$".
This yields the assertion.
\ssk

\noindent
2) \ By similar computations we get, for $0<j\le r$, the following
equality (for the rest of the proof we use the variables $a_i$,
set $a:=a_{r+1}$ and adapt a standard $\lambda$-ring notation):
$$
\tr\!\!\bigl(a(1-a)^{r-j}s_j(-A_r)\bigr)=
\tr\!\!\bigl((r-j)a(-a)^{r-j-1}s_j(-A_r)\bigr)=(-1)^{j+1}(r-j)
$$
Also, 
$$
\tr\!\!\bigl(a(1-a)^r\bigr)=a_1+\ldots+a_{r+1}-r.
$$
Taking these equalities into account, we compute
$$\aligned
\tr\!\!\bigl((a_1 &+\ldots +a_r)\cdot F_r\bigr)=  \\
&=(a_1+\ldots+a_{r+1})\tr\!\!(F_r)-\tr\!\!(a\cdot F_r) \\
&=(a_1+\ldots+a_{r+1})\tr\!\!(F_r) \\
&\qquad
-\tr\!\!\Bigl(a(1-a)^r+a(1-a)^{r-1}s_1(-A_r)+\ldots+a(1-a)^0s_r(-A_r)\Bigr)
F_{r+1}  \\
&=(a_1+\ldots+a_{r+1})\tr\!\!(F_r) \\
&\qquad\qquad
+\Bigl(r-(r-1)+(r-2)-\ldots\pm 1-(a_1+\ldots+a_{r+1})\Bigr)F_{r+1}  \\
& = \left\{{r\over 2}F_{r+1} \hfill\hbox{if $r$ is even,}
\atop \bigl({r+1\over 2}-(a_1+\ldots+a_{r+1})\bigr)F_{r+1}\quad
\hbox{if $r$ is odd.}\right.
\endaligned
$$
\ssk

\noindent
3) \ We have 
$$\aligned
H_r & =H_{r+1}\cdot(1-a-a_1)\ldots(1-a-a_{r+1})  \\
& =H_{r+1}\cdot\bigl((1-a)^{r+1}+(1-a)^rs_1(-A_{r+1})+\ldots+(1-a)^0
s_{r+1}(-A_{r+1})\bigr).
\endaligned
$$
By applying $\tr$ to $(1-a)^{r+1-j}s_j(-A_{r+1})$, $0\le j\le r+1$, one gets
a nonzero result only for $j=0,1:$
$$
\tr\!\!\bigl((1-a)^{r+1}\bigr)=-(a_1+\ldots+a_{r+1})+(r+1)
$$
and
$$
\tr\!\!\bigl((1-a)^rs_1(-A_{r+1})\bigr)=-(a_1+\ldots+a_{r+1}).
$$
Hence $\tr\!\!(H_r)=\bigl((r+1)-2(a_1+\ldots+a_{r+1})\bigr)\cdot H_{r+1}.$
\endproof

It will be now convenient to use the following notation:
for a partition $I=(i_1,\ldots,i_r)\,$, \ $s_I(A_r)=:s(\rho_{r-1}+I;A_r)$.\ 
We define, for the use of this appendix, the numbers $(\!(J)\!)$ and $[J]$
by
$$\aligned
F_r & =\sum \ [J] \ s(J;A_r)  \\
H_r & =\sum \ (\!(J)\!) \ s(J;A_r)
\endaligned$$
(thus $J=(j_1>j_2>\ldots>j_{r-1}>j_r\ge 0)$\ ).

\proclaim{Proposition 2}\par \parindent=\wd1
\j i;The numbers $[J]$ satisfy the following
recursive linear relations:
\itemitem{\rm 1)} $[1,0]=1$
\itemitem{\rm 2)}
$p[j_1,\dots,j_{2p}]-\sum\limits_k[j_1,\ldots,j_k-1,\ldots,j_{2p}]=$ 
$$=\cases
0 & \hbox{if \ }(j_{2p-1},j_{2p})\ne(1,0),\\
[j_1,\ldots,j_{2p-2}]& \hbox{if \ }(j_{2p-1},j_{2p})=(1,0).\\ \endcases
$$

\itemitem{ } We assume the terms with $j_{k+1}=j_k-1$ in the above sum 
to be zero.
If $r$ is odd then
$$
[j_1,\ldots,j_r]=\cases
[j_1,\ldots,j_{r-1}]& \hbox{if \ } j_r=0\\
0&\hbox{if \ } j_r>0.\\ \endcases
$$
\j ii;The numbers $(\!(J)\!)$ satisfy the following recursive linear relations:
\itemitem{\rm 1)} $(\!(1,0)\!)=1$
\itemitem{\rm 2)} $r(\!(j_1,j_2.\ldots,j_r)\!)-2\sum\limits_k
(\!(j_1,\ldots,j_k-1,\ldots,j_r)\!)=$
$$=\cases
0 &\hbox{if \ }j_r>0,\\
(\!(j_1-1,\ldots,j_{r-1}-1)\!)\, &\hbox{if \ }j_r=0.\\ \endcases $$

\itemitem{} We assume the terms with $j_{k+1}=j_k-1$ in the above sum to be zero.
\par
\endproclaim

\Proof
(i) Since, by the lemma, $\tr\!\!(F_{2p})=F_{2p+1}$, the coefficient of 
$s(J;A_{2p})$ in $F_{2p}$ is the same as the one of 
$s(J;A_{2p+1})$ in $F_{2p+1}$.
In particular, $[j_1,\ldots,j_{2p},j_{2p+1}]=0$ if $j_{2p+1}>0$, and
$[j_1,\ldots,j_{2p},j_{2p+1}]=[j_1,j_2,\ldots,j_{2p}]$ if $j_{2p+1}=0$.
To get 2) we invoke the Pieri formula (we need it for $k=1$):
$$
s(j_1,\ldots,j_r;A_r)\cdot s_k(A_r)=\sum s(h_1,\ldots,h_r;A_r)
$$
where the sum is over $h_1\ge j_1>h_2\ge j_2>\ldots>h_r\ge j_r$,
$\sum h_i=\sum j_i+k$ (see [M1] and [L-S1]).
Now writing,
$$
F_{2p-1}=\sum[i_1,i_2,\ldots,i_{2p-2},0]\,
s(i_1,i_2,\ldots,i_{2p-2},0;A_{2p-1}),
$$
we express $\tr\!\!\!\bigl(F_{2p-1}\cdot s_1(A_{2p-1})\bigr)$ in two different
ways: 
first, using the Pieri formula and then applying $\tr$, and secondly, 
using assertion 2) of the lemma.
The comparison of the coefficients of Schur polynomials in the two 
expressions so obtained gives
$$\displaylines{
\quad p[j_1,j_2,\ldots,j_{2p}]-\sum[j_1,\ldots,j_k-1,\ldots,j_{2p}]=\hfill \cr
\hfill =\cases
0 & \hbox{if \ }(j_{2p-1},j_{2p})\ne(1,0),\\
[j_1,j_2,\ldots,j_{2p-2}]& \hbox{if \ }(j_{2p-1},j_{2p})=(1,0).\\
\endcases \quad \cr}
$$

\noindent
(ii) Using the Pieri formula and the equality 
$\tr\!\!(H_{r-1})=\bigl(r-2\,s_1(A_r)\bigr)\,H_r$ from the lemma, one
immediately gets the assertion.
\endproof

\proclaim{Example 3}\rm
1) $2[5,3,2,1]-[5,3,2,0]-[4,3,2,1]= 2\cdot4-7-1=0$;

$2[5,4,3,1]-[5,4,3,0]-[5,4,2,1]=2\cdot6-5-7=0$;

$2[6,3,1,0]-[6,2,1,0]-[5,3,1,0]=2\cdot25-10-12=[6,3]=28$.

\ssk\noindent
2) $F_2=1+s_1+\bigl(s_2+s(2,1)\bigr)+\bigl(s_3+2\,s(3,1)\bigr)+
\bigl(s_4+3\,s(4,1)+2\,s(3,2)\bigr)+\ldots$ \ ;

\noindent
$F_4=s(3,2,1,0)+3\,s(4,2,1,0)+\bigl(6\,s(5,2,1,0)+4\,s(4,3,1,0)\bigr)+$

\noindent
$\ +\bigl(10\,s(6,2,1,0)+12\,s(5,3,1,0)+2\,s(4,3,2,0)\bigr)+$

\noindent
$\ +\bigl(15\,s(7,2,1,0)+25\,s(6,3,1,0)+13\,s(5,4,1,0)+
       7\,s(5,3,2,0)+s(4,3,2,1)\bigr)+$

\noindent
$\ +\bigl(21\,s(8,2,1,0)+44\,s(7,3,1,0)+32\,s(6,4,1,0)+
       16\,s(6,3,2,0)+10\,s(5,4,2,0)+$     
       
\noindent
\quad $+4\,s(5,3,2,1)\bigr)+$

\noindent
\ $+\bigl(28\,s(9,2,1,0)+70\,s(8,3,1,0)+87\,s(7,4,1,0)+41\,s(6,5,1,0)+
       30\,s(7,3,2,0)+$
       
\noindent
\quad $+33\,s(6,4,2,0)+5\,s(5,4,3,0)+10\,s(6,3,2,1)+7\,s(5,4,2,1)\bigr)
+ \dots $ .

\ssk
\noindent
3) One has

$(\!(i,j)\!)-(\!(i-1,j)\!)-(\!(i,j-1)\!)=0$ for $j>0$\,;

$(\!(i,0)\!)-(\!(i-1,0)\!)-2^{i-1}=0$\,; hence $(\!(i,0)\!)=2^i-1$.
\ssk

\noindent
The number $(\!(i,j)\!)\quad(i=1,2,3,\ldots,\ j=0,1,2,\ldots)$ is given
in the $i$-th row and $j$-th column of the matrix:
$$\CD
1   &\quad & 0    &\quad & 0    &\quad & \ldots  \\
3   && 3          && 0          && \ldots  \\
7   && 10         && 10         && \ldots  \\
15  && 25         && 35         && \ldots  \\
31  && 56         && 91         && \ldots  \\
63  && 119        && 210        && \ldots  \\
\vdots && \vdots  && \vdots     &&  \\
\endCD$$
so, for example, $(\!(6,2)\!)=210$.\endproclaim

\head A.4. A Gysin map proof of the formula of Example 3.5.\endhead
Let $\tau:\Cal F=\Fl (B_1,\dots,B_k)\to X$ be the flag bundle parametrizing flags 
$V_1\subset V_2\subset\dots\subset V_k$ of vector bundles on $X$ such that
$\rank V_i=i$ and $V_i\subset B_i$ for $i=1,\dots,k$.
There is a tautological sequence of vector bundles $R_1\subset R_2\subset\dots
\subset R_k$ on $\Cal F$.
Let $Z\subset \Cal F$ be a subscheme defined by the vanishing of the
homomorphisms:
$$
R_1\to(A_1)_{\Cal F},\quad R_2\to(A_2)_{\Cal F},\quad\dots\quad,\ 
R_k\to (A_k)_{\Cal F},
$$
induced by the homomorphism $\varphi:B\to A$.
Of course, $\tau(Z)=\Omega$.
Note that $Z$ can be described by a smaller number of 
equations, i.e., defined by the
vanishing of the homomorphisms:
$$
R_1\to(A_1)_{\Cal F},\quad R_2/R_1\to(A_2)_{\Cal F},\quad\dots\quad,\ 
R_k/R_{k-1}\to (A_k)_{\Cal F}.
$$
Let us assume, for a moment, that $X$ is Cohen-Macaulay, 
$\codim_{\Cal F}Z=n_1+\dots+n_k$ and $\tau$ maps $Z$ birationally onto $\Omega$.
Then 
$$
[Z]=c_{n_1}\bigl((A_1)_{\Cal F}-R_1\bigr)\cdot 
c_{n-2}\bigl((A_2)_{\Cal F}-R_2/R_1\bigr)\cdot\dots\cdot
c_{n_k}\bigl((A_k)_{\Cal F}-R_k/R_{k-1}\bigr)\cap [\Cal F].
$$

We have (we omit writing pullback indices and brackets denoting
classes of vector bundles in the Grothendieck group, for brevity):

$$\displaylines{
\ c_{n_1}(A_1-R_1)\cdot c_{n_2}(A_2-R_2/R_1)\cdot\ldots\cdot 
c_{n_k}(A_k-R_k/R_{k-1})\hfill\cr
\noalign{\smallskip}
\hfill = \ \left|\ \CD
c_{n_1}(A_1-R_1) &\ & 0 &\ & 0 &\ &\ldots &\ & 0  \\
*  && c_{n_2}(A_2-R_2+R_1) &&  0  &&  \ldots  &&  0  \\
*  &&  *  && c_{n_3}(A_3-R_3+R_2)  && \ldots && 0 \\
\vdots && \vdots && \vdots && \ddots && \vdots \\
%\vbox{\offinterlineskip\vskip 8pt\hbox{.}\vskip 3pt\hbox{.}\vskip 3pt\hbox{.}} 
%&& \vbox{\offinterlineskip\vskip 8pt\hbox{.}\vskip 3pt\hbox{.}\vskip 3pt\hbox{.}}    
%&&      & \hskip 3pt\bullet &  &&  
%\vbox{\offinterlineskip\hbox{.}\vskip 3pt\hbox{.}\vskip 3pt\hbox{.}}  \\
%   &&     &&          && \bullet   &&    \\   
*  &&  *  &&      *    && \ldots && c_{n_k}(A_k-R_k+R_{k-1})
\endCD\ \right|\cr}
$$
where the places under the diagonal can be occupied by {\it arbitrary}
elements. Consider the $k\times k$ matrix

$$\left[\,\CD
c_{n_1}(A_1-R_1) &\quad & c_{n_1+1}(A_1-R_1) &\quad & c_{n_1+2}(A_1-R_1)
&\ & \ldots &\ & c_{n_1+k-1}(A_1-R_1)  \\
c_{n_2-1}(A_2-R_2) && c_{n_2}(A_2-R_2) && c_{n_2+1}(A_2-R_2)
&& \ldots && c_{n_2+k-2}(A_2-R_2)  \\
c_{n_3-2}(A_3-R_3) && c_{n_3-1}(A_3-R_3) && c_{n_3}(A_3-R_3)
&& \ldots && c_{n_3+k-3}(A_3-R_3)  \\
\vdots && \vdots && \vdots && \ddots && \vdots \\
%   &&    && \bullet &&&&  \\
%\vbox{\offinterlineskip\hbox{.}\vskip 3pt\hbox{.}\vskip 3pt\hbox{.}} 
%&& \vbox{\offinterlineskip\hbox{.}\vskip 3pt\hbox{.}\vskip 3pt\hbox{.}} 
%&& \hskip 44pt\bullet &&&& 
%\vbox{\offinterlineskip\hbox{.}\vskip 3pt\hbox{.}\vskip 3pt\hbox{.}} \\
%   &&      &&&& \bullet &&  \\
c_{n_k-(k-1)}(A_k-R_k) && c_{n_k-(k-2)}(A_k-R_k) && \ldots
&& \ldots &&c_{n_k}(A_k-R_k)
\endCD\,\right]$$
We claim that the determinant of this matrix equals the preceding
determinant.
To show it we record:

\proclaim{Lemma \rm [J-L]}
Let $A_1,\ldots,A_k$ and $B_1,\ldots,B_k$ be elements of a $\lambda$-ring
equipped with $\lambda$-operations $\lambda^i$;
assume that $\rank B_i\le i-1$.
Then the determinant of the matrix
$$
\bigl[\lambda^{i_p-p+q}(A_p)\bigr]_{1\le p,q\le k}
$$
remains unchanged if one replaces the argument $A_p$ by $(A_p+B_q)$ in the
$(p,q)$-place.
\endproclaim

(More precisely, this is the dual version of a result from [J-L], given
originally using Wronski's aleph functions; see the next Appendix
and also [L-S1, 7.5]).

\vs

Our claim now follows by adding $R_i$ to the argument of the $(i+1)$-th
column in the latter matrix.
\smallskip

In order to end the calculation we now invoke the following formula
for the Gysin map in $\Fl (B_1,\ldots,B_k)$ (see [K-L]): for arbitrary 
integers $p_1,\ldots,p_k$ and $\alpha\in A_*(X)$,
$$
\tau_*\Bigl(\Det\bigl[c_{p_i-i+j}(A_i-R_i)\bigr]_{1\le i,j\le k}\cap 
\tau^*\alpha \Bigr)
=
\Det\bigl[c_{p_i-m_i+j}(A_i-B_i)\bigr]_{1\le i,j\le k}\cap \alpha.
$$
Applying this formula, we infer
$$
\tau_*\Bigl(\Det\bigl[c_{n_i-i+j}(A_i-R_i)\bigr]_{1\le i,j\le k}\cap [\Cal F]
\Bigr)=
\Det\bigl[c_{n_i-m_i+j}(A_i-B_i)\bigr]_{1\le i,j\le k}\cap [X],
$$
as desired.
\vs

To get the formula in general, let $C=B\oplus A$. Then $B$ is embedded in $C$
via the graph of $\varphi$ and $A$ is a quotient of $C$ via the projection
onto the second summand. By a universality property, we have a section $s_{\varphi}:
X\to \overline X = \Fl ^{n_1,\dots,n_k}(C)$ such that the sequence 
$C\toto A_1 \toto \ldots \toto A_k$ is the pullback of the universal one
$C_{\overline X}\toto \overline A_1 \toto \ldots \toto \overline A_k$
on $\overline X$.
Let $\overline {\Omega}\subset \overline X$ be a subscheme defined by
the conditions $\dim \Ker (\overline B_i\to \overline A_i)\ge i$, 
$i=1,\ldots,k$, where $\overline B_i=(B_i)_{\overline X}$.
Then $\Omega = s_{\varphi}^{-1}\overline {\Omega}$.
Let $\overline {\Omega}\strut^{\circ}\subset \overline {\Omega}$ 
be an open subset defined
by the conditions $\dim \Ker (\overline B_i \to \overline A_i)=i$, 
$i=1,\ldots,k$.
Then $\overline {\Omega}\strut^{\circ}$ is nonempty provided
$m_i\ge i$, $m_i-i \le n_i$ for $i=1,\ldots,k$ and both the equalities
$\overline B_i=\overline B_{i+1}$ and $\overline A_i=\overline A_{i+1}$
do not hold simultaneously for $i=1,\ldots,k-1$. 
The latter condition can be rewritten as
$n_i-m_i\ge n_{i+1}-m_{i+1}+1$ \ for $i=1,\ldots,k-1$, or
equivalently
$$
n_1-m_1+1 \ge n_2-m_2+2 \ge \ \ldots \  \ge n_k-m_k+k > 0
$$
(we can assume, without loss of generality, the last inequality to be strict because the equation
$n_k-m_k+k=0$ corresponds to a redundant condition).
Under these conditions $\overline {\Omega}\strut^{\circ}$ is nonempty
and the morphism $\tau : \Cal F = \Fl (\overline B_1,\ldots,\overline B_k)\to X$
restricted to the subscheme $Z$ defined for the barred
data, induces an isomorphism over 
$\overline {\Omega}\strut^{\circ}$. Moreover, the sections of
$$
(\overline A_1)_{\Cal F}\otimes R_1\hak,\quad (\overline A_2)_{\Cal F}\otimes(R_2/R_1)\hak,
\quad\ldots\quad,\quad(\overline A_k)_{\Cal F}\otimes(R_k/R_{k-1})\hak
$$
defining $Z$, are independent, i.e. $\codim_{\Cal F}Z=n_1+\ldots+n_k$.
Then the above calculation establishes the formula for $[\overline {\Omega}]$.
The general assertion of Example 3.5 for a Cohen-Macaulay $X$ then follows
from the one just obtained because $\overline {\Omega}$ is Cohen-Macaulay
and $s_{\varphi}$ is a regular embedding; this implies (see [F1, Sect.~6,~7])
$[\Omega]=s_{\varphi}^*[\overline {\Omega}]$ provided $\codim_X\Omega =
\codim_{\overline X}\overline {\Omega}$.

\head A.5. An operator proof of the Jacobi-Trudi formula.\endhead
In the following proof of Proposition 4.4, we use the
notation from loc.cit.
We use the equality $\partial=\tr\!\!\circ\,\partial'$ and induction on $n$.
We have, by the induction assumption:

$$ \displaylines{
\qquad \partial'\bigl(a_1^{i_1+n-1}\ldots 
a_{n-1}^{i_{n-1}+1}a_n^{i_n}\bigr)\hfill\cr
\hfill =(-1)^{n-1} \left|\ \CD
0 &\quad& s_{i_1+1}(A_{n-1}) 
&\quad& \ldots& \quad & s_{i_1+n-1}(A_{n-1}) \\
\vdots && \vdots && \ddots && \vdots \\
%\vbox{\vskip -5pt}\\
%\vbox{\offinterlineskip\hbox{.}\vskip 3pt\hbox{.}\vskip 3pt\hbox{.}} 
%&\vbox{\offinterlineskip\hbox{.}\vskip 3pt\hbox{.}\vskip 3pt\hbox{.}} 
%&\vbox{\offinterlineskip\hbox{.}\vskip 3pt\hbox{\ \ .}\vskip 3pt
%                                                     \hbox{\ \ \ \ .}} 
%&\vbox{\offinterlineskip\hbox{.}\vskip 3pt\hbox{.}\vskip 3pt\hbox{.}} \\
%\vbox{\vskip -7pt}\\
0 && s_{i_{n-1}-(n-3)}(A_{n-1}) &&\ldots && s_{i_{n-1}+1}(A_{n-1}) \\
a_n^{i_n} && a_n^{i_n+1} && \ldots && a_n^{i_n+n-1}
\endCD\ \right|_{\ n\times n}\qquad}
$$

Consider now the following $n\times n$ matrix:
$$
\left[\ \CD
s_{i_1}(A_n)&\quad &s_{i_1+1}(A_n) &\quad& \ldots &\quad & s_{i_1+n-1}(A_n) \\
\vdots&& \vdots && \ddots && \vdots \\
%\vbox{\vskip -9pt}\\
%\vbox{\offinterlineskip\hbox{.}\vskip 3pt\hbox{.}\vskip 3pt\hbox{.}} 
%&\vbox{\offinterlineskip\hbox{.}\vskip 3pt\hbox{.}\vskip 3pt\hbox{.}} 
%&\quad\vbox{\offinterlineskip\hbox{.}\vskip 3pt\hbox{\ \ .}\vskip 3pt
%                                                     \hbox{\ \ \ \ .}} 
%&\vbox{\offinterlineskip\hbox{.}\vskip 3pt\hbox{.}\vskip 3pt\hbox{.}} \\
%\vbox{\vskip -11pt}\\
s_{i_{n-1}-(n-2)}(A_n) && s_{i_{n-1}-(n-3)}(A_n) &&\ldots 
&& s_{i_{n-1}+1}(A_n) \\
a_n^{i_n} &&a_n^{i_n+1} &&\ldots && a_n^{i_n+n-1}
\endCD\ \right]_{ .}
$$
We claim that the determinant of this matrix equals the previous determinant,
i.e. $(-1)^{n-1}\partial'(a^{I+\rho_{n-1}})$.
To show it we recall:

\proclaim{Jacobi-Lascoux Lemma \rm [J-L]}
Let $A_1,\ldots,A_k$ and $B_1,\ldots,B_k$ be elements of a $\lambda$-ring;
assume that $\rank B_i\le i-1$.
Then the determinant of the matrix 
$$\bigl[s_{i_p-p+q}(A_p)\bigr]_{1\le p,q\le k}$$ 
remains unchanged if one 
replaces the argument $A_p$ by $(A_p-B_q)$ in the $(p,q)$-place.
Here, $s_i(A)=(-1)^i \lambda^i(-A)$, where $\lambda^i$ is the $i$-th
$\lambda$-operation in the $\lambda$-ring.

\noindent
\rm (See also [L-S1, 7.5].)
\endproclaim

The above claim is proved by subtracting $a_n$ from the argument of the last 
$n-1$ columns (note that the last row then becomes:
$$
a_n^{i_n} \ 0 \ldots 0
$$
so we can fill up the first column in an arbitrary way without changing the
determinant).

Finally, since $\tr\!\!(a_n^k)=(-1)^{n-1}s_{k-(n-1)}(A_n)$, the final assertion
about $$\partial(a^{I+\rho_{n-1}})=\tr\!\!\biggl(\partial'\Bigl((-1)^{n-1}
(\hbox{last determinant})\Bigr)\biggr)$$ follows by the Laplace expansion
with respect to the last row and $\tr\!\!(f\cdot g)=f\cdot\!\!\tr\!\!(g)$ for 
$f\in S\Cal P(A_n)$.

\head A.6. A Schur complex proof of the Giambelli-Thom-Porteous
formula.\endhead 
First, we need a certain Schur complex constructed in [Ni].

The notation here is as in Sections 1 and 2 of the present paper.
It was proved in [Ni] that there exists a complex $C_r(\varphi).$ such that:
$$
C_r(\varphi)_i=\bigoplus\limits_{\textstyle{I\subset(m-r)^{n-r}\atop |I|=i}}
S^{I^{\sim}}(F)\otimes S^{(m-r)^{n-r}/I}(E) \leqno\hbox{\rm (i)}
$$
{(ii)} \ for every $\varphi:F\to E$ and $r\ge 0$, \  
$\Supp C_r(\varphi).=D_r(\varphi)$\footnote{Given a complex 
$C.$ of vector bundles on $X$, by $\Supp(C.)$ we understand the complement to
the set of points $x\in X$ for which $(C.)_x$ is an exact complex
of vector spaces.}.
\par\noindent
($C_r(\varphi).$ is the complex denoted by $T_{(m-r)^{n-r}}(\varphi)$ in
[Ni].) 

  Suppose that $X$ is a smooth scheme. We want to pass from the Chow ring
$A^*(X)$ to $\Gr K(X)$, where $\Gr K(X)$ is the graded Grothendieck ring 
associated with the topological filtration on $K(X)$ (see [F1, 15.1.5]).
By loc.cit. we have a functorial morphism of graded rings
$$
\phi:A^*(X)\to \Gr K(X),
$$
where, for a subvariety $V\subset X$, $\phi([V])=[\Cal O_V]$.
A fundamental property of this homomorphism is that $\phi_{\Bb Q}$ is an
isomorphism (loc.cit., 15.2.16).

\vs
Given a (finite) complex $C.$ of vector bundles on $X$, 
we denote by $[C.]$ the class $\sum (-1)^i[C_i]$ in $\Gr K(X)$.

In particular, we have:
$$[C_r(\varphi).]=\sum(-1)^{|I|}[S^{I^{\sim}}(F)]\cdot[S^{(m-r)^{n-r}/I}(E)],$$ 
the sum over $I\subset(m-r)^{n-r}$.

\proclaim{Claim}Let $D$ be an irreducible $($closed\/$)$ subscheme
of a smooth scheme $X$.
Let $C.$ be a finite complex of vector bundles on $X$ 
and let $P$ be a homogeneous element in $A^{*}(X)$ of degree $\codim_XD$.
If $\Supp(C.)\subset D$ and $\phi(P)=[C.]$, then $[D]=q\cdot P$ for some
$q\in \Bb Q$.
\endproclaim

Indeed, consider the following commutative diagram with the first row exact:
$$\CD
A_*(D) @>i_*>> A_*(X) @>j^*_A>> A_*(X\setminus D) @>>> 0  \\
            && @VV\phi V      @VV\phi V   \\
	    && \Gr K(X) @>j^*_K>> \Gr K(X\setminus D) 
\endCD$$
Here, $i:D\to X$ and $j:X\ba D\to X$ denote the inclusions. Since 
$\Supp(C.)\subset D$, we have $j^*_K([C.])=0$.
Then the equality $\phi(P)=[C.]$ implies $(j^*_A)_{\Bb Q}(P)=0$, \ 
$\phi_{\Bb Q}$ being an isomorphism.
Since $\deg P=\codim_XD$ and $D$ is irreducible, %we get 
$[D]=q\cdot P$
for some $q\in \Bb Q$.

\vs

To prove the formula in question, we apply the claim to
the triple $D=D_r(\varphi)$,
$C.=C_r(\varphi\hak).$ and $P=s_{(m-r)^{n-r}}(E-F)$.
By passing to a universal case, if necessary, we can assume that $X$ is
smooth and the assumptions of the claim are satisfied (see, e.g., 
[F1, Chap.~14]). We have $\phi(P)=[C.]$. 
To this end recall that if $L\to X$ is a line bundle then 
$\phi(c_1(L))=[\hbox{\bf 1}_X]-[L\hak]$.
Now, use the splitting principle and write formally
$E=\bigoplus^n_{i=1} L_i$ and 
$F=\bigoplus^m_{j=1} M_j$, where $\rank L_i=\rank M_j=1$.
Let $A=(a_1,\dots,a_n)$ with $a_i=c_1(L_i)$ and $B=(b_1,\dots,b_m)$ with
$b_j=c_1(M_j)$.
We must perform the transformation: \ $a_i\mapsto [\hbox{\bf 1}_X]
-[L_i\hak]$ and
$b_j\mapsto [\hbox{\bf 1}_X]-[M_j\hak]$ to the element $s_{(m-r)^{n-r}}(A-B)$.
We have by the addition/linearity formula, 
$$
s_{(m-r)^{n-r}}(A-B)=\sum(-1)^{|I|}s_{I^{\sim}}(B)\cdot 
s_{(m-r)^{n-r}/I}(A),
$$
the sum over $I\subset (m-r)^{n-r}$.
According to Schur's Thesis,
if $S^IE=\bigoplus L_1^{t_1}\otimes\dots\otimes 
L_n^{t_n}$, then $s_I(A)=\sum a_1^{t_1}\dots a_n^{t_n}$,
where both the sums are taken over the same multiset of sequences
$(t_1,\dots,t_n)$.
Hence, using (i) for $C_r(\varphi\hak).$, we see that the transformation: 
\ $a_i\mapsto -[L_i\hak]$ and
$b_j\mapsto -[M_j\hak]$,  sends the element $s_{(m-r)^{n-r}}(A-B)$
to $[C_r(\varphi\hak).]$.
Consequently it suffices to show,
that by the change $a_i\mapsto a_i+1$,
$b_j\mapsto b_j+1$, the element $s_{(m-r)^{n-r}}(A-B)$ remains unchanged.
This is the key point of the argument 
which follows, e.g., from the fact that the minimal component of the ideal
$\Cal T_r$ generalizing the resultant (see Theorem 1.5)
is generated by $s_{(m-r)^{n-r}}(c./c.')$ (for more about that, 
consult [P4, Sect.~5]).
Observe, moreover, that by (ii), $\Supp C_r(\varphi\hak).=D_r(\varphi)$.
Hence, by the claim, $[D_r(\varphi)]=q\,s_{(m-r)^{n-r}}(E-F)$ for
some $q\in \Bb Q$.
To prove $q=1$, consider a variety $X$ and a morphism $\varphi':F'\to E'$
where $\rank E'=n-r$, $\rank F'=m-r$ such that 
$\codim_X D_0(\varphi')=(m-r)(n-r)$, $A^{(m-r)(n-r)}(X)$ is a nonzero
free abelian group and $s_{(m-r)^{n-r}}(E-F)\ne 0$ (such an example is
easily constructable with the help of a Grassmannian).
Let $E=E'\oplus \hbox{\bf 1}_X^r$, $F=F'\oplus \hbox{\bf 1}_X^r$ and 
$\varphi=\varphi'\oplus \id$.
Then
$$\aligned
q\,s_{(m-r)^{n-r}}(E-F)&=[D_r(\varphi)]=[D_0(\varphi')]
=c_\top (E'{}\hak \otimes F')=s_{(m-r)^{n-r}}(E'-F')  \\
& =s_{(m-r)^{n-r}}(E-F).
\endaligned
$$
This implies $q=1$ and the proof is complete.

\head A.$\tau$. Corrigenda and addenda to some former author's papers.\endhead
\vs
\vbox{
\line {\hfill {\sevenrm ``Every (good) paper must contain an error."}}
\line {\hfill {\sevenrm T. Mostowski}}
}
\vs

We finish this article with a corrigenda and an addenda to the author's
former papers
[A-L-P], [DC-P], [L-P], [P-P1], [P1-4] and [P-R 2,3].
We apologize for all inconveniences which the misprints, inaccuracies and 
errors corrected below have caused\footnote{%
We do not correct,
however, the errors of English because of two reasons. At first, 
this would make this paper too long, and secondly, this will lead,
undoubtedly, to \dots some new errors.}.

\vs
\noindent
[A-L-P]: 
\underbar{Misprints} --- should be:
\ssk
\noindent
p.~511$_{11}$ --- `` Let $X$ be a nonsingular " /\!/ \ \ 
p.~514$^{15}$ --- `` $\ldots < i_k \le n$, " /\!/

\ssk
\noindent
p.~514$_{16}$ --- `` $[i_1,\ldots,i_k]$ " /\!/ \ \
p.~516$^{15}$ --- `` polynomial in $x_1,\ldots,x_k$, " /\!/

\ssk
\noindent
p.~517$_{11}$ --- `` $\pi: G/B \to G/P$ " /\!/ \ \ \ 
p.~517$_{10}$ --- `` $\pi^{-1}(Y_{\sigma})$ " /\!/

\ssk
\noindent
p.~517$_{3}$ --- ``  on $(Y_{\sigma}\cap Z) \times_{Y_{\sigma}} 
X_{\sigma \tau}$  "
\ssk
\centerline{\hbox{\vrule height .4pt depth 0pt width 3cm}}
\medskip

\noindent
[DC-P]:
\underbar{Addenda} --- Here is an elegant argument, pointed out 
to us by W. Fulton, justifying the footnote on page 688: 

If $\Cal L_0$ is a representative of the Poincar\'e
bundle such that  
$\Cal L_0\,|\,_{\Pic ^{2g-2}(\widetilde C)\times\{c\}}\in$ \break %% JKK
$\Pic ^0 \bigl(\Pic ^{2g-2}(\widetilde C)\bigr)$ 
and for any $L\in P^{\pm}$,
$\Nm \Cal L_0\,|\,_{\{L\}\times C} \cong \omega_C$, then there exists
$M\in \Pic ^0\bigl(\Pic ^{2g-2}(\widetilde C)\bigr)$ such that
$\Nm \Cal L_0\,|\,_{P^{\pm}\times C} \cong p^*M\otimes 
q^*\Omega_C\,|\,_{P^{\pm}\times C}$ where $p:\Pic ^{2g-2}\widetilde C\times C\to
\Pic ^{2g-2}\widetilde C$ is the projection.
Since $\Pic ^{2g-2}(\widetilde C)$ is an abelian variety, there exists
$L\in \Pic ^0\bigl(\Pic ^{2g-2}(\widetilde C)\bigr)$ such that $M=L^{\otimes 2}$.
Then $\Cal L=\Cal L_0\otimes(p')^*(L\hak)$, where
$p':\Pic ^{2g-2}(\widetilde C)\times\widetilde C\to \Pic ^{2g-2}(\widetilde C)$
is the projection, does the job.

\ssk
\centerline{\hbox{\vrule height .4pt depth 0pt width 3cm}}
\medskip

\noindent
[L-P]:
\underbar{Misprints} --- should be:
\ssk
\noindent
p.~209$_3$ --- `` $1\le i_1<\ldots <i_k\le n, \ k=1,\ldots$"/\!/
p.~210$^4$ --- `` $\{f(b_1),\dots,f(b_m)\}$ "/\!/
\ssk
\noindent
p.~210$_{14}$ --- ``$\partial_i=\partial_i^A$"/\!/
p.~211$^2$ --- ``$F,G \in \Bb Z[A]$"/\!/
p.~211$_8$ --- ``(ii) stems from [5]."/\!/
\ssk
\noindent
p.~212$^5$ --- `` $\langle X_{\mu}\omega,X_{\mu \omega} \rangle$ " \ \ /\!/
p.~212$^{13}$ --- \ `` \ $\sum \limits_{i=0}^{\infty}$ " \ \ /\!/
p.~212$_{15}$ --- \ `` $Q_i(A)$ " /\!/
\ssk
\noindent
p.~215$^8$ --- \ `` $\sigma_j^{(i)}$: \ "   \ \ \ \ \  /\!/
p.~215$_6$ ---  \ \ \ \ \ `` S\'er. I \ \ "
\ssk
\centerline{\hbox{\vrule height .4pt depth 0pt width 3cm}}
\medskip
%\eject
\noindent
[P-P1]\footnote{The following remark 
applies to [P-P1] and to the earlier
papers [P3] and [P4]. Given partitions $I,J$ such that $l(I)\le k$ and
$l(J^{\sim})\le i$, by \ \ $(i)^k+I,J$ \ \ we denote the partition
$(i+i_1,\ldots,i+i_k,j_1,j_2,\ldots)$, i.e., \ $\bigl((i)^k+I\bigr),J$ \ 
following the literal convention of [P-P1] and [P3,4]
(as it was explained on the example of the factorization formula
in [P3, Lemma 1.1] and [P4, Proposition 2.2], and then
used without further comment).}
\underbar{Revision} --- p.~192: The definition of $(\!(J)\!)$ and $[J]$ was 
incorrectly reproduced from [P3] and [La-La-T]. For a correct definition,
see Section 2 of the present paper.
\ssk
\noindent
\underbar{Addenda} --- p.~194: Theorem 5 holds under the assumption
that $\Cal Z$ is a {\it Whitney stratification}; 
see the paper by A. Parusi\'nski
and the author, {\it A formula for the Euler characteristic of the singular
hypersurfaces}, J. Algebraic Geom. 4 (1995), 337--351. 

\ssk
\centerline{\hbox{\vrule height .4pt depth 0pt width 3cm}}
\medskip

\noindent
[P1]:
\underbar {Misprint} 
\ssk
\noindent
p.~250$^{9}$ --- instead of `` $(-1)^kc_k(A)$ ", \ should be `` $c_k(a)$ ".

\ssk
\centerline{\hbox{\vrule height .4pt depth 0pt width 3cm}}

\noindent
[P2]: 
\underbar {Misprints} --- should be:
\ssk
\noindent
p.~217$_1$ --- `` $\mapsto (k$-th elementary \ " /\!/
p.~218$^8$ --- `` $\mapsto (i$-th elementary \ " /\!/
\ssk
\noindent
p.~218$^{11}$ ---  `` \ $,\ldots,a_{r+1}-b_{r+1}))_{0}$ " /\!/
p.~218$_{4}$ --- `` But by (1.3), $s_I(B'-A')$ " /\!/

\ssk
\noindent
p.~219$^{13}$ --- ``linear ordering" /\!/
p.~219$^{16}$ --- ``as the $\Bb Z[c.(a)]$-combination of $s_J(B)$'s" /\!/ 

\ssk
\noindent
p.~219$_{6}$ --- `` $\sum m_I(A'+C)s_I(B'-A')$ " /\!/
p.~221$_{7}$ --- `` diviseur de deux " /\!/

\ssk
\noindent
p.~221$_{4}$ --- `` Porteous "

\ssk
\noindent
\underbar {Revisions}: 
\ssk
\noindent
 --- In this note, by mimicking some (probably not worth recommending)
manner, we identify both mathematically
and linguistically a polynomial $a(x)$ with the equation $`` a(x)=0 "$ it
determines; consequently, we write exchangeably 
``polynomial $a(x)$" or ``equation $a(x)$".

\noindent
 --- In all formulas of the paper, $c_0(-)$ is to be understood 
to be equal to $1$. 

\noindent
 --- in the proof of Lemma 2.5, $c_k(a')$ (resp. $c_k(C)$) means in the 
analogy to $c_k(a)$ and $c_k(b)$ the $k$-th elementary symmetric
polynomial in $A'$ (resp. $C$). 

\ssk 
\centerline{\hbox{\vrule height .4pt depth 0pt width 3cm}}
\medskip

\noindent
[P3]: Here is a list of misprints and revisions (apart of those 
in [P4, pp.~185--186]):  
\ssk
\noindent
\underbar{Misprints} --- should be:

\ssk
\noindent
p.~423$^{20}$ --- \ `` $n\choose 2$ less " /\!/ \ 
p.~426$^{10}$ --- `` ... $=[R_F\hak \otimes R_E]+[F\hak\otimes Q_E].$"  \  /\!/

\ssk
\noindent
p.~426$^{13}$ -- `` $=\pi_*[s_I(Q_E)$ " /\!/  
p.~426$_{13,15,17}$ -- All three instances of `` $s_{_{\scs (m-r)^n}}$" 
should be 
\ssk
\noindent
``$s_{_{\scs (m)^{n-r}}}$" /\!/  
p.~431$_{5,10}$ --`` Lemma 3.6 " /\!/
pp.~431$_1$, 435$_{15}$ --`` $s_{(m-r)^{n-r}+I}(E-F)$ "/\!/

\ssk
\noindent
p.~435$_{10}$ --- \ ``$s_{(m-r)^{n-r},(2)}(E-F)$ " \ /\!/
p.~440$_8$ --"$[E\otimes Q]-[\Lambda^2 Q] = [R\otimes Q]+[S_2 Q]$" /\!/ 

\ssk
\noindent
p.~442$_8$ --- \ `` $[E]=[Im\varphi]+[C]$ " \ /\!/ 
p.~444$_6$ -- `` $A(X)$-module structure "/\!/ 

\ssk
\noindent
p.~447$^5$ --- `` \dots+ ${{j_p+j_q}\choose {j_q+1}}\bigr]$ " \  /\!/
p.~447$^7$ --- 
``$\sum\limits^a_{p=1}(-1)^{p-1}2^{j_p}$ " /\!/
p.~448$^3$ --- `` by (30)"/\!/ 

\ssk
\noindent 
p.~449$_1$ --- \ `` $(E_i-F_i)]$  " \ /\!/ \ 
p.~451$_6$ ---  \ ``$\sum\limits^i_{r=0}(-t)^rs_($  \  \ "/\!/

\ssk
\noindent
p.~452$^9$ ---  \ `` \ $(a_r-a_s)(a_r+a_s)^{-1}$ \ " \ /\!/ \  
p.~452$_5$ --- \ `` $\sum\limits^k_{p=1}(-1)^{p-1}w($ \ \ ".

\ssk
\noindent
\underbar{Revisions:}
\ssk
\noindent
p.~417$_7$ --- before ``We \dots" insert a footnote
{\sevenrm We assume here that $A,B$ is a sequence of independent variables, 
which can be then specialized in a commutative ring.}

\ssk
\noindent
pp.~422$_3$--423$^2$ --- should be: ``The sequence (8) allows us to treat 
the following classes of polynomials: symmetric in $A$, symmetric 
separately in $A^k$ and in $A_{n-k}$, and finally symmetric 
in $A_{n-k}$, as operators 
respectively on $A(X)$, $A(G^k(E))$ and finally on $A(\Fl ^k(E)).$"

\ssk
\noindent
p.~424 --- Proposition 2.8 is stated incorrectly for $l(I)=q-1$.
(For $l(I)=q$, it is correct as well as the proof given.)
The correct formulation for \ $l(I)=q-1$ \ is: \ \
``~$\pi_*[c_\top (R\otimes Q)\ P_I(Q) \cap\pi^*\alpha] = P_I(E)\cap \alpha$
if $\rank R$ is even, and $0$ if $\rank R$ is odd ".
Since in Section 7 we use precisely this formula when $\rank R$ is even,
the correction affects no other results and proofs in the paper.
More generally one has for $l(I)=k\le q$
$$
\pi_*\bigl[c_\top (R\otimes Q)\ P_I(Q)\cap\pi^*\alpha\bigr] = d\ P_I(E)\cap\alpha,
$$
where $d=0$ if $(q-k)(n-q)$ is odd and 
$[(n-k)/2]\choose [(q-k)/2]$ --- otherwise.
For details consult Proposition 1.3 (in the present paper).

\ssk
\noindent
p.~427$^{10}$ --- replace the given reference by: ``(cf. [F], Theorem 6.2(a))"
\ssk
\noindent
p.~446 and 447 --- By quoting [L-L-T], we were sure that its authors would
present a divided-differences  proof of Proposition 7.11 independent of the
formulas of Proposition 7.12. In the final version of [L-L-T], the authors,
however, give the proof (of Proposition 7.11) which makes use of the formulas
of Proposition 7.12. For an original, selfcontained, divided-differences
proof of Proposition 7.11 due to Lascoux, see Appendix A.3 (in the present
paper).
\ssk
\noindent
p.~449 --- Example (8.3) is revised in Example 3.5 and Appendix A.4
(in the present paper).    
\ssk
\centerline{\hbox{\vrule height .4pt depth 0pt width 3cm}}
\medskip

\noindent
[P4]: Here is a list of misprints and revisions:

\ssk
\noindent
\underbar{Misprints} --- should be\footnote{Some readers 
reported that ``7.$\tau$" on p.~185$_{13}$ is a misprint.
Actually, it is not.
In the old Mediterranean tradition, the letter $\tau$ means:
``to recognize one's error".}.

\ssk
\noindent
p.~133$_1$ --- ``\ \   $\le {i_b+b-1} \}$  " \ \  /\!/ \ \ \ \
p.~137$^8$ --- `` \ \   $\operatorname{sign}(w)w[$ \ \ \ \ " /\!/

\ssk
\ssk
\noindent
p.~176$^6$ ---  \ \ the sum is over: \ \  
`` \ \ \ $1\le i_1<\dots<i_k\le n \atop k=0,1,\dots,n$ \ \ \ " \ \ /\!/

\ssk
\noindent
pp.~181$_1$, 182$^3$ ---`` $P$ homogen(e)ous symmetric, " /\!/
p.~182$_{11}$ --- \ `` [G-Z Lemma 8] "  /\!/

\ssk
\noindent
p.~185$_4$ --- `` $\cap p^*_D(d_k)) =$ ".

\ssk
\noindent
\underbar{Revisions:}

\ssk
\noindent
p.~136$^7$ --- should be: ``Move all zeros to the right-hand end, keeping
them in order."

\ssk
\noindent
p.~137$^{12}$ ---  should be: `` was illuminated in [B-G-G] and [D]. "

\ssk
\noindent
p.~154 --- Theorem 3.3 (ii) and its proof are valid if $k=q$.
The general case $k\le q$ requires the following correction:
$$
(\pi_E)_*\bigl[c_\top (R_E\otimes Q_E)\ P_J(R_E)\ P_I(Q_E) \cap 
\pi^*_E\alpha \bigr] = dP_{I,J}(E)\cap \alpha ,
$$
where $d=0$ if $(q-k)(n-q-h)$ is odd, and $d=(-1)^{(q-k)r}
{[(n-k-h)/2]\choose [(q-k)/2]}$ --- 
otherwise.
For more details consult Proposition 1.3 and Appendix A.1 
(in the present paper).

\ssk
\noindent
p.~179, Remark 6.11; p.~181, Remark 6.16 --- : replace the Chern classes and the 
Schur polynomials in $R$ by the ones in $R\hak$.
\ssk
\centerline{\hbox{\vrule height .4pt depth 0pt width 3cm}}
\medskip

\noindent
[P-R2]:
\underbar{Revision}: In the formula of Proposition 3.5, the factor 
$\overline \partial_{\mu}^{D}(f_{\lambda})$ (equal to $1$) can be omitted.
Then add at the end of the proof of the proposition:
\ssk
\noindent
``We get
$$
m_{\mu} = \sum \overline \partial_{\mu}^D(f_{\lambda})\cdot \underline \partial
_{\mu}^D(e_p)
$$  
and the assertion follows by the properties of $f_{\lambda}$."
\ssk
\centerline{\hbox{\vrule height .4pt depth 0pt width 3cm}}
\medskip

\noindent
[P-R3]:
\underbar{Misprints} --- should be:
\ssk
\noindent
p.~1036$^{20}$ ---  ``  $w=(\tau,$  \ \ " /\!/ \ \ \ \ \
p.~1036$_{11}$ --- ``     $+d_r$ , \ $r=1,\dots,m-n$, " \ /\!/

\noindent
p.~1039, the bottom picture should have a dot ``$\bullet$" in the last 
row, i.e., the last two rows of this picture should look like:
\ssk
\centerline {$\square \ \square$}
\ssk
\centerline { \ \ $\bullet$}
\ssk
\noindent
p.~1040$^{13}$ --- `` (iii) One has $m_{\frak b}=2$. "
\bigskip

\noindent
\underbar{Revision}: p.~1039$^{6}$ --- better is: ``The {\it roof} of a deformed
(nonextremal) component is its row of highest boxes"

\enddocument
\end